\documentclass{amsart}
\usepackage{amssymb,amsmath,accents}
\usepackage{amsthm,amsfonts}
\usepackage{amscd,verbatim,psfrag,dsfont,enumerate,comment,version,url}

\newtheorem{lemma}{Lemma}[section]
\newtheorem{proposition}[lemma]{Proposition}
\newtheorem{theorem}[lemma]{Theorem}
\newtheorem{corollary}[lemma]{Corollary}
\newtheorem{definition}[lemma]{Definition}
\newtheorem{remark}[lemma]{Remark}

\newtheorem*{introtheorem1}{Theorem A}

\DeclareMathOperator{\supp}{supp}
 \DeclareMathOperator{\Hom}{Hom}
\DeclareMathOperator{\re}{Re}

\DeclareMathOperator{\dist}{dist} \DeclareMathOperator{\tr}{tr}
 \DeclareMathOperator{\Area}{Area}
\DeclareMathOperator{\Length}{Length}\DeclareMathOperator{\Genus}{Genus}
\DeclareMathOperator{\inj}{inj}

 \DeclareMathOperator{\Int}{Int}    \DeclareMathOperator{\im}{Im}
 \DeclareMathOperator{\order}{ord} \DeclareMathOperator{\cl}{cl} \DeclareMathOperator{\pr}{pr}
\DeclareMathOperator{\End}{End}
\DeclareMathOperator{\Rank}{Rank}

\newcommand{\R}{\mathbb{R}}

\newcommand{\two}{B}

\newcommand{\comp}{{\mspace{2mu}\scriptstyle \circ\mspace{2mu}}}
\newcommand{\la}{\langle}
\newcommand{\ra}{\rangle}

\begin{document}
\title{Target-local Gromov compactness}
\author{Joel W. Fish}
\address{Department of Mathematics, Stanford University, Stanford, CA 94305}
\email{joelfish@math.stanford.edu}
\urladdr{http://www.stanford.edu/~joelfish}
\date{\today}
\thanks{Research partially supported by NSF grants DMS-0602191 and DMS-0802927}
\begin{abstract}
We prove a version of Gromov's compactness theorem for pseudo-holomorphic curves which holds locally in the target symplectic manifold.  This result applies to sequences of curves with an unbounded number of free boundary components, and in families of degenerating target manifolds which have unbounded geometry (e.g. no uniform energy threshold). Core elements of the proof regard curves as submanifolds (rather than maps) and then adapt methods from the theory of minimal surfaces.
\end{abstract}
\maketitle
\tableofcontents
\section{Introduction}\label{sec:Introduction}

In his seminal 1985 paper \cite{Gm85}, Gromov introduced the notion of a ``pseudo-holomorphic curve'' and established the fundamental notion of compactness for families of such $J$-curves.  Since then, the majority of modern proofs of Gromov's compactness theorem (and its generalizations) have all followed the same basic recipe, namely
to study $J$-curves as a type of special harmonic map.  This essentially reduces the compactness problem to applying Deligne-Mumford compactness to the underlying Riemann surfaces and then applying bubbling analysis.
However, there are a growing number of examples in which this approach badly breaks down -- for instance, $J$-curves in a family of symplectic manifolds which lacks a uniform energy threshold, or sequences of $J$-curves with bounded area but unbounded topology.  Such a case was considered in the author's Ph.D. thesis \cite{Fj07}, in which a compactness result was proved for $J$-curves in the connected sum of two contact manifolds for which the connecting handle collapsed to a point. More generally, the author is interested in studying $J$-curves in symplectic cobordisms between non-compact and/or degenerate contact manifolds (i.e. manifolds for which the contact form vanishes along a submanifold). Additionally, the author is analyzing the behavior of contact homology under subcritical surgeries, and attempting to develop a more general  ``sideways stretching'' operation in Symplectic Field Theory.
A key difficulty which is common in each of these research directions is the lack of a uniform energy threshold.  The lack of this quantity is so fundamental that it necessitates an alternate approach to the compactness problem: namely, to regard $J$-curves as submanifolds,
and then by incorporating elements from minimal surface theory to prove a compactness result which holds locally in the target. Indeed, in this article and \cite{Fj09b}, the author takes precisely this approach; the main arguments for this target-local version of Gromov compactness are provided here, and the author develops supporting analysis for these arguments in \cite{Fj09b}.  In successive papers, the author will extend the following results to some non-compact cases, and refine the notion of Gromov compactness near nodes and critical points.

\subsection{Statement of main result}
The main result of this article is Theorem \ref{thm:Main2} from Section \ref{sec:Main}.  We state a simplified version (in fact an immediate corollary) as Theorem A below.

\begin{introtheorem1}
Let $(M,J,g)$ be a compact almost Hermitian\footnote{That is, for $(M,J,g)$ we require that $g$ is a Riemannian metric, and the almost complex structure $J$ is an isometry.} manifold with boundary.  Let $(J_k,g_k)$ be a sequence of almost Hermitian structures which converge to $(J,g)$ in $C^\infty(M)$, and let $(u_k,S_k,j_k,J_k)$ be a sequence of compact $J_k$-curves (possibly disconnected, but having no constant components) satisfying the following:
\begin{enumerate}
\item $u_k:\partial S_k \to \partial M$,
\item $\Area_{u_k^*g_k}(S_k)\leq C_{A}$,
\item $\Genus(S_k)\leq C_{G}$.
\end{enumerate}
Then there exists a subsequence (still denoted with subscripts $k$) of the $\mathbf{u}_k$, an $\epsilon>0$, and an open dense set $\mathcal{I}\subset [0,\epsilon)$ with the following significance.  For each $\delta \in \mathcal{I}$, define $\tilde{S}_k^\delta:=\{\zeta\in S_k: \dist_g\big(u_k(\zeta),\partial M)\geq \delta\}$; then the $J_k$-curves $(u_k,\tilde{S}_k^\delta,j_k,J_k)$ converge in a Gromov sense\footnote{For a precise formulation of Gromov convergence, see Definition \ref{def:GromovConvergence} below.}.
\end{introtheorem1}

Note that we have \emph{not} assumed that the $u_k(\partial S_k)$ lie in a Lagrangian submanifold, and we have \emph{not} assumed that the $S_k$ have bounded topological type.  Indeed, the number of connected components of either the $\partial S_k$ or the $S_k$ may not be bounded.  It is for this reason that the above result is significantly different from all other versions of Gromov compactness.

To see the relevance of Theorem A, we shall consider its application in a couple of examples. Observe that in the case that $M$ is closed, the above result only recovers the usual Gromov compactness theorem, however the strength the Theorem A becomes more apparent when considering target manifolds with rather arbitrary smooth boundary.  Furthermore, this latter scenario occurs quite naturally when considering closed symplectic manifolds for which the almost complex structure degenerates in a small region. We explore such a case at present.

\textbf{Example 1.} Consider a closed symplectic manifold $(M,\omega)$, fix $p\in M$, and let $\Phi:\mathcal{O}(p)\to \R^{2n}$ be Darboux coordinates around $p$.  Locally define the complex structure $\tilde{J}$ near $p$ by $\tilde{J} \partial_{x^i}=\partial_{y^i}$, and let $\tilde{M}$ be the manifold obtained by performing a $\tilde{J}$-complex blowup at $p$.  Recall that $\tilde{M}$ can be equipped with a family of closed two-forms $\omega_\epsilon$ which are symplectic for $\epsilon>0$.  Furthermore the $\omega_\epsilon$-volume of the divisor $\mathcal{D}\subset\tilde{M}$ tends to zero as $\epsilon$ does, and the $\omega_\epsilon$ converge in $C_{loc}^\infty(\tilde{M}\setminus \mathcal{D})$ to $\omega_0$ which has the property that $(\tilde{M}\setminus \mathcal{D},\omega_0)$ and $(M\setminus\{p\},\omega)$ are symplectomorphic. In other words, we have performed symplectic blowups of weight $\epsilon$ at $p$.  Lastly, equip $\tilde{M}$ with a family $J_\epsilon$ of $\omega_\epsilon$-compatible almost complex structures which also converge in $C_{loc}^\infty(\tilde{M}\setminus \mathcal{D})$.
We
now consider the following question: given a sequence $\epsilon_k\to 0$ and a sequence of pseudo-holomorphic curves $u_k:(\mathbb{S}^2,i)\to (\tilde{M},J_{\epsilon_k})$ with uniformly bounded $\omega_{\epsilon_k}$-energy, does there exist a subsequence which converges in a reasonable sense (e.g. in a Gromov sense)?

There are some obvious tricks if the $J_\epsilon$ are integrable in a neighborhood of $\mathcal{D}$ or if the $u_k(\mathbb{S}^2)$ have empty intersection with $\mathcal{D}$, however answering the more general question is non-trivial.  Indeed, one key point here is that by construction, this family of symplectic forms and almost complex structures \emph{lacks a uniform energy threshold}.  That is, as $\epsilon\to 0$, there exist symplectic spheres of arbitrarily small symplectic area.  This is a serious problem since almost all proofs of Gromov compactness rely on an energy threshold in a critical way: energy thresholds guarantee that only finitely many bubbles develop in the limit.  Indeed, a priori it might be the case that for the above example the gradient blows up at arbitrarily many points in $\mathbb{S}^2$.

Despite these difficulties, we see that the $J_\epsilon$ converge in $C_{loc}^\infty(\tilde{M}\setminus \mathcal{D})$ by construction, so it seems reasonable that the portion of the $J_{\epsilon_k}$-curves which have image in the complement of a neighborhood of $\mathcal{D}$ should converge in a reasonable sense.  Thus a natural attempt to prove some sort of compactness would be to fix a neighborhood $\mathcal{U}$ of $\mathcal{D}$, and define the curves $u_k:u_k^{-1}(\tilde{M}\setminus \mathcal{U})\to \tilde{M}$, and attempt to prove Gromov convergence for a subsequence of these domain-restricted curves.  The boon here is that $\tilde{M}\setminus \mathcal{U}$ is compact, it has a uniform energy threshold, and the $J_\epsilon$ converge in $C^\infty(\tilde{M}\setminus\mathcal{U})$.  However one now faces a new problem, namely that the surfaces $\tilde{S}_k:=u_k^{-1}(\tilde{M}\setminus \mathcal{U})$ have no a priori bound on the number of connected components, nor an a priori bound on the number of boundary components.  This is seriously problematic for standard proofs of Gromov compactness because a lack of a topology bound on the underlying Riemann surfaces precludes one from applying Deligne-Mumford compactness to the domain curves.  This is in turn problematic because it is the Deligne-Mumford compactness (together with a uniformization theorem) which yields convenient reparameterizations of the given pseudo-holomorphic curves.

It is at this point that we see the utility of Theorem A above.  Indeed, it is not difficult to choose $\mathcal{U}$ so that $\tilde{M}\setminus \mathcal{U}$ and the $u_k^{-1}(\tilde{M}\setminus\mathcal{U})\subset S_k$ have the structures of compact manifolds with smooth boundary.  Furthermore the restricted curves $(u_k,u_k^{-1}(\tilde{M}\setminus \mathcal{U}),j_k,J_{\epsilon_k})$ certainly satisfy the hypotheses of Theorem A.  We can then conclude that after passing to a subsequence, we have convergence of our pseudo-holomorphic curves ``away from $\mathcal{U}\supset \mathcal{D}$.''  Let us make this more precise.  From the above result, we can deduce the following: for each open set $\mathcal{U}\supset \mathcal{D}$, there exists an open set $\mathcal{V}$ such that $\mathcal{D}\subset \mathcal{V}\subset \mathcal{U}$, and there exists a subsequence of the above curves such that for $\tilde{S}_k:=u_k^{-1}(\tilde{M}\setminus \mathcal{V})$, the domain-restricted curves $\big(u_k, \tilde{S}_k,j_k,J_{\epsilon_k})$ converge in a Gromov sense.  We now make two important observations.  First, for these curves to Gromov-converge it must be the case that the $(\tilde{S}_k,j_k)$ converge in a Deligne-Mumford sense, which guarantees that the ``trimmed'' surfaces $\tilde{S}_k$ have fixed topological type for all sufficiently large $k$.  The second important point is that even though we have ``trimmed away'' some portion of our original curves $(u_k,S_k,j_k,J_{\epsilon_k})$ to obtain convergence, we have only trimmed away portions of the curves which have image in the ``small'' region $\mathcal{V}\subset \mathcal{U}$.  This latter point can be stated more concisely as $u_k(S_k\setminus \tilde{S}_k)\subset \mathcal{V}\subset \mathcal{U}$.

\textbf{Example 2.} The previous example was somewhat simplistic, so we now consider a larger class of similar, but much more general, examples.  Fix a symplectic manifold $(M,\omega)$, and consider a compact embedded submanifold $N\subset M$ with $\dim N < \dim M$.  Consider a sequence of almost complex structures $J_k$ which converge in $C_{loc}^\infty(M\setminus N)$ and which degenerate along $N$.  Again, we ask if energy bounds and genus bounds for closed curves are sufficient to obtain a convergent subsequence, and again Theorem A guarantees convergence away from $N$ for some subsequence.  Also note that this example is significantly less artificial than the previous one since it contains both the contact-type neck-stretching construction from Symplectic Field Theory, as well as the degenerating symplectic-connected sums setup arising in the symplectic sum formula for Gromov-Witten invariants.  Furthermore the condition that $N$ be a compact embedded sub-manifold is easily relaxed to the condition that $N$ be a compact set of zero measure, so one expects Theorem A to play a role in a wide variety of degeneration problems in symplectic geometry.

The above examples hopefully illuminate the flexibility and generality of Theorem A, so we now take a moment to point out certain things that it does \emph{not} guarantee.  Firstly note that that Theorem A makes no claims about curves with Lagrangian boundary condition, however in light of estimates proved in \cite{Fj09b}, it appears that such a generalization is quite probable.  Secondly, Theorem A does not guarantee convergence up to the boundary of $M$.  Indeed, since the $(u_k,\tilde{S}_k^\delta, j_k,J_k)$ converge for each $\delta\in \mathcal{I}$ where $\mathcal{I}\subset [0,\infty)$ is an open dense set in a neighborhood of $0$, one is tempted to consider a sequence $\{\delta_k\}_{k\in \mathbb{N}}\subset \mathcal{I}$ such that $\delta_k\to 0$ and then conclude from Theorem A that the subsequence $(u_k,\tilde{S}_k^{\delta_k}, j_k,J_k)$ converges in a Gromov-sense, however in general this is false.  Indeed, the key point is that after passing to the subsequence guaranteed by Theorem A, we find that for each fixed $\delta\in \mathcal{I}$ the topological type of the $\tilde{S}_k^\delta$ is bounded as $k$ varies over $\mathbb{N}$, but the topological type of $\tilde{S}_k^\delta$ is not necessarily bounded as $k$ varies over $\mathbb{N}$ \emph{and} $\delta$ varies over $\mathcal{I}$.

In discussing the limitations of Theorem A, we return to our previous examples from symplectic geometry, and make the important observation that Theorem A does not guarantee any sort of convergence along the region in which $J$ degenerates.  In other words, in the symplectic blow-up example we do not obtain convergence in collapsingly small neighborhoods of the symplectic divisor $\mathcal{D}$; in the neck stretching example we do not obtain convergence of multi-level buildings which fall in to the contact-type hyper-surface; in the degenerating symplectic-connected sums example we do not capture curves falling into the collapsing handle.
The reason that the above theorem makes no claims about compactness in these regions is that the behavior of curves in these regions is critically dependent on the manner in which $J$-degenerates -- something not specified in the hypotheses of Theorem A.  However, for all of those examples, and a wide variety of others, there exist diffeomorphisms of neighborhoods of the set $N$ along which $J$ degenerates to some long/wide/vast region $\mathcal{N}$ (e.g. $\mathcal{N}:=\R\times N$ in the neck-stretching case) on which $J$ is standard.  If the original curves were closed and of bounded topological type, then one can apply Theorem A on compact domains (e.g. $[a,b]\times N$ in the contact case) contained in the ``long'' region $\mathcal{N}$ provided one has a uniform area bound in this compact domain.  Indeed, such bounds occur quite often, and in such cases one can then use the above result to build-up a variety of compactness results in non-compact or degenerating target manifolds.  An example of both occurs in the author's Ph.D. thesis \cite{Fj07}, in which a Symplectic Field Theory type compactness result was proved for a sequence of finite energy $J$-curves in a degenerating connected sum of contact manifolds.

\subsection{Proof outline}
We begin by recalling a result on which our proof relies.  Indeed, if $\mathbf{u}_k=(u_k,S,j_k,J_k)$ is a sequence of compact pseudo-holomorphic curves with bounded area, fixed domain manifold $S$, varying conformal structures $j_k$, and has annular neighborhoods $\mathcal{A}_i$ of each component of $\partial S$ which have conformal modulus uniformly bounded away from zero, then there exists a subsequence which Gromov-converges after removing a neighborhood of small conformal modulus near the the boundary $\partial S$.  Indeed, such a result was proved in \cite{IsSv00} (stated there as Theorem 1), however the language in that article does not explicitly mention this conformal trimming since convergence there is understood on compact sets of the interior of $S$.
We mention this difference because this trimming is a subtle but critically important consideration for the results that follow.  The primary goal of this paper then becomes the following:  for pseudo-holomorphic curves as in the hypotheses of Theorem A, and each $\delta>0$, pass to a sub-sequence and find $\tilde{S}_k\subset S_k$ with the property that each of the $\tilde{S}_k$ are diffeomorphic to some $\tilde{S}$, and
\begin{equation*}
u_k(S_k\setminus\tilde{S}_k)\subset \{q\in M: \dist_g(q,\partial M)< \delta\},
\end{equation*}
and that each boundary component of $\tilde{S}_k$ has annular neighborhood $\mathcal{A}_{i,k}$ which has conformal modulus bounded away from zero, and
\begin{equation*}
u_k(\mathcal{A}_{i,k})\subset \{q\in M: \dist_g(q,\partial M)< \delta\}.
\end{equation*}
This essentially reduces the problem to the result proved in \cite{IsSv00}, and after the conformal trimming near the boundaries is taken, one is left with a subsequence of curves with no area loss in the deep interior $\{q\in M: \dist_g(q,\partial M)\geq \delta\}$ and which converges in a Gromov sense.  Theorem A can then be deduced by repeating the argument for a sequence $\delta_k\to 0$, and then passing to a diagonal subsequence.

Thus the primary difficulty addressed in this article is to find the desired trimmings.  To that end, we build the result up in three steps.  We begin by observing that pseudo-holomorphic curves satisfy a mean curvature equation of the form $H_\nu=\tr_S Q$ where $H_\nu$ is the mean curvature vector along the image of a $J$-curve $u:S\to M$, $Q$ is a $(1,2)$-tensor defined on $M$ which depends on $J$ and $g$, and by $\tr_S Q$ we mean the trace of $Q$ along planes tangent to the image of $u$.  We then incorporate elements of minimal surface theory as follows. The first step is to show that if a sequence of immersed $J$-curves has uniformly bounded area and uniformly $L^\infty$-bounded second fundamental forms $B_{u_k}$, then one can extract a convergent subsequence.  Of importance here is that boundedness of the topological type of the underlying Riemann surfaces is not assumed, but rather constructed for the subsequence in the proof. It is this result which allows one to obtain compactness without a priori knowledge of the domain topology.

In light of this result, we see that given a sequence of $J$-curves which satisfy the hypotheses of Theorem A, the goal becomes to pass to some subsequence, and to find some region of the form
\begin{equation}\label{eq:IntroSpecialRegion}
M^{\delta_0,\delta_1}:=\{q\in M: 0<\delta_0\leq \dist(q,\partial M) \leq \delta_1< \delta\}
\end{equation}
so that the portion of  the $J$-curves (in the subsequence) with image in $M^{\delta_0,\delta_1}$ are immersed and have $L^\infty$-bounded curvature. Since $J$-curves can of course develop unbounded curvature, (consider the formation of the standard node, or a the formation of a critical point from immersed curves)  and may not be immersed, we temporarily impose two additional hypotheses on the curves in question, namely that the number of critical points is uniformly bounded and the total curvature $\int\|B\|^2$ is uniformly bounded.  In \cite{Fj09b}, it was shown that the square-length of the second fundamental form of a $J$-curve satisfies an $\epsilon$-regularity result similar to the result shown in \cite{ChSr85} for minimal surfaces.  This guarantees that after passing to a subsequence, the curvature of the $J$-curves can only point-wise blow-up at finitely many points in the interior of $M$.  Consequently after passing to a subsequence, one finds a region of the form (\ref{eq:IntroSpecialRegion}) on which the $J$-curves are immersed with $L^\infty$-bounded curvature.

In light of this result, the goal then becomes to verify that neither the total curvature $\int\|B\|^2$ nor the number critical points can increase without bound on the deep interior of $M$.  The first step here is to employ a desingularization result which reduces the problem of arbitrarily many critical points to the problem of unbounded total curvature of immersed curves.  To exclude the possibility of unbounded total curvature we first argue that if $\zeta\in S$ and $\inj_{u^*g}(\zeta)$ is very small and $u(\zeta)$ is in the deep interior of $M$ and $\Genus(S)$ is zero then there exists a short closed loop the removal of which disconnects $S$ into two pieces which each contain a threshold amount of area. Iterating this argument shows that the curves in question cannot develop too many nodes in the deep interior of $M$ -- even in the case of non-zero genus. We conclude that after passing to a subsequence, the injectivity radius can only be arbitrarily small in a neighborhood of a finite number of points in $M$, so by restricting our attention to complementary regions, we may assume the injectivity radius is uniformly bounded away from zero.  Then by employing a covering argument, it is sufficient to show that on an intrinsic disk $\mathcal{D}_r(\zeta_0):=\{\zeta\in S: \dist_{u^*g}(\zeta_0,\zeta) < r\}$ a $J$-curve with a uniformly bounded area cannot have arbitrarily large total curvature. This is proved by recalling that $J$-curves have Gaussian curvature uniformly bounded from above, and recalling a differential equation relating the area and curvature of such intrinsic disks. In particular we show that if the total curvature on $\mathcal{D}_{r/2}(\zeta_0)$ is arbitrarily large, then so too is the area of the disk $\mathcal{D}_r(\zeta_0)$.  Since the $J$-curves in question have a priori bounded area, this is sufficient to conclude that the total curvature of the $J$-curves with image in the interior of $M$ is not arbitrarily large, and the proof of Theorem A is then immediate.

It should be noted that the techniques used to prove Theorem A are sufficiently strong to develop a more refined version of Gromov-convergence which neither relies on bubbling-analysis of harmonic maps nor relies on Deligne-Mumford compactness.  This approach will be addressed in future work, and for now we suffice to prove Theorem A as stated and outlined above.

\subsection{Acknowledgements}
The following is an extension of ideas developed in my Ph.D. thesis \cite{Fj07} at New York University, and as such I would like to thank my advisor, Helmut Hofer, for his encouragement, support, and for creating a vibrant symplectic and contact research group at the Courant Institute.  I would also like to thank my pseudo-advisor Richard Siefring, who always appeared happy to answer my endless list of questions about pseudo-holomorphic curves, and symplectic geometry.  Lastly, I would like to thank Chris Wendl for a variety of fruitful conversations, and for his detailed comments regarding many elements of this manuscript.

\section{Preliminaries}\label{sec:Preliminaries}

We begin by providing some pertinent definitions.  For instance, let $M$ be a compact real $2n$-dimensional manifold (possibly with boundary) equipped with a smooth section $J\in \Gamma\big(\End(TM)\big)$ for which $J^2=-\mathds{1}$; we call $(M,J)$ an \emph{almost complex manifold}, and $J$ the \emph{almost complex structure}.  Note that $J$ need not be integrable; that is, it need not be induced from local complex coordinates.  Indeed, this will only be true if the Nijenhuis tensor $N_J$ associated to $J$ vanishes identically, and do not make such an assumption.

If $(M,J)$ is equipped with a smooth Riemannian metric $g$ for which $J$ is an isometry (i.e. $g(x,y)=g(Jx,Jy)$ for all $x,y\in TM$), then we call $(M,J,g)$ an \emph{almost Hermitian manifold}.  Observe that any almost complex manifold can be given an almost Hermitian structure $(J,g)$ by choosing an arbitrary Reimannian metric $\tilde{g}$, and defining $g(x,y):=\frac{1}{2}\big(\tilde{g}(x,y)+ \tilde{g}(Jx,Jy)\big)$.

To an almost Hermitian manifold $(M,J,g)$ one can associate a fundamental two form (c.f. \cite{KsNk96b}) $\omega \in \Gamma\big(\Lambda^2 TM \big)$ given by $\omega(x,y):=g(Jx,y)$.  We call $\omega$ the \emph{almost symplectic form} associated to $(J,g)$, where the ``almost'' refers to the fact that in general $d\omega\neq 0$.  Indeed, $\omega$ is non-degenerate by definition, so if $\omega$ is closed then it is a symplectic form, and in such case $J$ is an $\omega$-compatible almost complex structure.  Again, we do not make this additional assumption.

We also consider pseudo-holomorphic curves, or more concisely $J$-curves, which for our purposes will be four-tuples $\mathbf{u}=(u,S,j,J)$, with entries defined as follows. Given a target manifold $M$, $J$ will be a smooth almost complex structure on $M$, $S$ will be a smooth manifold of real dimension two, $j$ will be a smooth almost complex structure on $S$, and $u:S\to M$ will be a smooth map for which $J\cdot Tu = Tu\cdot j$. Unless otherwise specified, we will allow for $S$ to be non-compact, to have smooth boundary, and to have unbounded topology (i.e. countably infinite connected components, boundary, and genus).
We will say that a $J$-curve $\mathbf{u}$ is \emph{compact} provided $S$ has the structure of a compact manifold with smooth boundary, and we will say $\mathbf{u}$ is \emph{closed} provided $S$ has the structure of a compact manifold without boundary. Note that we do \emph{not} assume $(S,j)$ is a Riemann surface; that is we do not assume $(S,j)$ has the structure of a complex (and hence analytic) manifold, but only that $S$ is a smooth manifold.  The reason for this non-standard assumption (or lack thereof) is that in what follows it will be absolutely necessary to parameterize $J$-curves as smooth maps from smooth surfaces which are not analytic.  Indeed, requiring $(S,j)$ to have the structure of a complex manifold is not only unnecessary in what follows, but unnecessarily cumbersome.

Since $S$ can be quite complicated, we will need to make the notion of "genus" precise.  We do this in definition \ref{def:genus} below, but first we introduce the notion of a compact region.

\begin{definition}[compact region]\label{def:CompactRegion}
Let $M$ be a manifold.  Suppose $\mathcal{U}\subset M$ is an open set for which its closure $\cl(\mathcal{U})$ inherits from $M$ the structure of a smooth compact manifold possibly with boundary. Then we call $\cl(\mathcal{U})$ a \emph{compact region} in $M$.
\end{definition}
\begin{definition}[genus]\label{def:genus}
Let $S$ be a connected compact two-dimensional manifold with boundary.  We define $\Genus (S)$ to be the genus of the surface obtained by capping off the boundary components of $S$ by disks.  If $S$ is disconnected but compact, then we define $\Genus(S):=\sum_{k=1}^n\Genus(S_k)$ where the $S_k$ are the connected components of $S$.  If $S$ is non compact (but with at most countably infinite connected components), we define $\Genus(S):=\lim_{k\to\infty} \Genus(S_k)$, where $S_1\subset S_2 \subset S_3\subset \cdots$ is an exhausting sequence of compact regions in $S$.
\end{definition}

This raises an important point, namely that we will often abuse notation by referring to the genus of $\mathbf{u}$ or $u$, when we actually mean the genus of $S$.  We will similarly abuse language by saying that $\mathbf{u}$ is connected or compact by which will we mean $S$ has these properties.

We now turn our attention to some less standard definitions, which have some flavor of geometric measure theory, and are necessary for later proofs.

\begin{definition}[$\mathcal{K}$-proper and $\mathcal{K}$-convergence]\label{def:KproperAndConvergence}
Consider a sequence of maps $u_k:S_k\to M$ to and from manifolds which possibly have boundary and may be non-compact. Let $\mathcal{K}\subset\Int(M)$ be a compact set in the interior of $M$.  We call this a \emph{robustly $\mathcal{K}$-proper} sequence provided there exists another compact set $\widetilde{\mathcal{K}}\subset \Int(M)$ for which $\mathcal{K}\subset\Int(\widetilde{\mathcal{K}})$ and if $u_k^{-1}(\widehat{\mathcal{K}})\setminus \partial S_k$ is compact for every compact set $\widehat{\mathcal{K}}\subset\widetilde{\mathcal{K}}$.  Similarly a single map $u:S\to M$ is robustly $\mathcal{K}$-proper provided the constant sequence $u,u,u,\ldots$ is robustly $\mathcal{K}$-proper.

Furthermore, we say the above sequence of maps \emph{robustly $\mathcal{K}$-converge in $C^\infty$} provided there exists an auxiliary manifold $\widetilde{S}$ and diffeomorphisms $\psi_k:\widetilde{S}\to \psi_k(\widetilde{S})\subset S_k$ with the property that $u_k\big(S_k\setminus \psi_k(\widetilde{S})\big)\subset M\setminus \widetilde{\mathcal{K}}$, and the "trimmed" reparameterizations $u_k\circ\psi_k:\widetilde{S}\to M$ converge in $C^\infty$.
\end{definition}

\begin{definition}[uniformly robust $\mathcal{K}$-covers]
Let $M$ be a manifold, and $\mathcal{K}\subset \Int(M)$ a compact set. Suppose $u:S\to M$ is a smooth robustly $\mathcal{K}$-proper map. Then we say $(u,S)$ is $\mathcal{K}$-covered by maps $\phi_{i}:\mathcal{D}_r\to S$ for $i=1,\ldots,n$ provided that
\begin{equation*}
u\big(S\setminus \cup_{i=1}^n\phi_i(\mathcal{D}_r)\big)  \subset M \setminus \mathcal{K};
\end{equation*}
here $\mathcal{D}_r:=\{X\in\mathbb{R}^{\dim S}: \|X\|< r\}$.

We say a sequence of robustly $\mathcal{K}$-proper maps $u_k:S_k\to M$ is \emph{uniformly} $\mathcal{K}$-covered provided $\dim S_k$ is independent of $k$ and each $(u_k,S_k)$ is $\mathcal{K}$-covered by $\phi_{i,k}$ with $i=1,\ldots, n$; in other words, the number of maps needed to $\mathcal{K}$-cover each $u_k$ is independent of $k$.  Furthermore, we say a uniformly $\mathcal{K}$-covered sequence is a \emph{uniformly robust $\mathcal{K}$-covered sequence} provided there exists $\epsilon>0$ and a compact set $\widetilde{\mathcal{K}}\subset M$ with the properties that $\mathcal{K}$ is contained in the interior of $\widetilde{\mathcal{K}}$,  and
\begin{equation*}
u_k\big(S_k\setminus \cup_{i=1}^n\phi_{i,k}(\mathcal{D}_{r'})\big)  \subset M \setminus \widetilde{\mathcal{K}};
\end{equation*}
for all $r'\in (r-\epsilon,r)$.  We call the $(\phi_{i,k},\mathcal{D}_r)$ \emph{uniformly robust $\mathcal{K}$-covers} of the sequence $(u_k,S_k)$.
\end{definition}

\begin{remark}
Note that a uniformly robust $\mathcal{K}$-covered sequence $(u_k,S_k)$ with $\mathcal{K}$-covers $\phi_{i,k}:\mathcal{D}_r\to S_k$, has two convenient properties: first the $(u_k,S_k)$ are a robustly $\mathcal{K}$-proper sequence, and second for all sufficiently large $r'< r$ (independent of $k$), the restricted maps $\phi_{i,k}:\mathcal{D}_{r'}\to S_k$ again form  uniformly robust $\mathcal{K}$-covers for the sequence $(u_k,S_k)$.
\end{remark}

\begin{definition}[$\mathcal{K}_{loc}$-convergence]\label{def:KlocConvergence}
Given a uniformly robust $\mathcal{K}$-covered sequence $u_k:S_k\to M$, we say that the $u_k$ converge in a smooth $\mathcal{K}_{loc}$ sense provided there exists a sequence of uniformly robust $\mathcal{K}$-covers $\phi_{i,k}:\mathcal{D}_r\to S_k$ with the property that for each $i=1,\ldots,n$ the maps $u_k\circ\phi_{i,k}:\mathcal{D}_r\to M$ converge in $C^\infty(\mathcal{D}_r,M)$.  We say the limit is immersed provided each $\tilde{u}_{i,\infty}:=\lim u_k\circ \phi_{i,k}$ is immersed.
\end{definition}

It is instructive to point out that smooth $\mathcal{K}_{loc}$-convergence in general does not imply smooth $\mathcal{K}$-convergence.  This is due to the fact that $\mathcal{K}_{loc}$-convergence does not guarantee any sort of topological convergence of the underlying $S_k$.  Consider for instance, a sequence of double covers of $\mathbb{S}^1$, for which the domains alternate between being connected and disconnected.  Nevertheless, given $\mathcal{K}_{loc}$-convergence, one expects that after passing to a subsequence $\mathcal{K}$-convergence can be obtained.  Indeed, this is the content of Proposition \ref{prop:KlocImpliesK} below.

\begin{proposition}\label{prop:KlocImpliesK}
Let $M$ be a manifold, and $\mathcal{K}\subset \Int(M)$.  Let $u_k:S_k\to M$ be a uniformly robust $\mathcal{K}$-covered sequence which smoothly $\mathcal{K}_{loc}$ converge to an immersed limit. Then a subsequence robustly $\mathcal{K}$ converges in $C^\infty$.
\end{proposition}

The proof of Proposition \ref{prop:KlocImpliesK} is provided in Section \ref{sec:KlocImpliesK}.  We now return to establishing some notation, and discussing some elementary properties of $J$-curves which will be exploited in later sections.

If $M$ is a manifold and $A\subset M$, then we will use the notation  $\mathcal{O}(A)$ to denote some open set containing $A$.  Furthermore, if $M$ is equipped with a metric $g$, then we will use the notation $\mathcal{O}_\delta^g(A):=\{p\in M: \dist_{g}(p,A)<\delta\}$
to denote a $\delta$-neighborhood of $A$. In the case that $A=p\in M$ is just a point, and $\delta>0$ is sufficiently small so that a $\delta$ neighborhood of $p$ is a ball, then we will use the notation $\mathcal{B}_\delta^g(p)=\mathcal{O}_\delta^g(p)$.  Also, recall that for a map $F:\mathcal{O}(0)\subset\R^{m}\to\R^{n}$, we say that $F(x)=O_\ell(|x|^k)$ provided $|D^\alpha F(x)|=O\big(|x|^{k-|\alpha|}\big)$ for all multi-indices $\alpha$ with $|\alpha|=0,\ldots \ell$.

\begin{definition}[generally immersed]
We shall say a smooth map $u:S\to M$ between smooth manifolds (which may have boundary and corners, be disconnected, or be non-compact) is a \emph{generally immersed} provided that for each point $z\in S$ for which $T_z u\neq 0$ we have $\Rank(T_z u)=\dim S$, and the set of critical points, which we henceforth denote as $\mathcal{Z}_u:=\{z\in S: T_z u =0\}$, has no accumulation points. Furthermore if $M$ is equipped with a Riemannian metric $g$, then we require that the conformal structure $[u^*g]$ on $S\setminus \mathcal{Z}_u$ admits a smooth extension across $\mathcal{Z}_u$.
\end{definition}

\begin{lemma}[local model]\label{lem:CritPointLocalModel}
Let $(M,J,g)$ be an almost Hermitian manifold, with $\mathcal{K}\subset \Int(M)$ a compact set.  Suppose $(u,S,j,J)$ is a robustly $\mathcal{K}$-proper generally immersed $J$-curve in $M$, and fix $z\in u^{-1}(\mathcal{K})$.  Then there exists a local holomorphic coordinate chart $\phi_z:\mathcal{O}(z)\to\mathcal{O}(0)\subset\mathbb{C}\simeq \mathbb{R}^2$, geodesic polar coordinates $\Phi_z:\mathcal{O}\big(u(z)\big)\to \mathcal{O}(0)\subset\mathbb{C}^n\simeq \mathbb{R}^{2n}$, and unique $k_{z}\in \mathbb{N}$ such that $\phi_z(z)=0$, $\Phi_z\big(u(z_0)\big)=0$, and such that
\begin{equation*}
\Phi_z\circ u\circ \phi_z^{-1} (\rho) = \big(\rho^{k_z},0,\ldots,0\big) + F_z(\rho),
\end{equation*}
where $F_z(\rho)=O_{k_z+1}(|\rho|^{k_z+1})$.
\end{lemma}

\begin{proof}
First, we will drop the $z$-dependence from our notation, and simply write $k$, $\phi$, $\Phi$, and $F$.  Next, let $\tilde{\phi}:\mathcal{O}(z)\to\mathcal{O}(0)\subset \mathbb{C}\simeq\R^2$ be a holomorphic coordinate chart for which $\tilde{\phi}(z)=0$, and let $\widetilde{\Phi}:\mathcal{O}\big(u(z)\big)\to\mathcal{O}(0)\subset\R^{2n}$ be polar geodesic coordinates for which $\widetilde{\Phi}\big(u(z)\big)=0$ and $(\widetilde{\Phi}_* J)(p) =:\widetilde{J}(p)= J_0 + O(|p|)$; here $J_0$ is the standard almost complex structure defined by $J_0\partial_{x_\alpha}=\partial_{y_\alpha}$ for $\alpha=1,\ldots,n$.  Next, recall a consequence of Aronszajn's theorem, which guarantees that if $\mathcal{O}\subset \R^2$ is open and connected, and $\tilde{u}:\mathcal{O}\to \R^m$ is a smooth map which satisfies
\begin{equation*}
|\Delta \tilde{u}| \leq C(|\tilde{u}|+|\partial_s \tilde{u}|+|\partial_t \tilde{u}|)
\end{equation*}
on $\mathcal{O}$, and $(D^\alpha \tilde{u})(0)=0$ for all multi-indices $\alpha$, then $\tilde{u}\equiv 0$ on $\mathcal{O}$; here we are using subscripts to denote partial differentiation.  Since $du + J(u)\cdot du \cdot j = 0$, it follows that for $\tilde{u}:=\widetilde{\Phi}\circ u\circ \tilde{\phi}^{-1}$ we have $\tilde{u}_s+\widetilde{J}(\tilde{u}) \tilde{u}_t = 0$, and hence
\begin{align*}
|\Delta \tilde{u}| &= \big|\big( -\widetilde{J}(\tilde{u})\tilde{u}_t\big)_s + \big(\widetilde{J}(\tilde{u})\tilde{u}_s\big)_t\big|\\
&=\big|-\big(\widetilde{J}(\tilde{u})\big)_s \tilde{u}_t + \big(\widetilde{J}(\tilde{u})\big)_t \tilde{u}_s\big|\\
&\leq C(|\tilde{u}_s| + |\tilde{u}_t|);
\end{align*}
here we have made use of the fact that the $C^1$ norms of $J$ and $u$ are uniformly bounded. By assumption $u$ is generally immersed, and hence $\tilde{u}$ is not a constant map, so it follows that $\tilde{u}(s,t) = P(s,t)+F(s,t)$, where $F(s,t)=O(|s+it|^{k+1})$ and $P$ is a homogeneous polynomial of degree $k\in\mathbb{N}$.

Next define the following linear maps
\begin{align}
&\ell_\alpha:\mathbb{R}^2\to\mathbb{R}^2&\; &\ell_\alpha x=\alpha x \label{eq:L_alpha}\\
&L_\alpha:\mathbb{R}^{2n}\to\mathbb{R}^{2n}& \; & L_\alpha(x)=\alpha x.\label{eq:ell_alpha}
\end{align}
Observe that $(L_\epsilon^* \widetilde{J}) (p)= \widetilde{J}(\epsilon p) = J_0 + \epsilon O(|p|)$, and thus as $\epsilon\to 0$ we have $(L_{\epsilon}^* \widetilde{J}) \to J_0$ in $C^\infty$. Let us also define the maps $v_\epsilon:= L_{\epsilon^{-k}}\circ \tilde{u} \circ \ell_\epsilon$ so that
\begin{align*}
v_\epsilon(s,t)&=P(s,t)+\epsilon^{-k} F(\epsilon s,\epsilon t)\\
&=P(s,t) + \epsilon O(|s+it|^{k+1}),
\end{align*}
and thus $v_\epsilon \to P$ in $C^\infty$.  Also observe that
\begin{equation*}
\partial_s v_\epsilon + (L_{\epsilon^k}^* \widetilde{J})(v_\epsilon) \partial_t v_\epsilon =  \epsilon^{1-k} \tilde{u}_s\circ \ell_\epsilon +\epsilon^{1-k}\widetilde{J}(\tilde{u}\circ\ell_\epsilon)\tilde{u}_t\circ \ell_\epsilon =0;
\end{equation*}
this together with the fact that $v_\epsilon \to P$ and $(L_{\epsilon^k}^*\widetilde{J})\to J_0$, it follows that $P_s + J_0 P_t =0$, and hence $P$ has the form
\begin{equation*}
P(s,t)= c H \cdot \big(\re((s+it)^k),\im ((s+it)^k),0,\ldots,0\big),
\end{equation*}
where $c\in\R^+$, and $H\in \R^{2n\times 2n}$ is a real matrix for which $H^T H =\mathds{1}$ and $J_0 H = H J_0$. Consequently for $\Phi:=H^{-1}\cdot\widetilde{\Phi}$ and $\phi:=c^{1/k}\tilde{\phi}$ the lemma is proved.
\end{proof}

In light of Lemma \ref{lem:CritPointLocalModel}, it will be convenient to make the following definition.

\begin{definition}\label{def:ord}
Let $\mathbf{u}$ be a generally immersed $J$-curve.  Then for any interior point $z_0$, we define the order of $z_0$ to be the following:
\begin{equation*}
\order(z_0)=k-1
\end{equation*}
where $k$ is the integer guaranteed by Lemma \ref{lem:CritPointLocalModel}.
\end{definition}

Since much of the analysis that follows will regard $J$-curves as sub-manifolds, we take a moment to establish some convenient notation for certain pull-back bundles associated to a given immersion $u:S\to M$ with image in a Riemannian manifold $(M,g)$:
\begin{align*}
u^*TM&:=\{(\zeta,X)\in S\times T_pM: u(\zeta)=p\}\\
\mathcal{T}&:=\{(\zeta,X)\in u^*TM: X \in  Tu (T_\zeta S)\}\\
\mathcal{N}&:=\{(\zeta,X)\in u^*TM: \la X, Y\ra_g =0\;\; \forall\; (\zeta,Y)\in \mathcal{T}\}.
\end{align*}
We also define the second fundamental form $\two_u$ along the image of $u$ by the following.
\begin{equation*}
\two_u\in \Gamma\big(\Hom(\mathcal{T}\times\mathcal{T}, \mathcal{N})\big)\qquad\text{given by}\qquad \two_u(X,Y) = \big(\nabla_X Y)^\bot,
\end{equation*}
where $X,Y$ are sections of $\mathcal{T}\subset u^*TM$, and $\nabla$ is the Levi-Cevita connection on $u^*TM$ induced from $TM$, and $X\mapsto X^\bot$ is the $g$-orthogonal projection from $u^*TM$ to $\mathcal{N}$.  Recall that the mean curvature vector $H_\nu\in\Gamma (\mathcal{N})$ of an immersion $u:S\to M$ is given by
\begin{equation*}
H_\nu:=\sum_{i=1}^{\dim S}B_u(e_i,e_i)
\end{equation*}
where $\{e_1,\ldots,e_{\dim S}\}$ is any orthonormal frame in $\mathcal{T}$.  Recall that if the almost symplectic form $\omega=g\circ(J\times \mathds{1})$ is actually symplectic (i.e. $d\omega = 0$), then $J$-curves are minimal surfaces -- or more precisely generalized minimal immersions.  However, when $\omega$ is not closed, then immersed $J$-curves satisfy the mean curvature equation
\begin{equation}\label{eq:meanCurvatureEquation}
H_{\nu} = \tr_S Q
\end{equation}
where $Q:=J \nabla J$, and $\tr_S Q$ is the trace  $\tr_S Q:= Q(e,e) + Q(f,f)$ where $e,f\in \mathcal{T}$ form an orthonormal frame.  Consequently, we can recall the Gauss equations for two-dimensional immersions $u:S\to M$ are
\begin{equation}\label{eq:gaussEquationsGeneral}
K_{sec}(\mathcal{T}_{u(\zeta)}) = K_g(\zeta) -\la B_u(e,e),B_u(f,f)\ra_g + \|B_u(e,f)\|_g^2,
\end{equation}
which reduce to the following when $(u,S,j,J)$ is an immersed $J$-curve:
\begin{equation}
K_{sec}(\mathcal{T}_{u(\zeta)}) + {\textstyle \frac{1}{2}}\|\tr_S Q\|_g^2 = K_{u^*g}(\zeta) + {\textstyle \frac{1}{2}}\|B_u\|_g^2.
\end{equation}

Next we wish to define Gromov convergence of $J$-curves, however to do this we need some preliminary definitions; here we will essentially follow Sections 4 and 7 in \cite{BEHWZ03}. To that end, we define a \emph{marked} $J$-curve to be a pair $(\mathbf{u},\mu)$ where $\mathbf{u}=(u,S,j,J)$ is a $J$-curve and $\mu\subset S\setminus \partial S$ is a finite set of points called marked points.

A \emph{nodal} $J$-curve is a triple $(\mathbf{u},\mu,D)$ where $(\mathbf{u},\mu)$ is a marked $J$-curve, and $D$ is an unordered finite set of pairs of distinct points $D=\{\overline{d}_1,\underline{d}_1,\ldots,\overline{d}_\delta,\underline{d}_\delta\}\subset S\setminus \partial S$ with the property that $u(\overline{d}_i)=u(\underline{d}_i)$ for $i=1,\ldots,\delta$ and $\mu\cap D=\emptyset$.  As in Section 4.4 of \cite{BEHWZ03}, we define $S^D$ to be the oriented blow-up of $S$ at the points $D$, and we let $\overline{\Gamma}_i:=\big( T_{\overline{d}_i}(S)\setminus \{0\}\big)\setminus \R_+^*\subset S^D$ and $\underline{\Gamma}_i:=\big( T_{\underline{d}_i}(S)\setminus \{0\}\big)\setminus \R_+^*\subset S^D$
denote the newly created boundary circles over the $d_i$. Furthermore, we say a nodal $J$-curve is \emph{stable} provided that
for each connected component $\tilde{S}$ of $S$ we have $3\leq 2\Genus(\tilde{S})+\#(\tilde{\mu}\cup \tilde{D})$ where $\tilde{\mu}=\tilde{S}\cap \mu$ and $\tilde{D}:=\tilde{S}\cap D$.  Note that in the case that $\tilde{S}$ is compact, then this condition is equivalent to $\chi(\tilde{S})-\#(\tilde{\mu}\cup \tilde{D}) < 0$, so that there exists a unique complete finite area hyperbolic metric of constant curvature $-1$ on $S':=S\setminus(\mu\cup D)$ which is in the same conformal class as $j$ and for which each connected component of $\partial S$ is a geodesic; we denote this metric by $h^{j,\mu\cup D}$.

A \emph{decorated nodal} $J$-curve $(\mathbf{u},\mu,D,r)$ is a quadruple for which $(\mathbf{u},\mu,D)$ is a nodal $J$-curve and $r$ is a set of orientation reversing orthogonal maps $r_i:\overline{\Gamma}_i\to \underline{\Gamma}_i$, which we call \emph{decorations}.
We also define $S^{D,r}$ to be the smooth surface obtained by gluing the components of $S^D$ along the boundary circles $\{\overline{\Gamma}_1,\underline{\Gamma}_1,\ldots,\overline{\Gamma}_\delta,\underline{\Gamma}_\delta\}$ via the decorations $r_i$.  We will let $\Gamma_i$ denote the special circles $\overline{\Gamma}_i=\underline{\Gamma}_i\subset S^{D,r}$.  Observe that the smooth map $u:S\to M$ then lifts to a continuous map $u:S^{D,r}\to M$.

\begin{definition}[Gromov convergence]\label{def:GromovConvergence}
A sequence $(u_k,S_k,j_k,J_k)$ of compact $J_k$-curves (i.e. potentially with boundary)
is said to converge in a Gromov-sense to a nodal $J$-curve with boundary $(\mathbf{u},D)$ with $\mathbf{u}=(u,S,j,J)$,
provided the following are true for all sufficiently large $k\in \mathbb{N}$.
\begin{enumerate}
\item $J_k\to J$ in $C^\infty$.
\item There exist sets of marked points $\mu_k\subset S_k\setminus \partial S_k$ and $\mu\subset S\setminus (\partial S \cup D)$ with the property that $\#\mu = \# \mu_k$ for all $k$, and the marked $J$-curves $(\mathbf{u}_k,\mu_k)$ and the nodal curve $(\mathbf{u},\mu,D)$ are all stable. We further require that if $\tilde{S}$ is a connected component of $S$ and $u:\tilde{S}\to M$ is a constant map, then $3\leq 2\Genus(\tilde{S})+\# (D\cap \tilde{S})$.
\item There exist a decoration $r$ for $(\mathbf{u},D)$ and sequences of diffeomorphisms $\phi_k: S^{D,r}\to S_k$ such that $\phi_k(\mu) = \mu_k$ and for each $i=1,\ldots,\delta$ the curve $\phi_k(\Gamma_i)$ is a $h^{j_k,\mu_k}$-geodesic in $S_k'$.
\item $\phi_k^* h^{j_k,\mu_k} \to h^{j,\mu\cup D}$ in $C_{loc}^\infty\big(S^{D,r}\setminus (\mu\cup_i \Gamma_i) \big)$; here we have abused notation by letting $h^{j,\mu\cup D}$ also denote its lift to $S^{D,r}$.
\item $\phi_k^*u_k \to u$ in $C^0(S^{D,r})$.
\item $\phi_k^*u_k \to u$ in $C_{loc}^\infty(S^{D,r}\setminus \cup_i\Gamma_i)$.
\end{enumerate}
\end{definition}

With this definition in hand, we finish this section by defining the notion of robust $\mathcal{K}$-convergence in a Gromov sense.

\begin{definition}[robust $\mathcal{K}$-convergence in Gromov sense]
Consider an almost Hermitian manifold given by $(M,J,g)$ and a sequence of almost Hermmitian structures $(J_k,g_k)$ for which $(J_k,g_k)\to (J,g)$ in $C^\infty$, and a compact set $\mathcal{K}\subset \Int(M)$, and a robustly $\mathcal{K}$-proper sequence of generally immersed $J_k$-curves $\mathbf{u}_k=(u_k,S_k,j_k,J_k)$. We say that the $\mathbf{u}_k$ \emph{robustly $\mathcal{K}$-converge in a Gromov sense} provided there exists a compact set $\widetilde{\mathcal{K}}\subset \Int(M)$ for which $\mathcal{K}\subset \Int(\widetilde{\mathcal{K}})$, and there exist
compact regions $\tilde{S}_k\subset S_k$ with the property that $u_k(S_k\setminus \tilde{S}_k)\subset M\setminus \widetilde{\mathcal{K}}$, and the domain restricted $J_k$-curves $(u_k,\tilde{S}_k,j_k,J_k)$ converge in a Gromov sense. We additionally require that the sequence of marked points added to the $(\tilde{S},j_k)$ to obtain Gromov convergence are chosen so that lengths of each connected component of $\partial \tilde{S}_k$ (computed with respect to the associated Poincar\'{e} metric) are uniformly bounded away from zero and infinity. Moreover we require that each component of the limit curve with non-empty boundary is non-constant, and $u_\infty$ restricted to some neighborhood of the boundary is an immersion.
\end{definition}

\section{Target-Local Compactness}\label{sec:Main}

The goal of this section is to prove Theorem \ref{thm:Main2} below, which is the main result of this article.  Also of importance in this section is the proof of Corollary \ref{cor:Main} below, which is a restatement of Theorem A from the introduction.

\begin{theorem}\label{thm:Main2}
Let $(M,J,g)$ be an almost Hermitian manifold, and let $(J_k,g_k)$ be a sequence of almost Hermitian structures which converge in $C^\infty$ to $(J,g)$. Also let $\mathcal{K}\subset \Int(M)$ be a compact region, and let $\mathbf{u}_k$ be a sequence of generally immersed $J_k$-curves which are robustly $\mathcal{K}$-proper and satisfy
\begin{enumerate}
\item $\Area_{u_k^*g_k}(S_k)\leq C_A <\infty$
\item $\Genus(S_k)\leq C_G <\infty$.
\end{enumerate}
Then a subsequence robustly $\mathcal{K}$-converges in a Gromov sense.
\end{theorem}

The proof of Theorem \ref{thm:Main2} consists of three main steps. The first step is to prove Theorem \ref{thm:Main2} with the additional assumptions that the curves are immersed and $\|\two_{u_k}\|_{L^\infty}$ is uniformly bounded, but without the assumption of bounded topology; this is the content of Section \ref{sec:CompactnessWithLInftyBounds}.  The second step is to use this result to prove  Theorem \ref{thm:Main2} with the additional assumptions that $\|\two_{u_k}\|_{L^2}$ is uniformly bounded and that the number of critical points of the $u_k$ are uniformly bounded; this is the content of Section \ref{prop:compactnessFromL2Bounds}.  Finally, the third step is to use this result to prove Theorem \ref{thm:Main2} with no additional assumptions.

\subsection{Compactness with $\|\two\|_{L^\infty}$ bounds}\label{sec:CompactnessWithLInftyBounds}

In this section we prove the following result.
\begin{proposition}\label{prop:compactnessFromLInfinityBounds}
Let $(M,J,g)$ be a compact almost Hermitian manifold with boundary, and let $(J_k,g_k)$ be a sequence of almost Hermitian structures which converge in $C^\infty$ to $(J,g)$. Also let $\mathcal{K}\subset \Int(M)$ be a compact region, and let $\mathbf{u}_k$ be a sequence of immersed compact $J_k$-curves which are robustly $\mathcal{K}$-proper.  Suppose further that
\begin{enumerate}
\item $\Area_{u_k^*g_k}(S_k) \leq C_A <\infty$
\item $\sup_{\zeta\in S_k}\|\two_{u_k}^{g_k}(\zeta)\|_{g_k} \leq C_{\two} <\infty$.
\end{enumerate}
Then a subsequence robustly $\mathcal{K}$-converges.  Here $B_{u_k}^{g_k}$ denotes the second fundamental form the the immersions $u_k:S_k\to M$ computed with respect to the metrics $g_k$ on $M$.
\end{proposition}
\begin{proof}
We note that as a consequence of Proposition \ref{prop:KlocImpliesK}, it is sufficient to show that a subsequence robustly $\mathcal{K}_{loc}$-converges.  Consequently, we need some convenient local parameterizations.  In particular we will consider local graphical parameterizations over coordinate tangent planes.  We make this precise with the following.

\begin{proposition}[Uniform Local Graphs]\label{prop:uniformLocalGraphs}
Let $M$ be a be a compact manifold of dimension $2n$ and possibly with boundary.  Let $(J_k,g_k)\to (J,g)$ be a sequence of almost Hermitian structures on $M$ which converge in $C^\infty$.  Fix a compact set $\mathcal{K}\subset\Int(M)$, and a constant $C_{\two}>0$. Then there exist positive constants $r_0, C_0, C_1,C_2, C_3\ldots,$ depending only on $C_B$, $\dist_g(\mathcal{K},\partial M)$,  and the geometry of $(M,J,g)$ with the following significance.
For each proper immersed $J_k$-curve denoted by $(u_k,S_k,j_k,J_k)$ for which
\begin{equation*}
\partial S_k= u_k^{-1}(\partial M),\qquad \text{and} \qquad \sup_{\zeta \in S_k} \|\two_{u_k}^{g_k}(\zeta)\| \leq C_{\two},
\end{equation*}
and each $\zeta\in S_k$ such that $u_k(\zeta)\in \mathcal{K}$
there exists a map $\phi:\mathcal{D}_{r_0}\to S_k$ and geodesic normal coordinates $\Phi:\mathcal{B}_{2r_0}^{g_k}\big(u(\zeta)\big)\to\mathbb{R}^{2n}$ with the following properties.
\begin{enumerate}
\item $\tilde{u}(s,t):=\Phi\circ u\circ \phi (s,t) = \big(s,t,\tilde{u}^3(s,t),\ldots,\tilde{u}^{2n}(s,t)\big)$,
\item $\phi(0)=\zeta$, $\tilde{u}(0,0) = 0$, and $D_\alpha \tilde{u}^i (0,0) = 0$ for $|\alpha|=1$ and $i=3,\ldots,2n$,
\item $\sum_{|\alpha|=1}\sum_{i=3}^{2n}\|D_\alpha \tilde{u}^i\|_{C^0(\mathcal{D}_{r_0})}^2 \leq 10^{-20}$ and $\|\tilde{u}\|_{C^k(\mathcal{D}_{r_0})} \leq C_k$ for $k\in \mathbb{N}$.

\item For Euclidian coordinates $\rho=(s,t)$, on $\mathcal{D}_{r_0}$, we have
\begin{equation*}
{\textstyle \frac{1}{2}}|\rho|\leq \dist_{(u_k\circ \phi)^*g_k}(0,\rho) \leq 2 |\rho|\quad \text{and} \quad {\textstyle \frac{1}{2}}|\rho|\leq \dist_{g_k}\big(u_k(\phi(\rho)),u_k(\phi(0))\big)\leq 2|\rho|.
\end{equation*}
\end{enumerate}
\end{proposition}
A proof of Proposition \ref{prop:uniformLocalGraphs} can be found in \cite{Fj09b}; the idea of the proof goes as follows.  First one shows that $J$-curves satisfy an inhomogenous mean curvature equation of the form $H=\tr_S Q$ with $Q$ a tensor on $M$. Next one writes this equation in local coordinates on $M$ to see that locally the $u$ solve a second order partial differential equation. The uniform curvature bound guarantees that in geodesic normal coordinates, in a small disk centered at $\zeta$ tangent planes don't deviate too much from being ``horizontal.''  One concludes the existence of a graphical parameterization, in which case the partial differential equation that the graphically (but not holomorphically) parameterized $J$-curves solve is uniformly elliptic.  One readily sees that uniform curvature bounds then guarantee uniform $C^2$ bounds, in which case the uniform bounds on the $\|D_\alpha \tilde{u}^i\|$ with $|\alpha| > 2$ then follows from the usual elliptic regularity theory.

In order to prove robust $\mathcal{K}_{loc}$-convergence we must now show that the parameterizations of Proposition \ref{prop:uniformLocalGraphs} can be used to construct a uniformly robust $\mathcal{K}$-cover.  The desired convergence will then follow essentially from the Arzel\`{a}-Ascoli theorem.  To construct the desired  $\mathcal{K}$-cover, we first recall the extrinsic monotonicty of area lemma.

\begin{proposition}[monotonicity of area]\label{prop:monotonicityOfArea}
Let $(M,J,g)$, be a compact almost Hermitian manifold possibly with boundary.  Then for all $(J',g')$ sufficiently close to $(J,g)$ in a $C^2$-sense, the following holds.  Let $(u,S,j,J')$ be a compact generally immersed pseudo-holomorphic curve for which $u(\partial S)\cap \mathcal{O}_r^{g'}\big(u(\zeta)\big)=\emptyset$ for some $r>0$ satisfying
\begin{equation*}
8r< \min\big(C^{-1},\inj_{g}\big(u(\zeta)\big)\big),
\end{equation*}
where
\begin{equation*}
\sup_{p\in M} |K_{sec}^g(p)| \leq {\textstyle \frac{1}{4}}C^2\qquad\text{and}\qquad \sup_{p\in M}\|\nabla J\|_g \leq {\textstyle\frac{1}{4}}C;
\end{equation*}
here $\nabla$ is the Levi-Civita connection associated to $g$, and $|K_{sec}^g(p)|$ defined by
\begin{equation}\label{eq:defOfAbsKsec}
|K_{sec}^g(p)|:=\sup\{|K_{sec}^g(X,Y)| : X,Y\in T_p M\text{ and } X\wedge Y\neq 0\}.
\end{equation}
Then for all $0<a<b\leq r$ we have
\begin{equation*}
\frac{1}{a^2} \Area_{u^*g'}\big(S_a(\zeta)\big)\leq  \frac{2}{b^2}\Area_{u^*g'}\big(S_b(\zeta)\big),
\end{equation*}
where $S_a(\zeta)$ the connected component of $u^{-1}\big(\mathcal{O}_a^{g'}(u(\zeta))\big)$ which contains $\zeta$. In particular, letting $a\to 0$ and $b=r$ yields the familiar result
\begin{equation}\label{eq:MonotonicityStandard}
\frac{\pi r^2}{2}\leq  \Area_{u^*g'}\big(S_r(\zeta)\big).
\end{equation}
\end{proposition}

A proof of the above proposition can be found in \cite{Fj09b}; it is a modification of the well known result for minimal surfaces. Also note that the weaker version of monotonicity given in (\ref{eq:MonotonicityStandard}), is a very well known result for $J$-curves (c.f. \cite{Gm85} , \cite{Hc97}, \cite{Mmp94}), and it is sufficient for our purposes the remainder of this article.
We now prove a fairly standard covering result.

\begin{lemma}\label{lem:coveringLemma}
Let $\mathcal{K}$,  $M$, $(J_k,g_k)$, and $\mathbf{u}_k$, be as in the hypotheses of Proposition \ref{prop:compactnessFromLInfinityBounds}.  Then after passing to a subsequence, a robust uniform $\mathcal{K}$-cover can by obtained using only the graphical parameterizations $\phi$ given by Proposition \ref{prop:uniformLocalGraphs}.
\end{lemma}

\begin{proof}
We begin by fixing two auxiliary compact regions\footnote{See Definition \ref{def:CompactRegion}.} $\widehat{\mathcal{K}},\widetilde{\mathcal{K}}\subset M$ for which $\mathcal{K}\subset \Int(\widehat{\mathcal{K}})$, $\widehat{\mathcal{K}}\subset \Int(\widetilde{\mathcal{K}})$, and $\widetilde{\mathcal{K}}\subset\Int(M)$, and for which the $\mathbf{u}_k$ are robustly $\widetilde{\mathcal{K}}$-proper.  Observe that the functions defined by $f_k(p):=\big(\dist_{g_k}(p,\widehat{\mathcal{K}})\big)^2$, are all smooth in a neighborhood of the form $\mathcal{O}(\widehat{\mathcal{K}})\setminus \widehat{\mathcal{K}}$. Observe that the functions $f_k\circ u_k$ have critical values which are the compliment of an open dense set in $(0,\epsilon)$.  It follows that there exists some $\epsilon_0>0$ which is not a critical value of any of the $u_k$, and hence the $J_k$ curves $\big(u_k,u_k^{-1}(\{f\leq \epsilon_0\}), j_k,J_k)$ satisfy the hypotheses of Proposition \ref{prop:uniformLocalGraphs}.  Without loss of generality, we will henceforth assume that $S_k=u_k^{-1}(\{f\leq \epsilon_0\})$. Then for each $\zeta\in u_k^{-1}(\widehat{\mathcal{K}})$ we let $\phi_{\zeta,k}:\mathcal{D}_{r_0}\to S_k$ and  $\Phi_\zeta:\mathcal{O}_{2r_0}^{g_k}\big(u_k(\zeta)\big)\to\mathbb{R}^{2n}$ denote the maps guaranteed by Proposition \ref{prop:uniformLocalGraphs}.  Next, for each $k$, choose $\{\zeta_{1,k},\zeta_{2,k},\ldots \zeta_{m_k,k}\}\subset u_k^{-1}(\widehat{\mathcal{K}})$ so that the open sets $\mathcal{O}_{r_0/8}^{u_k^*g_k}(\zeta_{1,k}),\mathcal{O}_{r_0/8}^{u_k^*g_k}(\zeta_{2,k}),\ldots,\mathcal{O}_{r_0/8}^{u_k^*g_k}(\zeta_{m_k,k})$ are maximally disjoint.  For clarity, we define $\phi_{i,k}:=\phi_{\zeta_{i,k},k}$ so that $\phi_{i,k}(0)=\zeta_{i,k}$. We now observe that to complete the proof of Lemma \ref{lem:coveringLemma}, it is sufficient to prove the following two claims: firstly there exists an $m\in \mathbb{N}$ such that $m_k \leq m$ for all sufficiently large $k$, and secondly
\begin{equation}\label{eq:coveringArgument1}
u_k^{-1}(\widehat{\mathcal{K}})\subset \bigcup_{i=1}^{m_k}\phi_{i,k}(\mathcal{D}_{r_0/2}).
\end{equation}
We prove the former statement first.  Indeed, recall that by Proposition \ref{prop:uniformLocalGraphs}, we have
\begin{equation*}
u_k\circ \phi_{i,k} (\partial \mathcal{D}_{r_0}) \cap \mathcal{O}_{r_0/2}^{g_k}\big(u_k(\zeta_{i,k})\big) =\emptyset,
\end{equation*}
so by Lemma \ref{prop:monotonicityOfArea} it follows that
\begin{equation*}
\pi (r_0/16)^2/2\leq \Area_{u_k^*g_k}\big(S_{r_0/16}(\zeta_{i,k})\big),
\end{equation*}
where $S_{r_0/16}(\zeta_{i,k})$ is the connected component of $u_k^{-1}\big(\mathcal{B}_{r_0/16}^{g_k}(u_k(\zeta_{i,k}))\big)$ containing $\zeta_{i,k}$.
Again by Proposition \ref{prop:uniformLocalGraphs}, one finds that
\begin{equation*}
S_{r_0/16}(\zeta_{i,k})\subset \mathcal{O}_{r_0/8}^{u_k^*g_k}\big(\phi_{i,k}(0)\big).
\end{equation*}
Since the these latter sets are disjoint, it follows that
\begin{equation*}
 2^{-9}m_k\pi r_0^2 \leq \sum_{i=1}^{m_k}\Area_{u_k^*g_k}\big( \mathcal{O}_{r_0/8}^{u_k^*g_k}\big(\phi_{i,k}(0)\big)\big) \leq \Area_{u_k^*g_k}(S_k) \leq C_A,
\end{equation*}
and thus the $m_k$ are uniformly bounded.  To prove the latter statement, namely the containment (\ref{eq:coveringArgument1}), we first observe that as a consequence  of Proposition \ref{prop:uniformLocalGraphs}, it follows that $\mathcal{O}_{r_0/4}^{u_k^*g_k}(\zeta_{i,k})\subset \phi_{i,k}(\mathcal{D}_{r_0/2})$, and thus to prove (\ref{eq:coveringArgument1}) it is sufficient to show that
\begin{equation}\label{eq:coveringArgument2}
u_k^{-1}(\widehat{\mathcal{K}})\subset \bigcup_{i=1}^{m_k}\mathcal{O}_{r_0/4}^{u_k^*g_k}(\zeta_{i,k}).
\end{equation}
To see this, we suppose not.  Then there exists $\zeta\in u_k^{-1}(\widehat{\mathcal{K}})$ such that
\begin{equation*}
\min_{1\leq i \leq m_k} \dist_{u_k^*g_k}\big(\zeta, \zeta_{i,k}\big)\geq r_0/4.
\end{equation*}
However, it then follows that for all $i=1,\ldots,m_k$ we have $\mathcal{O}_{r_0/8}^{u_k^*g_k}(\zeta)\cap \mathcal{O}_{r_0/8}^{u_k^*g_k}(\zeta_{i,k})=\emptyset$, but by the maximality of the $\mathcal{O}_{r_0/8}^{u_k^*g_k}(\zeta_{1,k}), \mathcal{O}_{r_0/8}^{u_k^*g_k}(\zeta_{2,k}),\ldots,\mathcal{O}_{r_0/8}^{u_k^*g_k}(\zeta_{n_k,k})$, we must have $\mathcal{O}_{r_0/8}^{u_k^*g_k}(\zeta)\cap \mathcal{O}_{r_0/8}^{u_k^*g_k}(\zeta_{i,k})\neq\emptyset$ for some $i\in \{1,\ldots,m_k\}$.  This contradiction proves (\ref{eq:coveringArgument2}), and hence completes the proof of Lemma \ref{lem:coveringLemma}.
\end{proof}

With Lemma \ref{lem:coveringLemma} in hand, we now complete the proof of Proposition \ref{prop:compactnessFromLInfinityBounds}.  To that end, we note that it is sufficient to prove that for each $i=1,\ldots, m$ a subsequence of the maps $u_k\circ \phi_{i,k}$ converges in $C^\infty(\mathcal{D}_{r_0-\epsilon'},M)$ for some small $\epsilon'\in (0,r_0/2)$.  To see this, we first note that after passing to a subsequence we arrange that for each $i=1,\ldots,m$ the sequence of points $u_k\circ\phi_{i,k}(0)$ converges, as well as linear maps $T_0\Phi_{i,k}^{-1}:T_0\mathbb{R}^{2n}\to T_{u_k\circ \phi_{i,k}(0)}M$; here, as before with $\phi$, we have let $\Phi_{i,k}:=\Phi_{\zeta_{i,k},k}$ be the geodesic polar coordinates guaranteed by Proposition \ref{prop:uniformLocalGraphs}.  By that same proposition, all the derivatives of the maps $\Phi_{i,k}\circ u_k\circ \phi_{i,k}$ are uniformly bounded, and hence by the Arzel\`{a}-Ascoli theorem, it follows that after passing to a further subsequence the $u_k\circ\phi_{i,k}$ converge in $C^\infty(\mathcal{D}_{r_0-\epsilon'}, M)$.  Furthermore we have shown that the $\phi_{i,k}:\mathcal{D}_{r_0-\epsilon'}\to S_k$ form a uniform robust $\mathcal{K}$-cover, and hence we have passed to a subsequence for which the $u_k$ robustly $\mathcal{K}_{loc}$-converge.
The proof of Proposition \ref{prop:compactnessFromLInfinityBounds} now follows from Proposition \ref{prop:KlocImpliesK}.
\end{proof}

\subsection{Compactness with $\|\two\|_{L^2}$ bounds}\label{prop:compactnessFromL2Bounds}
The purpose of this section is to prove Theorem \ref{thm:Main1} below.

\begin{theorem}\label{thm:Main1}
Let $(M,J,g)$, $(J_k,g_k)$, and $\mathcal{K}$, be as in the hypotheses of Proposition \ref{prop:compactnessFromLInfinityBounds}. Let $\mathbf{u}_k=(u_k,S_k,j_k,J_k)$ be a robust $\mathcal{K}$-proper sequence of compact $J_k$-holomorphic curves which satisfy the following
\begin{enumerate}
\item $\Genus(S_k) \leq C_G < \infty$,
\item $\Area_{u_k^*g_k}(S_k)\leq C_A<\infty$,
\item $\#\mathcal{Z}_{u_k}\leq C_{\mathcal{Z}}<\infty$,
\item $\int_{S_k}\|\two_{u_k}^{g_k}\|_{g_k}^2 \leq C_{Total}<\infty$,
\end{enumerate}
where $C_G$, $C_A$, $C_{\mathcal{S}}$, and $C_{Total}$ do not depend on $k$.  Then a subsequence robustly $\mathcal{K}$-converges in a Gromov sense.
\end{theorem}

\begin{proof}
We begin by letting $\widetilde{\mathcal{K}}\subset \Int(M)$ be a compact set for which $\mathcal{K}\subset \Int(\widetilde{\mathcal{K}})$, and for which the $\mathbf{u}_k$ are robustly $\widetilde{\mathcal{K}}$-proper.  Next we observe that since $\#\mathcal{Z}_{u_k}\leq C_{\mathcal{Z}}$, it follows that after passing to a subsequence, we may assume that $\mathcal{Z}_{u_k}=\{z_{1,k},z_{2,k},\ldots,z_{n_0,k}\}$, and for each $i=1,\ldots,n_0$ either the sequence $u_k(z_{i,k})$ converges to a point in $\widetilde{\mathcal{K}}$ or else $\dist_{g_k}\big(\widetilde{\mathcal{K}},u_k(z_{i,k})\big)\geq \delta >0$; we denote the associated limit set in $\widetilde{\mathcal{K}}$ by $\mathcal{S}_1$.  Next we claim the following.

\begin{lemma}[finite points of curvature blowup]\label{lem:finiteCurvatureBlowup}
Let $(M,J,g)$, $(J_k,g_k)$, $\mathcal{K}$, and $\mathbf{u}_k=(u_k,S_k,j_k,J_k)$ be as in Theorem \ref{thm:Main1}. Fix a compact set $\mathcal{K}_0\subset M$ such that $\mathcal{K}\subset \Int(\mathcal{K}_0)$ and $\mathcal{K}_0\subset \Int(\widetilde{\mathcal{K}})$.  Then after passing to a subsequence, there exists a finite set $\mathcal{S}_2\subset \mathcal{K}_0$ with the following significance. For each $\epsilon>0$, there exists $C>0$ and $k_0\in\mathbb{N}$ such that for all $k\geq k_0$, the following holds:
\begin{equation*}
\sup_{\zeta\in u_k^{-1}\big(\mathcal{K}_0\setminus \mathcal{O}_\epsilon^{g}(\mathcal{S}_2)\big)}\|\two_{u_k}^{g_k}(\zeta)\|^2 \leq C
\end{equation*}
\end{lemma}

\begin{proof}
The proof of Lemma \ref{lem:finiteCurvatureBlowup} has one major technical component, which we now state.

\begin{proposition}[Curvature Threshold]\label{prop:CurvatureThreshold}
Let $(M,J,g)$ be a compact almost Hermitian manifold possibly with boundary, and let $\epsilon>0$. Then
for all $(J',g')$ sufficiently close to $(J,g)$ in a $C^3$ sense, there exists an $\hbar>0$ depending on $\epsilon$ and the geometry of $(M,J,g)$ with the following significance.  If $(u,S,j,J')$ is a compact immersed $J'$-curve, with $\zeta\in S$ satisfying $\dist_{g'}\big(u(\zeta),\partial M\big)\geq \epsilon$, and $u(\partial S) \cap \mathcal{O}_\hbar^{g'}\big(u(\zeta)\big)=\emptyset$, and for some $0<r<\hbar$
\begin{equation*}
\|B_u^{g'}(\zeta)\|_{g'} \geq \frac{1}{r}
\end{equation*}
then
\begin{equation*}
\int_{S_r(\zeta)}\|B_u^{g'}\|^2 \geq \hbar;
\end{equation*}
here integration is taken with respect to $u^*\omega'$ where $\omega':=g'\circ (J'\times\mathds{1})$, and $S_r(\zeta)$ is the connected component of $u^{-1}\big(\mathcal{O}_r^{g'}(u(\zeta))\big)$ which contains $\zeta$.
\end{proposition}

A proof of this result can be found in \cite{Fj09b}; it is a modification of the proof of the $\epsilon$-regularity of the second fundamental form of a minimal surface in a Riemannian three-manifold (c.f. \cite{ChSr85}).  We proceed with the proof of Lemma \ref{lem:finiteCurvatureBlowup}.

Next we define an iterative procedure to construct the desired set $\mathcal{S}_2$.  Begin by defining $\mathcal{S}_{2,0}:=\mathcal{S}_1$.  Then either it's the case that there exists a sequence $\zeta_{1,k}\in u_k^{-1}(\mathcal{K}_0)$ such that $\limsup \|\two_{u_k}^{g_k}(\zeta_{1,k})\|=\infty$ and $\dist_g\big(\mathcal{S}_{2,0},\zeta_{1,k})\geq \epsilon$ for some $\epsilon >0$, or else we define $\mathcal{S}_2:=\mathcal{S}_{2,0}$ and we are done; we suppose the former.  In this case we pass to a subsequence so that $\|\two_{u_k}^{g_k}(\zeta_{1,k})\|\to\infty$,  and $u_k(\zeta_{1,k})$ converges to a point $p_1\in\mathcal{K}_0$, and we define the finite set $\mathcal{S}_{2,1}:= \mathcal{S}_{2,0}\cup \{p_1\}$. Again either it's the case that there exists a sequence $\zeta_{2,k}\in u_k^{-1}(\mathcal{K}_0)$ such that $\limsup \|\two_{u_k}^{g_k}(\zeta_{2,k})\|=\infty$ and $\dist_g\big(\mathcal{S}_{2,1},\zeta_{2,k})\geq \epsilon$ for some $\epsilon >0$, or else we define $\mathcal{S}_2:=\mathcal{S}_{2,1}$ and we are done; we again suppose the former and pass to a further subsequence so that $\|\two_{u_k}^{g_k}(\zeta_{2,k})\|\to\infty$,  $u_k(\zeta_{2,k})$ converges to a point $p_2\in\mathcal{K}_0$, and we define $\mathcal{S}_{2,2}:= \mathcal{S}_{2,1}\cup \{p_2\}$. We now iterate this procedure to construct a collection of sets: $\mathcal{S}_{2,0}\subset \mathcal{S}_{2,1}\subset \cdots$.

We now claim that this process must terminate after a finite number of iterations.  Indeed, fix $\epsilon>0$ such that $\mathcal{O}_{2\epsilon}^g(\mathcal{K}_0) \subset \Int(\widetilde{\mathcal{K}})$, and let $\hbar>0$ be the constant guaranteed by Proposition \ref{prop:CurvatureThreshold} and is associated to $(M,J,g)$ and $\epsilon$; also fix $n_0\in\mathbb{N}$, and suppose that $k$ is sufficiently large so that for some $\delta\in (0,\hbar)$ the following conditions hold
\begin{enumerate}
\item $\dist_g(\mathcal{S}_1,\mathcal{S}_{2,n_0}\setminus\mathcal{S}_{1})\geq \delta$,
\item $\dist_{g}\big(u_k(\zeta_{i,k}),u_k(\zeta_{j,k})\big)\geq \delta$ for all $i,j\in\{1,\ldots,n_0\}$ for which $i\neq j$,
\item for all $i\in \{1,\ldots,n_0\}$ we have $\|B_{u_k}^{g_k}(\zeta_{i,k})\|\geq 2/\epsilon$.
\end{enumerate}
Then by Proposition \ref{prop:CurvatureThreshold}, it follows that
\begin{equation*}
n_0 \hbar\leq \sum_{i=1}^{n_0}\int_{S_{\epsilon/2}^{g_k}(\zeta_{i,k})}\|\two_{u_{k}}^{g_k}(\zeta_{i,k})\|_{g_k}^2\leq \int_{S_k}\|\two_{u_k}^{g_k}\|^2 \leq C_{Total},
\end{equation*}
and thus $n_0$ is bounded.  This completes the proof of Lemma \ref{lem:finiteCurvatureBlowup}.
\end{proof}

We now continue with the proof of Theorem \ref{thm:Main1}.  As a consequence of Lemma \ref{lem:finiteCurvatureBlowup}, it follows that after passing to a subsequence, there exist compact sets $\mathcal{K}_i\subset M$ for $i=0,\ldots,4$ such that  $\mathcal{K}_{i+1}\subset\Int(\mathcal{K}_i)$, $\mathcal{K}_0\subset \Int(\widetilde{\mathcal{K}})$, $\mathcal{K}\subset\Int(\mathcal{K}_4)$, and such that
\begin{equation*}
\sup_{\zeta\in u_k^{-1}(\mathcal{K}_1\setminus\mathcal{K}_4)}\|\two_{u_k}^{g_k}(\zeta)\|_{g_k}\leq C <\infty,
\end{equation*}
for all $k$.  In this case we define $\widehat{M}:=\Int(\mathcal{K}_1)\setminus \mathcal{K}_4$ and $\widehat{\mathcal{K}}:= \mathcal{K}_2\setminus \Int(\mathcal{K}_3)$, and $\hat{\mathbf{u}}_k:=(u_k,\widehat{S}_k,j_k,J_k)$ where $\widehat{S}_k:=u_k^{-1}(\widehat{M})$.  Observe that $\widehat{M}$, $\widehat{\mathcal{K}}$, and the $\hat{\mathbf{u}}_k$ satisfy the hypotheses of Proposition \ref{prop:compactnessFromLInfinityBounds}, and thus after passing to a further subsequence, there exist a compact manifold $\widehat{S}$ with boundary, and there exist maps $\hat{\psi}_k:\widehat{S}\to \widehat{S}_k$ which are diffeomorphic with their images and satisfy $u_k\big(\widehat{S}_k\setminus \hat{\psi}_k(\widehat{S})\big)\subset \widehat{M}\setminus \widehat{\mathcal{K}}$, and additionally the maps $\hat{u}_k\circ\hat{\psi}_k:\widehat{S}\to \widehat{M}$ converge in $C^\infty$ to an immersion.  Consequently, we may define the set of boundary circles
\begin{equation*}
\Gamma_{-}:=(u_k\circ\hat{\psi}_k)^{-1}(\mathcal{K}_3)\cap \partial \widehat{S} \qquad\text{and}\qquad\Gamma_{+}:= (\hat{u}_k\circ\hat{\psi}_k)^{-1}\big(M\setminus \Int(\mathcal{K}_2)\big)\cap \partial \widehat{S}
\end{equation*}
so that $\Gamma_{-}\cap\Gamma_{+}=\emptyset$ and $\partial{\widehat{S}}=\Gamma_{-}\cup\Gamma_{+}$.  We can also define
\begin{equation*}
\widetilde{S}_k:= u_k^{-1}\big(\Int(\mathcal{K}_2)\big) \cup \psi_k(\widehat{S}),
\end{equation*}
so that $\partial \widetilde{S}_k = \psi_k(\Gamma_+)$, and we have $\psi_k:\widehat{S}\to\widetilde{S}_k$.  Observe that by construction we have
\begin{equation}\label{eq:trimmedCurveProperty1}
u_k(S_k\setminus \widetilde{S}_k)\subset M \setminus \mathcal{K}_2\qquad\text{and}\qquad u_k(\widetilde{S}_k) \subset \mathcal{K}_0
\end{equation}
for all $k$; also recall that $\mathcal{K}\subset\Int(\mathcal{K}_2)$ and $\mathcal{K}_0\subset\Int(\widetilde{\mathcal{K}})$. Next we observe that the number of boundary components of $\widetilde{S}_k$ is equal to the number of connected components of $\Gamma_+$ which is independent of $k$; furthermore since $\Genus(S_k)\leq C_G$ it follows that $\Genus(\widetilde{S}_k)\leq C_G$. Also note that the number of connected components of $\widetilde{S}_k$ must also be bounded; this follows from  monotonicity of area\footnote{See Proposition \ref{prop:monotonicityOfArea}.}, which guarantees that the image of each closed connected component of $\widetilde{S}_k$ captures a threshold amount of area.  As a consequence of these these facts, it follows that after passing to a further subsequence the $\widetilde{S}_k$ are all diffeomorphic; we denote these diffeomorphisms $\varphi_k:\widetilde{S}\to \widetilde{S}_k$.  Thus we define
\begin{equation*}
\tilde{\mathbf{u}}_k:=(\tilde{u}_k,\widetilde{S},\tilde{j}_k,J_k)=(u_k\circ \varphi_k,\widetilde{S},\varphi_k^*j_k,J_k),
\end{equation*}
and observe that by construction these $J_k$-curves have uniformly bounded area, and their images are contained in $\mathcal{K}_0\subset\Int(M)$.  We would like to claim that a subsequence converges in a Gromov sense, however some care must be taken near $\partial \widetilde{S}$, a matter to which we now attend.

Let $\mathcal{A}:=\cup_{i=1}^{n_0}\mathcal{A}_i$ be the union of pair-wise disjoint annular neighborhoods of $\Gamma_+\subset \widehat{S}$, and let $\theta_i:\mathcal{A}_i\to \mathbb{S}^1\times [0,1)$ be diffeomorphisms.  Next define the diffeomorphisms
\begin{equation*}
\psi_{i,k}:=\varphi_k^{-1}\circ\hat{\psi}_k\circ \theta_i^{-1}: \mathbb{S}^1\times [0,1)\to \psi_{i,k}\big(\mathbb{S}^1\times [0,1)\big) \subset\widetilde{S},
\end{equation*}
and observe that the $\psi_{i,k}$ satisfy the following.
\begin{enumerate}
\item For each fixed $k$, the images of the maps $\psi_{i,k}$ are pairwise disjoint.
\item Each $\psi_{i,k}$ is a diffeomorphism with its image.
\item For each fixed $k$ we have $\partial \widetilde{S} = \bigcup_{i=1}^{n_0} \psi_{i,k}(\mathbb{S}^{1}\times\{0\})$
\item For each fixed $i$ the maps $\tilde{u}_k\circ \psi_{i,k}$ converge in $C^\infty$ to an immersion.
\end{enumerate}
By construction, for all $k$ remaining in our subsequence we have $u_k\big(S_k\setminus \varphi_k(\widetilde{S})\big) \subset M\setminus \mathcal{K}_2$ with $\mathcal{K}_2\subset\Int(\widetilde{\mathcal{K}})$, and by Proposition \ref{prop:GromovWithNiceBdry} below, it follows that a subsequence of the $\tilde{\mathbf{u}}_k$ converge in a Gromov sense.  This completes the proof of Theorem \ref{thm:Main1}.

\end{proof}

\begin{proposition}\label{prop:GromovWithNiceBdry}
Let $(M,J,g)$ be a compact almost Hermitian manifold with boundary, and let $(J_k,g_k)$ be a sequence of almost Hermitian structures which converge in $C^\infty$ to $(J,g)$.  Suppose $\mathbf{u}_k=(u_k,S,j_k,J_k)$ is a sequence of compact $J_k$-curves which satisfy the following conditions.
\begin{enumerate}
\item $\Area_{u_k^*g_k}(S) \leq C < \infty$.
\item $\dist_g\big(u_k(S),\partial M\big) \geq \delta>0$.
\item For $i=1,\ldots,n$ (where $n$ is the number of connected connected components of $\partial S$) and for all $k\in \mathbb{N}$, there exist maps $\psi_{i,k}:\mathbb{S}^1\times[0,\epsilon)\to S$ with the following properties.
\begin{enumerate}
\item For each fixed $k$, the images of the maps $\psi_{i,k}$ are pairwise disjoint.
\item Each $\psi_{i,k}$ is a diffeomorphism with its image.
\item For each fixed $k$ we have $\partial S = \bigcup_{i=1}^n \psi_{i,k}(\mathbb{S}^1\times\{0\})$
\item For each fixed $i$ the maps $u_k\circ \psi_{i,k}$ converge in $C^\infty$ to an immersion.
\end{enumerate}
\end{enumerate}
Then a subsequence of the $\mathbf{u}_k$ converges in a Gromov sense.  Furthermore, the sequence of marked points added to the $(S,j_k)$ in order to obtain Gromov convergence can be chosen in such a way so that the lengths  of the connected components of $\partial S$ (with respect to the poincar\'{e} metrics) are uniformly bounded away from zero and infinity.  Moreover, each component of the limit curve with non-empty boundary is non-constant, and the map is an immersion in a neighborhood of the boundary.
\end{proposition}

A proof of this result can be found in Section \ref{sec:GromovWithNiceBdry}.

\subsection{Compactness without curvature bounds}

In this section we prove Theorem \ref{thm:Main2}, as stated in the beginning of Section \ref{sec:Main}.
Before providing the proof of this result, we assume its validity for the moment and state an immediate corollary (stated in the introduction as Theorem A).

\begin{corollary}\label{cor:Main}
Let $(M,J,g)$ be a compact almost Hermitian manifold with boundary.  Let $(J_k,g_k)$ be a sequence of almost Hermitian structures which converge to $(J,g)$ in $C^\infty(M)$, and let $(u_k,S_k,j_k,J_k)$ be a sequence of compact $J_k$-curves (possibly disconnected, but having no constant components) satisfying the following:
\begin{enumerate}
\item $u_k:\partial S_k \to \partial M$,
\item $\Area_{u_k^*g_k}(S_k)\leq C_{A}$,
\item $\Genus(S_k)\leq C_{G}$.
\end{enumerate}
Then there exists a subsequence (still denoted with subscripts $k$) of the $\mathbf{u}_k$, an $\epsilon>0$, and an open dense set $\mathcal{I}\subset [0,\epsilon)$ with the following significance.  For each $\delta \in \mathcal{I}$, define $\tilde{S}_k^\delta:=\{\zeta\in S_k: \dist_g\big(u_k(\zeta),\partial M)\geq \delta\}$; then the $J_k$-curves $(u_k,\tilde{S}_k^\delta,j_k,J_k)$ converge in a Gromov sense.
\end{corollary}

\begin{proof}
Begin by defining a function $f:M\to\R$ by $f(p):=\dist_g(p,\partial M)$, and define the sets $M^\delta:= f^{-1}([\delta,\infty))$.  Observe that by construction $M=M^0$, $M^{\delta_2}\subset \Int(M^{\delta_1})$ whenever $\delta_2 > \delta_1$, and for all sufficiently small $\delta>0$ the sets $M^\delta$ are compact regions, and $\|df\|$ is uniformly bounded away from zero near $\partial M$.  We then apply Theorem \ref{thm:Main2} to this sequence with $\mathcal{K}=M^1$, to obtain a subsequence.  Apply Theorem \ref{thm:Main2} to this subsequence with $\mathcal{K}=M^{1/2}$ to obtain a further subsequence. We iterate this procedure with $\mathcal{K}=M^{1/\ell}$ and $\ell\in \mathbb{N}$, and pass to further and further subsequences.  Taking a diagonal subsequence we are left with a subsequence of $J_k$-curves which robustly $\mathcal{K}$-converge in a Gromov sense for each $\mathcal{K}\subset \Int(M)$.  The regular values of $f$ composed with the limit curves are an open dense set $\mathcal{I}\subset (0,\delta)$ for some sufficiently small $\delta>0$.  The corollary is then immediate.
\end{proof}

We proceed with the proof of Theorem \ref{thm:Main2} momentarily, but first we state a result upon which the proof heavily relies.

\begin{proposition}[a priori total curvature bounds]\label{prop:Main1}
Let $(M,J,g)$, $(J_k,g_k)$, $\mathcal{K}$, and $\mathbf{u}_k$ be as in the hypotheses of Theorem \ref{thm:Main2}. Then for each compact set $\widetilde{\mathcal{K}}\subset \Int(M)$ for which $\mathcal{K}\subset \Int(\widetilde{\mathcal{K}})$ and for which the $u_k$ are robustly $\widetilde{\mathcal{K}}$-proper,  there exist positive constants $C_{\mathcal{Z}}$ and $C_{total}$ with the following significance. For $\widetilde{S}_k:=u_k^{-1}\big(\Int(\widetilde{\mathcal{K}})\big)$, the $J_k$-curves defined by $\tilde{\mathbf{u}}_k:=(u_k,\widetilde{S}_k,j_k,J_k)$ are robustly $\mathcal{K}$-proper, have uniformly bounded area and genus, and satisfy
\begin{enumerate}
\item $\#\mathcal{Z}_{u_k} \leq C_{\mathcal{Z}}<\infty$,
\item $\int_{\widetilde{S}_k}\|\two_{u_k}^{g_k}\|^2 \leq C_{Total} < \infty$.
\end{enumerate}
\end{proposition}

Postponing the proof of Proposition \ref{prop:Main1} for the moment, we now use it to prove Theorem \ref{thm:Main2}.

\begin{proof}[Proof of Theorem \ref{thm:Main2}]  We begin by applying Proposition \ref{prop:Main1} to obtain the compact set $\widetilde{\mathcal{K}}$ and associated $J_k$-curves $\tilde{\mathbf{u}}_k$.  However these curves satisfy the hypotheses of Theorem \ref{thm:Main1}, and so a subsequence robustly $\mathcal{K}$-converges in a Gromov sense.
\end{proof}

The proof of Proposition \ref{prop:Main1} relies on two main technical results, which we now state.

\begin{proposition}[Desingularization]\label{prop:Desingularization}
Consider $(M,J,g)$ an almost Hermitian manifold, and a compact generally immersed
$J$-curve $\mathbf{u}=(u,S,j,J)$ with immersed boundary. Then for each
$\epsilon>0$, there exists an $0<\epsilon_0<\epsilon$ and an immersion $\hat{u}:S\to
M$ such that the following properties hold.

\begin{enumerate}[({D}1)]
  \item The sets $\mathcal{B}_{\epsilon_0}(z):=\{\zeta\in S:\dist_{u^*g}(z,\zeta)<\epsilon_0\}$ for $z\in \mathcal{Z}_u:=\{z\in S: T_z u =0\}$ are pairwise disjoint. Also $\hat{u}(\zeta)=u(\zeta)$ whenever
  $\zeta\in S\setminus\big(\bigcup_{z\in \mathcal{Z}_u} \mathcal{B}_{\epsilon_0}(z)\big)$. \label{en.D1}
  \item $\sup_{\zeta\in S}\dist_g\big(u(\zeta),\hat{u}(\zeta)\big)\leq
\epsilon_0$\label{en.D2}
  \item For every vector $X$ tangent to the image of $\hat{u}$ we have
\begin{equation*}\label{eq:DesingJProperty}
\|(J X)^\bot\|_g\leq \epsilon_0\|X\|_g,
\end{equation*}
where the map $Y\mapsto Y^\bot$ is the $g$-orthogonal projection to the normal bundle over the immersion $\hat{u}:S\to M$.\label{en.D3}
  \item For any open set $\mathcal{U}\subset S$,
\begin{equation*}
\big|\Area_{u^*g}(\mathcal{U})-\Area_{\hat{u}^*g}(\mathcal{U})\big|\leq \epsilon_0
\end{equation*}\label{en.D4}
  \item The following point-wise estimate holds
  \begin{equation*}
    \sup_{\zeta\in S} K_{\hat{u}^*g}(\zeta)\leq \epsilon_0+\sup_{q\in M}|K_{sec}(q)| + \sup_{\substack{q\in M\\ J \mathcal{P}_q = \mathcal{P}_q}}{\textstyle \frac{1}{2}}\|\tr_{\mathcal{P}_q} J\nabla J\|_g^2 \label{L:ImmersedApproximation5A}
  \end{equation*}\label{eq:DesingCurvatureEstimate}
  where $K_{\hat{u}^*g}(\zeta)$ is the Gaussian curvature of $S$ at the point $\zeta$ with respect to the metric $\hat{u}^*g$,
  $|K_{sec}(q)|$ is defined as in (\ref{eq:defOfAbsKsec}), and $\tr_{\mathcal{P}_q} J\nabla J$ is the trace of the $(1,2)$-tensor $J\nabla J$ along the $J$-invariant plane $\mathcal{P}_q = \R e \oplus \R J e\subset T_q M$.\label{en.D5}
\item Let $\mathcal{U}\subset S$ be an open set, and define the set
    \begin{equation*}
    \mathcal{U}^{\epsilon_0}:=\{\zeta\in \mathcal{U}: \dist_{u^*g}(\zeta,\partial \mathcal{U}) > \epsilon_0\}.
    \end{equation*}
    Then
\begin{equation*}
-\int_{\mathcal{U}}K_{u^*g}+(1-\epsilon_0)2\pi \sum_{z\in \mathcal{Z}_u\cap\mathcal{U}^{\epsilon_0}}\order(z)\leq
-\int_{\mathcal{U}}K_{\hat{u}^*g}+\epsilon_0
\end{equation*}
where $\order(z)$ is given as in Definition \ref{def:ord}.\label{en.D6}
\end{enumerate}
\end{proposition}
The proof of Proposition \ref{prop:Desingularization} is given Section \ref{sec:Desingularization}.  We take a moment to summarize the results of said proposition.  Roughly it guarantees that any $J$-curve can be perturbed a $C^0$ small amount only near its critical points, in such a way that it becomes immersed, and the resulting tangent planes are $C^0$-close to being $J$-invariant, the resulting area changes by only a small amount, the Gaussian curvature is uniformly bounded from above, and each original critical point is locally traded for a threshold amount of total curvature.

We now continue with the proof of Proposition \ref{prop:Main1}.  To that end, we turn our attention towards showing that it is not possible for too many nodes to develop, and that away from a finite set of points in $\widetilde{\mathcal{K}}$, the integral of the Gaussian curvature is bounded from below.  We make this precise with Proposition \ref{prop:nodesAndCurvature} below.

\begin{proposition}\label{prop:nodesAndCurvature}
Let $(M,J,g)$ be an almost Hermitian manifold possibly with boundary, suppose $(J_k,g_k)\to (J,g)$ in $C^\infty(M)$, and let  $\mathcal{K}\subset \subset M$ be a compact set.  Suppose further that $\mathbf{u}_k:=(u_k,S_k,j_k, J_k)$ is a sequence of compact generally immersed $J_k$-curves which are robustly $\mathcal{K}$-proper, and satisfy
\begin{enumerate}
\item $\Area_{u_k^*g_k}(S_k)\leq C_A < \infty$
\item $\Genus(S_k)\leq C_G < \infty$
\item $\mathcal{Z}_{u_k}\cap \partial S_k = \emptyset$.
\end{enumerate}
Furthermore, for a sequence of positive numbers $\epsilon_k\to 0$, let $v_k$ be the immersed approximations associated to $(u_k,\epsilon_k)$ yielded by Proposition \ref{prop:Desingularization}.  Then after passing to a subsequence, there exists a finite set $\mathcal{S}=\{\sigma_1,\ldots,\sigma_{n_0}\}\subset M$ and  $\delta_0>0$ with the following significance.  For each $0<\delta<\delta_0$, there exists $\epsilon>0$ and $k_0\in \mathbb{N}$ such that
\begin{equation*}
\text{if}\quad k\geq k_0\;\;\text{and}\;\;v_k(\zeta)\in \mathcal{K}\setminus \mathcal{O}_\delta^{g_k}(\mathcal{S}) \qquad\text{then}\qquad \inj_{S_k}^{v_k^*g_k}(\zeta) > \epsilon;
\end{equation*}
here $\inj_{S_k}^{v_k^*g_k}(\zeta)$ is the injectivity radius of $S_k$ at the point $\zeta$ computed with respect to the metric $v_k^*g_k$.  Furthermore, for each $0<\delta<\delta_0/2$ there exists a constant $C>0$ such that for all sufficiently large $k$ in the subsequence we have
\begin{equation*}
-\int_{\widehat{S}_k^{{\delta}}}K_{v_k^*g_k}\leq C,\qquad\text{where}\qquad \widehat{S}_k^{{\delta}}:=v_k^{-1}\big(\Int(\mathcal{K})\setminus \overline{\mathcal{O}_\delta^{g_k}(\mathcal{S})}\big);
\end{equation*}
here $K_{v_k^*g_k}:\widehat{S}_k^\delta\to \mathbb{R}$ is the Gaussian curvature associated to the metric $v_k^*g_k$.
\end{proposition}

The proof of Proposition \ref{prop:nodesAndCurvature} can be found Section \ref{sec:nodesAndCurvature}. Roughly, the idea is to show that if there were many locations in which the injectivity radius were very small, then one could remove many small loops and disconnect the $J_k$-curves into many connected components each of which has a threshold amount of area, which would yield a contradiction. Then one sees that in the absence of arbitrarily small injectivity radii, exceedingly negative Gaussian curvature results in exceedingly large area, which also yields a contradiction.  At present, we now provide the proof of Proposition \ref{prop:Main1}

\begin{proof}[Proof of Proposition \ref{prop:Main1}]
Suppose not.  Then there exists a compact set $\widetilde{\mathcal{K}}\subset \Int(M)$ with $\mathcal{K}\subset\Int(\widetilde{\mathcal{K}})$ for which the ${u}_k$ are robustly $\widetilde{\mathcal{K}}$-proper, and either the total curvature or the number of critical points is unbounded on $\widetilde{\mathcal{K}}$. Since the ${u}_k$ are robustly $\widetilde{\mathcal{K}}$-proper it follows that there exist compact regions $\widetilde{S}_k\subset S_k$ with the property that $\mathcal{Z}_{u_k}\cap  \partial \widetilde{S}_k =\emptyset$, and $u_k(S_k\setminus \widetilde{S}_k)\subset M\setminus \mathcal{K}_0$ for some compact set $\mathcal{K}_0\subset \Int(M)$ for which $\widetilde{\mathcal{K}}\subset \Int(\mathcal{K}_0)$.  By assumption, the restricted $J_k$-curves $u_k:\widetilde{S}_k\to M$ are again robustly $\widetilde{\mathcal{K}}$-proper, and this sequence of curves again has either an unbounded number of critical points or else unbounded total curvature on $\widetilde{\mathcal{K}}$.  Rather than expending notation to keep track of the $u_k$ restricted to the $\widetilde{S}_k$, we will (without much loss of generality) assume $\widetilde{S}_k\equiv S_k$ for all $k$.

Next we consider a sequence of positive numbers $\epsilon_k\to 0$ as $k\to\infty$, and consider the immersed approximations $v_k$ associated to $(u_k,\epsilon_k)$ and yielded by Proposition \ref{prop:Desingularization}.  We then apply Proposition \ref{prop:nodesAndCurvature} for some auxiliary compact set $\check{\mathcal{K}}\subset \Int(M)$ for which $\widetilde{\mathcal{K}}\subset \Int(\check{\mathcal{K}})$ and for which the $u_k$ are again robustly $\check{\mathcal{K}}$-proper.  Consequently, after passing to a subsequence there exists a finite set  $\mathcal{S}=\{\sigma_1,\ldots,\sigma_{n_0}\}\subset M$ with the properties guaranteed by that proposition.  As a further consequence of Proposition \ref{prop:nodesAndCurvature},
for each sufficiently small $\delta>0$ there exists a constant $C>0$ such that for all sufficiently large $k$ in our subsequence, we have
\begin{equation*}
-\int_{\widehat{S}_k^\delta}K_{v_k^*g_k}\leq C,\qquad\text{where}\qquad \widehat{S}_k^\delta:=v_k^{-1}\big(\mathcal{O}_\delta^{g_k}(\widetilde{\mathcal{K}})\setminus \overline{\mathcal{O}_\delta^{g_k}(\mathcal{S})}\big).
\end{equation*}
Fix $\delta>0$ sufficiently small so that the sets $\mathcal{O}_{3\delta}^{g_k}(\sigma_i)$ are pair-wise disjoint for all sufficiently large $k$, and $\mathcal{O}_{4\delta}^{g}(\widetilde{\mathcal{K}})\subset \Int(\check{\mathcal{K}})$.
Next define
\begin{equation*}
\widehat{M}':=\mathcal{O}_{3\delta}^g(\widetilde{\mathcal{K}})\setminus\overline{\mathcal{O}_{\delta}^{g}(\mathcal{S})},
\;\;
\widehat{M}:=\mathcal{O}_{2\delta}^g(\widetilde{\mathcal{K}})\setminus \overline{\mathcal{O}_{2\delta}^{g}(\mathcal{S})}
\;\;\text{and}\;\;
\widehat{\mathcal{K}}:=\overline{\mathcal{O}_{\delta}^g(\widetilde{\mathcal{K}})}\setminus \mathcal{O}_{3\delta}^{g}(\mathcal{S});
\end{equation*}
then for all sufficiently large $k$ we have
\begin{equation}\label{eq:mainProp1}
-\int_{v_k^{-1}(\widehat{M}')}K_{v_k^*g_k}\leq C.
\end{equation}
However, recall that the $\Area_{u_k^*g_k}(S_k)$ are uniformly bounded; by property (D\ref{en.D4}) of Proposition \ref{prop:Main1} the $\Area_{v_k^*g_k}(S_k)$ are also uniformly bounded. Furthermore by property (D\ref{en.D5}) of Proposition \ref{prop:Desingularization}, it follows that the Gaussian curvatures $K_{v_k^*g_k}$ uniformly point-wise bounded from above.  Consequently (\ref{eq:mainProp1}) and property (D\ref{en.D6}) of Proposition \ref{prop:Desingularization} allow us to conclude that for all sufficiently large $k$ we have
\begin{equation}\label{eq:mainProp2}
\#\big(\mathcal{Z}_{u_k}\cap u_k^{-1}(\widehat{M})\big) \leq C'\qquad\text{and}\qquad \int_{u_k^{-1}(\widehat{M})}K_{u_k^*g_k}\leq C',
\end{equation}
where $C'$ depends on $C$, $g$, $J$, and the uniform area and genus bounds on the $\mathbf{u}_k$.
Next we recall the Gauss equations for $J$-curves:
\begin{equation*}
K_{sec}+\textstyle{\frac{1}{2}}\|\tr_S Q\|^2=K_g + {\textstyle \frac{1}{2}}\|\two\|^2,
\end{equation*}
where $K_{sec}$ is the sectional curvature and $K_g$ is the Gaussian curvature.  Integrating these equations then yields
\begin{align*}
{\textstyle\frac{1}{2}}\int_{u_k^{-1}(\widehat{M})}\|\two_{u_k}^{g_k}\|^2 &= \int_{u_k^{-1}(\widehat{M})} K_{sec} - \int_{u_k^{-1}(\widehat{M})}K_{u_k^*g_k}+\textstyle{\frac{1}{2}}\int_{u_k^{-1}(\widehat{M})}\|\tr_S Q\|^2\\
&\leq \big(\|K_{sec}\|_{L^\infty}+2\|\nabla J\|_{L^\infty}^2\big) \Area_{u_k^*g_k}\big(u_k^{-1}(\widehat{M})\big) + C',
\end{align*}
which is uniformly bounded.  Combining this fact with the left-most statement of (\ref{eq:mainProp2}), then shows that the $J_k$-curves $\widehat{u}_k:=(u_k,u_k^{-1}(\widehat{M}),j_k,J_k)$ are robustly $\widehat{\mathcal{K}}$-proper, and satisfy the hypotheses of Theorem \ref{thm:Main1}.  We conclude that after passing to a further subsequence, there exists a compact manifold $\overline{S}$ with boundary, and diffeomorphisms $\phi_{k}:\overline{S}\to\phi_k(\overline{S})\subset u_k^{-1}(\widehat{M})\subset S_k$ with the following properties
\begin{enumerate}
\item $u_k\big(S_k \setminus \phi_k(\overline{S})  \big) \subset  \mathcal{O}_{3\delta}^g(\mathcal{S}) \cup \big(M\setminus \widetilde{\mathcal{K}}\big)$,
\item $u_k\circ\phi_k\big|_{\partial \overline{S}}$ converges in $C^\infty$ to an immersion,
\item $\int_{\overline{S}}\|\two_{u_k\circ\phi_k}^{g_k}\|_{g_k}^2$ is uniformly bounded,
\item $\#\mathcal{Z}_{u_k\circ \phi_k}$ is uniformly bounded.
\end{enumerate}
The first two properties are consequences of Theorem \ref{thm:Main1}, and the last two are by construction.  We now take a moment to recall our method of proof for Proposition \ref{prop:Main1}: we are assuming that either the number of critical points or total curvature in $\widetilde{\mathcal{K}}$ is unbounded.  However, after passing to a subsequence we see as a consequence of points 1, 3, and 4 above above, it is only possible for these quantities to blow up in the set $\mathcal{O}_{3\delta}^g(\mathcal{S})=\cup_{i=1}^{n_0}\mathcal{O}_{3\delta}^{g}(\sigma_i)$.  This leads us to define the following
\begin{equation*}
\check{S}_{i,k}:= \big(S_k\setminus\psi_k(\overline{S}) \big)\cap u_k^{-1}\big(\mathcal{O}_{3\delta}^g(\sigma_i)\big)\qquad\text{and}\qquad \check{\mathbf{u}}_{i,k}=(u_k,\check{S}_{i,k},j_k,J_k).
\end{equation*}
We note that these curves have uniformly bounded area and genus, and each has image in $\mathcal{O}_{3\delta}^g(\sigma_i)$.  Also note that for each $i=1,\ldots,n_0$ we have $\partial \check{S}_{i,k}\subset \psi_{k}(\partial \overline{S})$, and by point 2 above, $u_k\circ \phi_k\big|_{\partial \overline{S}}$ converge in $C^\infty$ to an immersion.  Consequently the geodesic curvature $\kappa_{u_k^*g_k}$ of $\partial \check{S}_{i,k} \subset \check{S}_{i,k}$ is uniformly bounded.  Thus we let $v_{i,k}:=u_{k}^{\epsilon_k}$ be the approximations of $u_k\big|_{\check{S}_{i,k}}$ guaranteed by Proposition \ref{prop:Desingularization}.  Arguing as before (i.e. making use of properties (D\ref{en.D4}), (D\ref{en.D5}), and (D\ref{en.D6}) of Proposition \ref{prop:Desingularization}) we see that to complete our proof by contradiction, it is sufficient to show that
\begin{equation}\label{eq:mainProp3}
-\int_{\check{S}_{i,k}}K_{v_{i,k}^*g_k} \leq C'' <\infty
\end{equation}
for all sufficiently large $k$.  However, at this point we invoke the Gauss-Bonnet theorem, and find that
\begin{equation*}
-\int_{\check{S}_{i,k}}K_{v_{i,k}^*g_k}= -\chi(\check{S}_{i,k}) + \int_{\partial \check{S}_{i,k}} \kappa_{v_{i,k}^*g_k}.
\end{equation*}
We have already argued that the last term on the right hand side is uniformly bounded, so it is sufficient to show that $-\chi(\check{S}_{i,k})$ is uniformly bounded.  However recall that $-\chi(\check{S}_{i,k}) = -2 + 2\Genus(\check{S}_{i,k}) + b$, where $b$ is the number of boundary components of $\check{S}_{i,k}$.  However $\Genus(\check{S}_{i,k})\leq \Genus(S_k) \leq C_G$, and $\partial \check{S}_{i,k} \subset \psi_k(\partial \overline{S})$ and $\overline{S}$ has finitely many boundary components (and no $k$-dependence), from which it follows that indeed $-\chi(\check{S}_{i,k})$ is uniformly bounded.  This shows inequality  (\ref{eq:mainProp3}) holds, which in turn provides the desired contradiction, which completes the proof of Proposition \ref{prop:Main1}.
\end{proof}

\section{Proofs}\label{sec:proofs}

Here we prove some of the more technical results from the previous sections.

\subsection{Proof of Proposition \ref{prop:KlocImpliesK}}\label{sec:KlocImpliesK}
\begin{proof}We begin by fixing some notation.  Let $\phi_{i,k}:\mathcal{D}_{r_0}^i\to S_k$ be the maps guaranteed by the definition of $\mathcal{K}_{loc}$-convergence, and define the maps
\begin{equation*}
\tilde{u}_{i,k}:=u_{k}\circ\phi_{i,k}:\mathcal{D}_{r_0}^i\to M\qquad\text{and}\qquad \tilde{u}_{i,\infty}:=\lim_{k\to\infty} \tilde{u}_{i,k}:\mathcal{D}_{r_0}^i\to M.
\end{equation*}
Note that $\mathcal{D}_{r_0}=\mathcal{D}_{r_0}^i$; in this case the superscript $i$ simply enumerates the domains the the maps $\phi_{i,k}$. Note that since the limit is immersed and since the sequence is uniformly and robustly covered, it is possible to construct a refined uniform and robust cover which has the additional property that each $\tilde{u}_{i,k}$ and $\tilde{u}_{i,\infty}$ is an embedding.

Next we fix a smooth auxiliary Riemannian metric $g$ on $M$ which has the property
\begin{equation}\label{eq:compareMetrics}
\|X\|_{\bar{g}}\leq \|X\|_{\tilde{u}_{i,k}^*{g}}
\end{equation}
for all $X=X_q\in T_q \mathcal{D}_{r_0}$;  here $\bar{g}$ is the canonical Euclidean metric on $\mathcal{D}_{r_0}$.  Recall our notation that if $(W,\tilde{g})$ is a Riemannian manifold, $p\in W$, and $\epsilon>0$, then
\begin{equation*}
\mathcal{O}_\epsilon^{\tilde{g}}(p):=\{q\in W: \dist_{\tilde{g}}(p,q)<\epsilon\}.
\end{equation*}
The proof of Proposition \ref{prop:KlocImpliesK} is now split into three main steps: constructing the auxiliary manifold $\widetilde{S}$, constructing (almost) reparameterizations $\psi_k:\widetilde{S}\to S_k$, and then showing that these maps have the desired properties.  We approach these steps in order, and begin with a rather technical result.

\begin{lemma}\label{lem:constructSurface}
Fix $r_1\in (0,r_0)$, and suppose $\rho\in \mathcal{D}_{r_1}^i$ and $\rho'\in\mathcal{D}_{r_1}^j$ and
\begin{equation}\label{eq:constructSurface0}
\dist_{u_k^*g}\big(\phi_{i,k}(\rho),\phi_{j,k}(\rho')\big)\to 0.
\end{equation}
Then for each $\delta\in\big(0,(r_0-r_1)/16\big)$ and all sufficiently large $k$, the maps
\begin{equation}\label{eq:constructSurface2}
\phi_{j,k}^{-1}\circ\phi_{i,k}:\mathcal{O}_{2\delta}^{\tilde{u}_{i,k}^*g}(\rho)\to\mathcal{O}_{2\delta}^{\tilde{u}_{j,k}^*g}(\rho')
\end{equation}
are well defined, and they are smooth diffeomorphisms.  Furthermore, with $\rho$ and $\rho'$ as above and for which (\ref{eq:constructSurface0}) holds, and for any points  $\tilde{\rho}\in\mathcal{O}_\delta^{\tilde{u}_{i,\infty}^*g}(\rho)$ and $\tilde{\rho}'\in\mathcal{O}_\delta^{\tilde{u}_{j,\infty}^*g}(\rho')$ with $\tilde{u}_{i,\infty}(\tilde{\rho})=\tilde{u}_{j,\infty}(\tilde{\rho}')$ we have
\begin{equation}\label{eq:constructSurface3}
\dist_{u_k^*g}\big(\phi_{i,k}(\tilde{\rho}),\phi_{j,k}(\tilde{\rho}')\big)\to 0,
\end{equation}
and the maps in (\ref{eq:constructSurface2}) with domains restricted to $\mathcal{O}_{\delta}^{\tilde{u}_{i,\infty}^*g}(\rho)$
converge in $C^\infty$ to the map
\begin{equation}\label{eq:constructSurface4}
\tilde{u}_{j,\infty}^{-1}\circ\tilde{u}_{i,\infty}:\mathcal{O}_{\delta}^{\tilde{u}_{i,\infty}^*g}(\rho)\to\mathcal{O}_{\delta}^{\tilde{u}_{j,\infty}^*g}(\rho').
\end{equation}
\end{lemma}
\begin{proof}
We begin by fixing $k_0\in\mathbb{N}$ so that for all $k\geq k_0$ we have
\begin{equation*}\label{eq:constructSurface5}
\dist_{u_k^*g}\big(\phi_{i,k}(\rho),\phi_{j,k}(\rho')\big)< \delta.
\end{equation*}
Next we note that
\begin{equation}\label{eq:constructSurface6}
\phi_{i,k}:\mathcal{O}_{2\delta}^{\tilde{u}_{i,k}^*g}(\rho)\to \mathcal{O}_{2\delta}^{u_k^*g}\big(\phi_{i,k}(\rho)\big)
\end{equation}
is a diffeomorphism; this follows as a consequence of inequality (\ref{eq:compareMetrics}), namely
\begin{equation*}\label{eq:constructSurface7}
\mathcal{O}_{2\delta}^{\tilde{u}_{i,k}^*g}(\rho)\subset \mathcal{O}_{2\delta}^{\bar{g}}(\rho)\subset \mathcal{D}_{r_0}^i;
\end{equation*}
a similar statement holds with $i$ and $\rho$ replaced with $j$ and $\rho'$ respectively.  Thus to prove the maps in (\ref{eq:constructSurface2}) are smooth diffeomorphisms, it is sufficient to prove that
\begin{equation}\label{eq:constructSurface8}
\phi_{i,k}\big(\mathcal{O}_{2\delta}^{\tilde{u}_{i,k}^*g}(\rho)\big) \subset \phi_{j,k}(\mathcal{D}_{r_0}^j).
\end{equation}
To that end, we fix $\tilde{\rho}\in \mathcal{D}_{r_0}^i$ such that $\dist_{\tilde{u}_{i,k}^*g}(\rho,\tilde{\rho})<2\delta$. It then follows that $\dist_{u_k^*g}\big(\phi_{i,k}(\rho),\phi_{i,k}(\tilde{\rho})\big)<2\delta$. By (\ref{eq:constructSurface0}) and the triangle inequality, it follows that $\dist_{u_k^*g}\big(\phi_{j,k}(\rho'),\phi_{i,k}(\tilde{\rho})\big)<3\delta$, or in other words $\phi_{i,k}(\tilde{\rho})\in\mathcal{O}_{3\delta}^{u_{k}^*g}\big(\phi_{j,k}(\rho')\big)$.
However, since $3\delta<(r_0-r_1)$ it again follows from  (\ref{eq:compareMetrics}) that
\begin{equation}\label{eq:constructSurface9}
\phi_{i,k}(\tilde{\rho})\in \mathcal{O}_{3\delta}^{u_k^*g}\big(\phi_{j,k}(\rho')\big)\subset \phi_{j,k}\big(\mathcal{O}_{3\delta}^{\bar{g}}(\rho')
\big)\subset \phi_{j,k}(\mathcal{D}_{r_0}^j).
\end{equation}
Since $\tilde{\rho}$ was an arbitrary point in $\mathcal{O}_{2\delta}^{\tilde{u}_{i,k}^*g}(\rho)$, we see that we have proved (\ref{eq:constructSurface8}), and thus the maps in (\ref{eq:constructSurface2}) are smooth diffeomorphims.

To prove the next part of the lemma we assume that $\tilde{\rho}\in\mathcal{O}_\delta^{\tilde{u}_{i,\infty}^*g}(\rho)$ and $\tilde{\rho}'\in\mathcal{O}_\delta^{\tilde{u}_{j,\infty}^*g}(\rho')$ with $\tilde{u}_{i,\infty}(\tilde{\rho})=\tilde{u}_{j,\infty}(\tilde{\rho}')$, and we will show that
(\ref{eq:constructSurface3}) holds.  Indeed, since the sequences of maps $\tilde{u}_{i,k}$ and $\tilde{u}_{j,k}$ converge in $C^\infty$, it follows that for all sufficiently large $k$, we have
\begin{equation*}
\dist_{u_k^*g}\big(\phi_{j,k}(\rho'),\phi_{j,k}(\tilde{\rho}')\big)<2\delta,
\end{equation*}
and thus by (\ref{eq:constructSurface2}) we can define
\begin{equation*}
\tilde{\rho}_k':=\phi_{i,k}^{-1}\circ\phi_{j,k}(\tilde{\rho}')\in \mathcal{O}_{2\delta}^{\tilde{u}_{i,k}^*g}(\rho)\subset\overline{\mathcal{D}}_{(r_0+3r_1)/4}^i\subset\mathcal{D}_{r_0}^i.
\end{equation*}
Passing to a subsequence, we may assume the $\tilde{\rho}_k'$ converge to $\tilde{\rho}_\infty'\in \mathcal{D}_{r_0}^i$. We then find
\begin{equation*}
\tilde{u}_{i,\infty}(\tilde{\rho}_\infty')=\lim_{k\to\infty} \tilde{u}_{i,k}(\tilde{\rho}_{k}') =\lim_{k\to\infty} u_k \big(\phi_{j,k}(\tilde{\rho}')\big)=\tilde{u}_{j,\infty}(\tilde{\rho}') = \tilde{u}_{i,\infty}(\tilde{\rho}).
\end{equation*}
However, since $\tilde{u}_{i,\infty}:\mathcal{D}_{r_0}^i\to M$ is an embedding, it follows that $\tilde{\rho}=\tilde{\rho}_\infty'$; consequently $\dist_{\tilde{u}_{i,k}^*g}\big(\tilde{\rho},\phi_{i,k}^{-1}\circ\phi_{j,k}(\tilde{\rho}')\big)\to 0$, and thus (\ref{eq:constructSurface3}) holds.

To prove the last part of the lemma, we observe that $\phi_{j,k}^{-1}\circ \phi_{i,k} = \tilde{u}_{j,k}^{-1}\circ\tilde{u}_{i,k}$ which converges in $C^\infty$ to $\tilde{u}_{j,\infty}^{-1}\circ\tilde{u}_{i,\infty}$.  This completes the proof of Lemma \ref{lem:constructSurface}
\end{proof}

For clarity, now we define
\begin{equation*}
\mathcal{U}_{ij}:=\mathcal{D}_{r_1}^i\cap\tilde{u}_{i,\infty}^{-1}\big(\tilde{u}_{j,\infty}(\mathcal{D}_{r_1}^j) \big)\qquad\text{and}\qquad \overline{\mathcal{U}}_{ij}:=\overline{\mathcal{D}}_{r_1}^i\cap\tilde{u}_{i,\infty}^{-1}\big(\tilde{u}_{j,\infty}(\overline{\mathcal{D}}_{r_1}^j) \big);
\end{equation*}
we note that $\overline{\mathcal{U}}_{ij}$ is closed and contains the closure of $\mathcal{U}_{ij}$ in $\mathcal{D}_{r_0}^i$.  Next for $\delta$ as above, and for each $i,j\in\{1,\ldots,n\}$ for which $\overline{\mathcal{U}}_{ij}\neq\emptyset$, and for $\ell=1,\ldots, m_{ij}$ we let $\rho_{ij\ell}\in \overline{\mathcal{U}}_{ij}$ be points such that
\begin{equation}\label{eq:constructSurfaceA}
\overline{\mathcal{U}}_{ij}\subset \bigcup_{\ell=1}^{m_{ij}}\mathcal{O}_{\delta}^{\tilde{u}_{i,\infty}^*g}(\rho_{ij\ell})\subset\mathcal{D}_{r_0}^i.
\end{equation}
Note that the finiteness of the $\{\rho_{ij1},\rho_{ij2},\ldots\}$ is a consequence of the fact that $\overline{\mathcal{U}}_{ij}\subset \mathcal{D}_{r_0}^i$ is compact and $\bigcup_{\rho\in\overline{\mathcal{U}}_{ij}}\mathcal{O}_\delta^{\widetilde{u}_{i,\infty}^*g}(\rho)$ is an open cover of $\overline{\mathcal{U}}_{ij}$.  Next we define $\rho_{ij\ell}'\in\overline{\mathcal{D}}_{r_1}^j$ to be the unique point for which $\tilde{u}_{j,\infty}(\rho_{ij\ell}')=\tilde{u}_{i,\infty}(\rho_{ij\ell})$.  Now, since the set of points $\{\rho_{ij\ell}\}$ is finite we may pass to a subsequence (still denoted with subscripts $k$) so that for each pair $(\rho_{ij\ell},\rho_{ij\ell}')$ one of the two statements holds:
\begin{enumerate}
\item ${\displaystyle\liminf_{k\to\infty}}\dist_{u_k^*g}\big(\phi_{i,k}(\rho_{ij\ell}),\phi_{j,k}(\rho_{ij\ell}')\big)> 0$
\item ${\displaystyle\lim_{k\to\infty}}\dist_{u_k^*g}\big(\phi_{i,k}(\rho_{ij\ell}),\phi_{j,k}(\rho_{ij\ell}')\big)=0$.
\end{enumerate}
Thus we may define $P_{ij}\subset\{\rho_{ij1},\ldots,\rho_{ijm_{ij}}\}$ to be those points which satisfy the second condition.  For convenience we also define $P_{ij}':= \tilde{u}_{j,\infty}^{-1}\circ\tilde{u}_{i,\infty}(P_{ij})$.
As a consequence of Lemma \ref{lem:constructSurface}, it follows that the sets
\begin{equation*}
\widetilde{\mathcal{U}}_{ij}:=\mathcal{D}_{r_1}^i\cap\bigcup_{\rho\in P_{ij}'}\tilde{u}_{i,\infty}^{-1}\circ\tilde{u}_{j,\infty}\Big(\mathcal{O}_\delta^{\tilde{u}_{j,\infty}^*g}(\rho)\cap \mathcal{D}_{r_1}^j\Big)
\end{equation*}
are open, and the maps denoted by $\tilde{u}_{j,\infty}^{-1}\circ\tilde{u}_{i,\infty}:\widetilde{\mathcal{U}}_{ij}\to\widetilde{\mathcal{U}}_{ji}$ are smooth diffeomorphisms.  We now provide a convenient characterization of the $\widetilde{\mathcal{U}}_{ij}$.

\begin{lemma}[characterization of $\widetilde{\mathcal{U}}_{ij}$]\label{lem:transitionCharacterization}
Having passed to the subsequence as above, we suppose that $\rho\in\mathcal{D}_{r_1}^i$, $\rho'\in\mathcal{D}_{r_1}^j$, and
\begin{equation}\label{eq:transitionCharacterization1}
\liminf_{k\to\infty}\dist_{u_k^*g}\big(\phi_{i,k}(\rho),\phi_{j,k}(\rho')\big) = 0.
\end{equation}
Then $\rho\in\widetilde{\mathcal{U}}_{ij}$, $\rho'\in\widetilde{\mathcal{U}}_{ji}$, and $\tilde{u}_{i,\infty}(\rho)=\tilde{u}_{j,\infty}(\rho')$.  Also, if $\rho\in\widetilde{\mathcal{U}}_{ij}$, $\rho'\in\widetilde{\mathcal{U}}_{ji}$, and $\tilde{u}_{i,\infty}(\rho)=\tilde{u}_{j,\infty}(\rho')$ then
\begin{equation}\label{eq:transitionCharacterization2}
\lim_{k\to\infty}\dist_{u_k^*g}\big(\phi_{i,k}(\rho),\phi_{j,k}(\rho')\big) = 0,
\end{equation}
\end{lemma}
\begin{proof}
As we shall see, this result follows quickly from the definition of the $\widetilde{\mathcal{U}}_{ij}$ and Lemma \ref{lem:constructSurface}.  We begin by noting that
if $\rho$ and $\rho'$ are as in the first part of the lemma, and if (\ref{eq:transitionCharacterization1}) holds, then since the sequences of points $\tilde{u}_{i,k}(\rho)$ and $\tilde{u}_{j,k}(\rho')$ converge, it follows that $\tilde{u}_{i,\infty}(\rho)=\tilde{u}_{j,\infty}(\rho')$, and thus $\rho\in\mathcal{U}_{ij}$ and $\rho'\in\mathcal{U}_{ji}$.  By construction, there exists $\rho_{ij\ell_0}\in\overline{\mathcal{U}}_{ij}$ and $\rho_{ij\ell_0}'\in\overline{\mathcal{U}}_{ji}$ which satisfy the following:
\begin{enumerate}
\item $\tilde{u}_{i,\infty}(\rho_{ij\ell_0})=\tilde{u}_{j,\infty}(\rho_{ij\ell_0}')$
\item $\dist_{\tilde{u}_{i,\infty}^*g}(\rho,\rho_{ij\ell_0})<\delta$
\item either
\begin{enumerate}
\item ${\displaystyle \lim_{k\to\infty} }\dist_{u_k^*g}\big(\phi_{i,k}(\rho_{ij\ell_0}),\phi_{j,k}(\rho_{ij\ell_0}')\big) =0$.
\item ${\displaystyle \liminf_{k\to\infty}} \dist_{u_k^*g}\big(\phi_{i,k}(\rho_{ij\ell_0}),\phi_{j,k}(\rho_{ij\ell_0}')\big)>0$
\end{enumerate}
\end{enumerate}
Since (\ref{eq:transitionCharacterization1}) holds it follow from Lemma \ref{lem:constructSurface} that property (3a) must hold.
It then follows from the definition of the $\widetilde{\mathcal{U}}_{ij}$, that $\rho\in \widetilde{\mathcal{U}}_{ij}$, and thus $\rho'\in\widetilde{\mathcal{U}}_{ji}$.

To prove the second part of the lemma we note that if $\rho\in\widetilde{\mathcal{U}}_{ij}$ and $\rho'\in\widetilde{\mathcal{U}}_{ji}$, with $\tilde{u}_{i,\infty}(\rho)=\tilde{u}_{j,\infty}(\rho')$, then again by construction there exist $\rho_{ij\ell_0}\in\overline{\mathcal{U}}_{ij}$ and $\rho_{ij\ell_0}'\in\overline{\mathcal{U}}_{ji}$ which satisfy properties 1, 2, and 3 above.  Observe that since $\rho\in\widetilde{\mathcal{U}}_{ij}$, it follows from the definition of the $\widetilde{\mathcal{U}}_{ij}$ that it must be the case that property (3a) holds. By Lemma \ref{lem:constructSurface} it then follows that (\ref{eq:transitionCharacterization2}) must also hold.
\end{proof}

Later it will be convenient to have the following corollary at our disposal.

\begin{corollary}\label{cor:reparamConstruction}
For each $r_4\in(0,r_1)$ there exists $k_0\in\mathbb{N}$ and $\epsilon>0$ such that the following holds. If $k\geq k_0$, $\rho\in\mathcal{D}_{r_4}^i$, $\rho'\in\mathcal{D}_{r_4}^j$, $i\neq j$, and
\begin{equation*}
\dist_{u_k^*g}\big(\phi_{i,k}(\rho),\phi_{j,k}(\rho')\big)<\epsilon,
\end{equation*}
then $\rho\in \widetilde{\mathcal{U}}_{ij}$ and $\rho'\in\widetilde{\mathcal{U}}_{ji}$.  Furthermore, for each compact set $\mathcal{V}\subset \widetilde{\mathcal{U}}_{ij}$, and each open set $\mathcal{O}\subset \mathcal{D}_{r}^j$ for which $\tilde{u}_{j,\infty}^{-1}\circ \tilde{u}_{i,\infty}(\mathcal{V})\subset \mathcal{O}$ there exists a $k_0\in \mathbb{N}$ such that $\phi_{i,k}(\mathcal{V})\subset \phi_{j,k}(\mathcal{O})$ for all $k\geq k_0$.
\end{corollary}
\begin{proof}
Suppose not. Then for some $r_4\in(0,r_1)$ (and after possibly passing to a subsequence) there exist some $i\neq j$ and sequences of points $\rho_k\in\mathcal{D}_{r_4}^i$ and $\rho_k'\in\mathcal{D}_{r_4}^j$, such that
\begin{equation}\label{eq:anotherCharacterization}
\dist_{u_k^*g}\big(\phi_{i,k}(\rho_k),\phi_{j,k}(\rho_k')\big)\to 0,
\end{equation}
and $\rho_k\notin\widetilde{\mathcal{U}}_{ij}$. After passing to a further subsequence, we assume that these sequences converge: $\rho_k\to\rho_\infty\in \overline{\mathcal{D}}_{r_4}^i$ and $\rho_k'\to\rho_\infty'\in\overline{\mathcal{D}}_{r_4}^j$. By the uniform convergence of the $\tilde{u}_{i,k}$ and $\tilde{u}_{j,k}$ it follows that $\tilde{u}_{i,\infty}(\rho_\infty)=\tilde{u}_{j,\infty}(\rho_\infty')$, or in other words $\rho_\infty\in \mathcal{U}_{ij}$ and $\rho_\infty'\in \mathcal{U}_{ji}$.  Furthermore,
\begin{align*}
\dist_{u_k^*g}&\big(\phi_{i,k}(\rho_\infty),\phi_{j,k}(\rho_\infty')\big)\\
&\leq\dist_{\tilde{u}_{i,k}^*g}(\rho_k,\rho_\infty)+\dist_{u_k^*g}\big(\phi_{i,k}(\rho_k),\phi_{j,k}(\rho_k')\big)+ \dist_{\tilde{u}_{j,k}^*g}(\rho_k',\rho_\infty')\\
&\to 0.
\end{align*}
Thus by Lemma \ref{lem:transitionCharacterization}, we have  $\rho_\infty\in\tilde{\mathcal{U}}_{ij}$, but the latter is an open set, so $\rho_k\in\tilde{\mathcal{U}}_{ij}$ for all sufficiently large $k$.  This contradiction completes the proof of the first part of the corollary.

To prove the second part, argue by contradiction.  Indeed, if the second part were not true, then there would exist a closed (and hence compact) set $\mathcal{V}\subset \widetilde{\mathcal{U}}_{ij}$ and open set $\mathcal{O}\subset \mathcal{D}_{r}^j$ which contains the image of $\mathcal{V}$ via $\tilde{u}_{j,\infty}^{-1}\circ \tilde{u}_{i,\infty}$ with the property that there exist arbitrarily large $k\in \mathbb{N}$ for which $\phi_{i,k}(\mathcal{V})\nsubseteq \phi_{j,\infty}(\mathcal{O})$.  After passing to a further subsequence, one constructs a sequence of points $\hat{\rho}_k\subset \mathcal{V}$ which converge to $\hat{\rho}_\infty\in \mathcal{V}$ and have the the property that
\begin{equation*}
\liminf \dist_{u_k^*g}\big(\phi_{i,k}(\hat{\rho}_k),\phi_{j,k}\circ\tilde{u}_{j,\infty}^{-1}\circ\tilde{u}_{i,\infty}(\mathcal{V})\big)>0.
\end{equation*}
However, this implies that
\begin{equation*}
\liminf \dist_{u_k^*g}\big(\phi_{i,k}(\hat{\rho}_\infty),\phi_{j,k}\circ \tilde{u}_{j,\infty}^{-1}\circ \tilde{u}_{i,\infty}(\hat{\rho}_\infty)>0,
\end{equation*}
but this of course is a contradiction since $\hat{\rho}_\infty\in \mathcal{V}\subset \widetilde{\mathcal{U}}_{ij}$.  This completes the proof of the corollary.
\end{proof}

We now define the  topological space $S_\infty:=\sqcup_{i=1}^n\mathcal{D}_{r_1}^j/\sim$ where $\rho\sim \rho'$ provided $\rho\in \widetilde{\mathcal{U}}_{i,j}$, $\rho'\in \widetilde{\mathcal{U}}_{ji}$, and $\tilde{u}_{i,\infty}(\rho)=\tilde{u}_{j,\infty}(\rho')$.  We now claim the following.

\begin{lemma}\label{lem:limitSurfaceAManifold}
$S_\infty$ is a smooth manifold.
\end{lemma}

The proof of Lemma \ref{lem:limitSurfaceAManifold} is elementary; we do not provide it here. We now turn our attention to constructing the desired reparameterizations of the $u_k$. To that end we define the local normal bundles $\mathcal{E}^i\overset{\pi_i}{\to}\tilde{u}_{i,\infty}(\mathcal{D}_{r_1}^i)$ with total space
\begin{equation*}
\mathcal{E}^i:=\{X \in T_p M: p=\tilde{u}_{i,\infty}(\rho)\;\;\text{and}\;\; 0=\la X,Y\ra_g\; \;\forall Y \in T_p\tilde{u}_{i,\infty}(\mathcal{D}_{r_1}^i)\}.
\end{equation*}
Given a smooth unitary trivialization $\Phi_i:\mathcal{E}^i\to\tilde{u}_{i,\infty}(\mathcal{D}_{r_1}^i)\times \mathbb{R}^{\dim M-\dim S}$ we can consider the map
\begin{align*}
&\Psi_i:\tilde{u}_{i,\infty}(\mathcal{D}_{r_1}^i)\times\mathcal{B}_\epsilon\to M\\
&\Psi_i(\rho,X):=\exp_{\pi_i\circ \Phi_{i}^{-1}(\rho,X)}^g\big(\Phi_i^{-1}(\rho,X)\big),
\end{align*}
which is a diffeomorphism with its image provided $\epsilon>0$ is sufficiently small; here $\mathcal{B}_\epsilon:=\{\mathbb{R}^{\dim M-\dim S}:\|X\|<\epsilon\}$, and $\exp^g$ is the exponential map associated to the metric $g$.  Since $i=1,\ldots,n$, let us suppose that $\epsilon$ is sufficiently small so that each of these maps is a diffeomorphism with its image, and let us denote these images by $\mathcal{N}^i$.  Recall that by construction the maps $\tilde{u}_{j,\infty}^{-1}\circ\tilde{u}_{i,\infty}:\widetilde{\mathcal{U}}_{ij}\to\widetilde{\mathcal{U}}_{ji}$ are diffeomorphisms.  Thus by construction of the $\Psi_i$ we see that the maps
\begin{equation*}
\Psi_j^{-1}\circ\Psi_i:\tilde{u}_{i,\infty}(\widetilde{\mathcal{U}}_{ij})\times\mathcal{B}_\epsilon \to \tilde{u}_{j,\infty}(\widetilde{\mathcal{U}}_{ji})\times\mathcal{B}_\epsilon
\end{equation*}
are bundle isomorphisms.  In other words these maps are homomorphisms (in fact diffeomorphims) for which $\pr_1 \circ \Psi_j^{-1}\circ\Psi_i =  \pr_1$, and the $\Psi_j^{-1}\circ\Psi_i$ are linear maps on the fibers;  here $\pr_1$ is the canonical projection to the first component of the cartesian product. It will also be convenient to define the maps
\begin{equation*}
\tilde{\pi}_i:\mathcal{N}^i\to\mathcal{D}_{r_1}^i\qquad\text{by}\qquad\tilde{\pi}_i:=\tilde{u}_{i,\infty}^{-1}\circ\pr_1\circ\Psi_i^{-1}.
\end{equation*}
We also recall that the maps $\tilde{u}_{i,k}:=u_k\circ\phi_{i,k}:\mathcal{D}_{r_1}^i\to M$ converge to $\tilde{u}_{i,\infty}$, and thus it follows that for each $r_2\in(0,r_1)$ we have $\tilde{u}_{i,k}(\mathcal{D}_{r_2}^i)\subset \mathcal{N}^i$ for all sufficiently large $k$. Furthermore for all sufficiently large $k$ the maps given by $\tilde{\pi}_i\circ \tilde{u}_{i,k}:\mathcal{D}_{r_2}^i\to \mathcal{D}_{r_1}^i$ are diffeomorphisms with their images.  In fact,
\begin{equation}\label{eq:constructReparams3}
\tilde{\pi}_i\circ \tilde{u}_{i,k} \to Id \qquad\text{in}\quad C^{\infty}(\mathcal{D}_{r_2}^i,\mathcal{D}_{r_1}^i).
\end{equation}
Consequently, for each $r_3\in (0,r_2)$, and for all sufficiently large $k$, we have
\begin{equation}\label{eq:constructReparams4}
\mathcal{D}_{r_3}^i\subset\tilde{\pi}_i\circ \tilde{u}_{i,k}(\mathcal{D}_{r_2}^i) \qquad\text{and}\qquad \tilde{\pi}_i\circ \tilde{u}_{i,k}(\mathcal{D}_{r_3}^i) \subset \mathcal{D}_{r_2}^i,
\end{equation}
in which case we can define the maps
\begin{equation*}
\psi_{i,k}:\mathcal{D}_{r_3}^i\to S_k\qquad\text{by}\qquad
\psi_{i,k}:=\phi_{i,k}\circ(\tilde{\pi}_i\circ \tilde{u}_{i,k})^{-1}.
\end{equation*}
Fix $r_4\in(0,r_3)$ and define the sets $\widehat{\mathcal{U}}_{ij}:=\mathcal{D}_{r_4}^i\cap\tilde{u}_{i,\infty}^{-1}\circ\tilde{u}_{j,\infty}\big(\mathcal{D}_{r_4}^j\cap\widetilde{\mathcal{U}}_{ji}\big)$,
and the smooth manifold
\begin{equation*}
\widetilde{S}_\infty:=\sqcup_{i=1}^n\mathcal{D}_{r_4}^i/\sim
\end{equation*}
where $\rho\sim\rho'$ provided $\rho\in\widehat{\mathcal{U}}_{ij}$ $\rho'\in\widehat{\mathcal{U}}_{ji}$ and $\tilde{u}_{i,\infty}(\rho)=\tilde{u}_{j,\infty}(\rho')$.  With these definitions made, we now claim the following.

\begin{lemma}\label{lem:descendanceOfPsi}
The maps $\psi_{i,k}$ defined above, with domains restricted to $\mathcal{D}_{r_4}$, descend to smooth maps $\psi_k:\widetilde{S}_\infty\to S_k$.
\end{lemma}
\begin{proof}
We begin by comparing the $\phi_{i,k}$ to the $\psi_{i,k}$.  Indeed, by applying $\psi_{i,k}$ to each side of the second containment of (\ref{eq:constructReparams4}) we see that
\begin{equation}
\phi_{i,k}(\mathcal{D}_{r_3}^i)\subset \psi_{i,k}(\mathcal{D}_{r_2}^i),
\end{equation}
and thus the maps $\psi_{i,k}^{-1}\circ \phi_{i,k}:\mathcal{D}_{r_3}^i\to\mathcal{D}_{r_2}^i$ are well defined.  Furthermore,
\begin{align}
\psi_{i,k}^{-1}\circ\phi_{i,k}&=\Big(\phi_{i,k}\circ(\tilde{\pi}_i\circ \tilde{u}_{i,k})^{-1}\Big)^{-1}\circ\phi_{i,k}\notag\\
&=\tilde{\pi}_i\circ \tilde{u}_{i,k}\notag\\
&\to Id\quad\text{in}\;\; C^\infty(\mathcal{D}_{r_3}^i,\mathcal{D}_{r_2}^i);\label{eq:constructReparams5}
\end{align}
here the convergence in the last line follows from (\ref{eq:constructReparams3}).  Consequently for all sufficiently large $k$, the maps $\psi_{i,k}:\mathcal{D}_{r_3}^i\to \psi_{i,k}(\mathcal{D}_{r_3}^i)\subset S_k$ are diffeomorphisms. Recall Lemma \ref{lem:transitionCharacterization} which guarantees that if $\rho\in \widehat{\mathcal{U}}_{ij}\subset \widetilde{\mathcal{U}}_{ij}$, $\rho'\in \widehat{\mathcal{U}}_{ji}\subset \widetilde{\mathcal{U}}_{ji}$, and  $\tilde{u}_{i,\infty}(\rho)=\tilde{u}_{j,\infty}(\rho')$ then
\begin{equation*}
\lim_{k\to\infty}\dist_{u_k^*g}\big(\phi_{i,k}(\rho),\phi_{j,k}(\rho')\big)=0,
\end{equation*}
and combining this with (\ref{eq:constructReparams5}) yields
\begin{equation*}
\lim_{k\to\infty}\dist_{u_k^*g}\big(\psi_{i,k}(\rho),\psi_{j,k}(\rho')\big)=0.
\end{equation*}
We now claim that for all sufficiently large $k$, we have
\begin{equation}\label{eq:constructReparams6}
\psi_{j,k}(\widehat{\mathcal{U}}_{ji})\subset \psi_{i,k}(\widetilde{\mathcal{U}}_{ij})\qquad\text{and}\qquad
\psi_{j,k}(\widehat{\mathcal{U}}_{ji})\subset \phi_{i,k}(\widetilde{\mathcal{U}}_{ij}).
\end{equation}
Indeed, to see this observe that $\cl (\widehat{\mathcal{U}}_{ji})\subset \widetilde{\mathcal{U}}_{ji}$, and $\phi_{jk}^{-1}\circ\psi_{jk}\to Id$, so that by the latter part of Corollary \ref{cor:reparamConstruction} we see that
for all sufficiently large $k$ we have
\begin{equation*}
\phi_{j,k}\circ\phi_{j,k}^{-1}\circ \psi_{j,k}(\widehat{\mathcal{U}}_{ji})\subset \phi_{i,k}(\widetilde{\mathcal{U}}_{ij}).
\end{equation*}
From this we conclude the second containment in (\ref{eq:constructReparams6}).  The first containment is similarly obtained from
\begin{equation*}
\phi_{j,k}\circ\phi_{j,k}^{-1}\circ \psi_{j,k}(\widehat{\mathcal{U}}_{ji})\subset \phi_{i,k}\circ\phi_{i,k}^{-1}\circ \psi_{i,k}(\widetilde{\mathcal{U}}_{ij}).
\end{equation*}
Now let $\rho\in \widehat{\mathcal{U}}_{ij}$ and $\rho'\in \widehat{\mathcal{U}}_{ji}$ with $\tilde{u}_{i,\infty}(\rho)=\tilde{u}_{j,\infty}(\rho')$; to prove the lemma, we must show that $\psi_{i,k}(\rho)=\psi_{j,k}(\rho')$.  As a consequence of (\ref{eq:constructReparams6}) there exist $\tilde{\rho}_k\in\mathcal{D}_{r_3}^i$ such that $\psi_{i,k}(\tilde{\rho}_k)=\psi_{j,k}(\rho')$.  However, observe that

\begin{align*}
\tilde{u}_{j,\infty}(\rho')&=\tilde{u}_{j,\infty}\circ\psi_{j,k}^{-1}\big(\psi_{j,k}(\rho')\big)\\
&=\tilde{u}_{j,\infty}\circ\tilde{\pi}_j\circ \tilde{u}_{j,k} \circ \phi_{j,k}^{-1}\big(\psi_{j,k}(\rho')\big)\\
&=\tilde{u}_{j,\infty}\circ\tilde{u}_{j,\infty}^{-1}\circ\pr_1\circ\Psi_j^{-1}\circ u_k \big(\psi_{j,k}(\rho')\big)\\
&=\pr_1\circ\Psi_j^{-1}\circ\Psi_i\circ\Psi_i^{-1}\circ u_k \big(\psi_{j,k}(\rho')\big)\\
&=\pr_1\circ\Psi_i^{-1}\circ u_k \big(\psi_{j,k}(\rho')\big)\\
&=\tilde{u}_{i,\infty}\circ \tilde{\pi}_i \circ u_k \big(\psi_{j,k}(\rho')\big)\\
&=\tilde{u}_{i,\infty}\circ \tilde{\pi}_i \circ u_k \big(\psi_{i,k}(\tilde{\rho}_k)\big)\\
&=\tilde{u}_{i,\infty}\circ \tilde{\pi}_i \circ \tilde{u}_{i,k}\circ \phi_{i,k}^{-1}\big(\psi_{i,k}(\tilde{\rho}_k)\big)\\
&=\tilde{u}_{i,\infty}(\tilde{\rho}_k).
\end{align*}
Recall that $\tilde{u}_{i,\infty}(\rho)=\tilde{u}_{j,\infty}(\rho')$, and thus $\tilde{u}_{i,\infty}(\tilde{\rho}_k)=\tilde{u}_{i,\infty}(\rho)$, however $\tilde{u}_{i,\infty}$ is an embedding.  Therefore $\tilde{\rho}_k=\rho$, and thus $\psi_{i,k}(\rho) = \psi_{i,k}(\tilde{\rho}_k) = \psi_{j,k}(\rho')$, and thus the $\psi_{i,k}$ do indeed descend to maps $\psi_k$ on $\widetilde{S}_\infty$.
\end{proof}

We are now ready to finish the proof of Proposition \ref{prop:KlocImpliesK}.  We begin by observing that the above results hold for all $r_i$ whenever $0<r_4<r_3<r_2<r_1<r$. Since the $\phi_{i,k}:\mathcal{D}_r^i\to S_k$ form a sequence of uniformly robust $\mathcal{K}$-covers, it follows that we may choose $r_6\in(0,r_4)$ so that the maps $\phi_{i,k}:\mathcal{D}_{r_6}^i\to S_k$ also form a sequence of uniformly robust $\mathcal{K}$-covers of the $(u_k,S_k)$. We let $\widetilde{\mathcal{K}}\subset \Int (M)$ be a compact set (the existence of which is guaranteed by the definition of a robust cover) whose interior contains $\mathcal{K}$ and for which
\begin{equation}\label{eq:construction}
u_{k}\big(S_k\setminus \cup_{i=1}^n\phi_{i,k}(\mathcal{D}_{r_6}^i)\big)\subset M\setminus \widetilde{\mathcal{K}}
\end{equation}
for all $k$.  Next we observe that $u_k\circ \phi_{i,k} \to \tilde{u}_{i,\infty}$ in $C^\infty(\mathcal{D}_{r_1}^i,M)$, and these limit maps descend to a smooth immersion $\tilde{u}_\infty: S_\infty\to M$.  We have also seen that the maps $\psi_{i,k}:\mathcal{D}_{r_4}^i\to S_k$ descend to $\psi_k:\widetilde{S}_\infty\to S_k$, and they have the property that the sequence $u_k\circ\psi_k:\widetilde{S}_\infty\to M$ converges in $C^\infty$.  We then fix $r_5\in (r_6,r_4)$ and define $\widetilde{S}\subset \widetilde{S}_\infty$ to be a compact region\footnote{See Definition \ref{def:CompactRegion}.} for which
\begin{equation*}
\cup_{i=1}^n \phi_{i,k}(\mathcal{D}_{r_6}^i) \subset \psi_k(\widetilde{S})\subset \cup_{i=1}^n \phi_{i,k}(\mathcal{D}_{r_5}^i),
\end{equation*}
for all sufficiently large $k$. We then note that by (\ref{eq:construction}) we have $u_k\big(S_k\setminus \psi_k(\widetilde{S})\big)\subset M\setminus \widetilde{\mathcal{K}}$.  Thus all that remains to finish the proof is to show the maps $\psi_k:\widetilde{S}\to S_k$ are diffeomorphims with their images.  Since  $\psi_{i,k}^{-1}\circ\phi_{i,k}\to Id$ in $C^\infty(\mathcal{D}_{r_3},\mathcal{D}_{r_2})$, it follows that the $\psi_k$ are immersions, and since $\dim\widetilde{S}=\dim S_k$ it follows that it is sufficient to show that the $\psi_k:\widetilde{S}\to S_k$ are one-to-one for all sufficiently large $k$.

To prove this, suppose not.  Then after possibly passing to a subsequence, there exist sequences of points $[\rho_k],[\rho_k']\in\widetilde{S}\subset \widetilde{S}_\infty$ with representatives $\rho_k\in\mathcal{D}_{r_5}^i$ and $\rho_k'\in\mathcal{D}_{r_5}^j$ for which neither $i$ nor $j$ depend on $k$, $[\rho_k]\neq [\rho_k']$, and $\psi_{i,k}(\rho_k)=\psi_{j,k}(\rho_k')$.  Next observe that the $\psi_{i,k}:\mathcal{D}_{r_3}^i\to \psi_{i,k}(\mathcal{D}_{r_3}^i)$ are diffeomorphisms  for all sufficiently large $k$ (and similarly for $j$), so if $i=j$, then $\rho_k=\rho_k'$ and hence $[\rho_k]=[\rho_k']\in \widetilde{S}$ which is a contradiction.  Thus we henceforth assume that $i\neq j$.

Since $\psi_{i,k}(\rho_k) = \psi_{j,k}(\rho_k') \in S_k$, and $\psi_{i,k}^{-1}\circ\phi_{i,k}\to Id$ in $C^\infty(\mathcal{D}_{r_3}^i)$ and similarly for $j$, we conclude that $\dist_{u_k^* g}\big(\phi_{i,k}(\rho_k),\phi_{j,k}(\rho_k') \big)\to 0$. We then apply the first part of Corollary \ref{cor:reparamConstruction}, and conclude that for all sufficiently large $k$, $\rho_k\in \widetilde{\mathcal{U}}_{ij}\cap \mathcal{D}_{r_5}^i$ and $\rho_k'\in \widetilde{\mathcal{U}}_{ji}\cap \mathcal{D}_{r_5}^j$.  Consequently for all sufficiently large $k$ there exist $\tilde{\rho}_k\in \widetilde{\mathcal{U}}_{ij}\subset \mathcal{D}_{r_1}^i$ with the property that $\tilde{u}_{i,\infty}(\tilde{\rho}_k)=\tilde{u}_{j,\infty}(\rho_k')$.  Our goal now is to show that for some large $k$, we have $\tilde{\rho}_k\in \mathcal{D}_{r_4}^i$.  Indeed, if this were true, then $\psi_{i,k}(\tilde{\rho}_k)$ is well defined; furthermore, $\tilde{\rho}_k\in \widetilde{\mathcal{U}}_{ij}\cap\mathcal{D}_{r_4}^i$ and $\tilde{u}_{j,\infty}^{-1}\circ \tilde{u}_{i,\infty}(\tilde{\rho}_k)=\rho_k'\in \widetilde{\mathcal{U}}_{ji}\cap \mathcal{D}_{r_4}^j$ (in other words, $\tilde{\rho}_k\in \widehat{\mathcal{U}}_{ij}$), so that $\psi_{i,k}(\tilde{\rho}_k)=\psi_{j,k}(\rho_k')$.  However by assumption $\psi_{i,k}(\rho_k)=\psi_{j,k}(\rho_k')$, so $\psi_{i,k}(\rho_k)=\psi_{i,k}(\tilde{\rho}_k)$; but for all sufficiently large $k$, $\psi_{i,k}$ is a diffeomorphism with its image, so we conclude that $\rho_k=\tilde{\rho}_k$.  Then we would have shown that $\rho_k\in \widetilde{\mathcal{U}}_{ij}\cap \mathcal{D}_{r_4}^i$, $\rho_k'\in \widetilde{\mathcal{U}}_{ji}\cap \mathcal{D}_{r_4}^j$, and $\tilde{u}_{i,\infty}(\rho_k)=\tilde{u}_{j,\infty}(\rho_k')$. In other words $[\rho_k]=[\rho_k']\in \widetilde{S}$, which is the desired contradiction.

We have so far shown that to complete the proof of Proposition \ref{prop:KlocImpliesK}, it is sufficient to show that for some large $k$, we have $\tilde{\rho}_k\in\mathcal{D}_{r_4}^i$. To that end, we pass to a further subsequence so that $\rho_k\to \rho_\infty\in \cl(\widetilde{\mathcal{U}}_{ij})\cap \cl(\mathcal{D}_{r_5}^i)$, and $\tilde{\rho}_k\to \tilde{\rho}_\infty \in \cl(\mathcal{D}_{r_1}^i)$.  By the definition of $\tilde{\rho}_k$ and $\tilde{\rho}_\infty$, and the uniform convergence of the $\tilde{u}_{i,k}$ and $\tilde{u}_{j,k}$ to an embedding, it follows that $\tilde{u}_{i,\infty}(\rho_\infty)=\tilde{u}_{i,\infty}(\tilde{\rho}_\infty)$.  Consequently, for all sufficiently large $k$ we have $\tilde{\rho}_k\in \mathcal{D}_{r_4}^i$, and due to the discussion in the previous paragraph, this shows that we have completed the proof of Proposition \ref{prop:KlocImpliesK}.
\end{proof}

\subsection{Proof of Proposition \ref{prop:GromovWithNiceBdry}}\label{sec:GromovWithNiceBdry}
\begin{proof}

In what follows, it will be convenient to have the following notation:
\begin{equation}\label{eq:standardCylinder}
\Sigma_{\mu}:=\big(\mathbb{R}/2\pi\mathbb{Z}\big)\times [0,\mu),
\end{equation}
with $\mu\in \mathbb{R}^+\cup\{\infty\}$.  We shall call the product metric $\bar{g}:=dx^2+dy^2$ the \emph{standard metric} on $\Sigma_\mu$, and $[\bar{g}]$ the \emph{standard conformal structure}. We also abuse notation by defining $\Sigma_0:=\big(\mathbb{R}/2\pi\mathbb{Z}\big)\times\{0\}$.  Also, by assumption $S$ has a finite number of connected components, so without loss of generality we shall assume that $S$ is connected.

Much of the proof is standard, so we focus primarily on the less standard aspects, namely showing the existence of reparameterizations of the $u_k:S\to M$ which have the desired boundary convergence. To that end we recall several  important results.

\begin{lemma}[uniformization]\label{lem:uniformization}
Let $(S,g)$ be a smooth connected compact Riemannian manifold of dimension two with boundary, and consider the finite set $\Gamma\subset S\setminus \partial S$.  If $\chi(S)-\#\Gamma < 0$, then there exists a unique smooth geodesically complete metric $h$ on $\dot{S}:=S\setminus\Gamma$ in the conformal class of $g$ such that $\Area_h(\dot{S}) <\infty$; furthermore the Gauss curvature of $h$ is identically $-1$, and the boundary components of $S$ are all $h$-geodesics.
\end{lemma}
\begin{proof}This is a well known result.  A proof via variational partial differential equation methods in the case that $\Gamma=\emptyset = \partial S$ case can be found in \cite{Ta92}.  The case with boundary can be treated by modifying the argument in \cite{Ta92} to consider an associated Neumann boundary value problem.  The case with punctures can be treated by removing disks of arbitrarily small radius centered at points in $\Gamma$ and taking limits.
\end{proof}
\begin{lemma}[conformal distance]\label{lem:conformalDistance}
Consider the half-cylinders $\Sigma_\mu$ endowed with the standard conformal structure and with moduli $\mu\in \mathbb{R}^+\cup\{\infty\}$. Then for each number $r>0$ there exists a number $\ell=\ell(r)>0$ with the following significance.  If $\mathcal{U}\subset \Sigma_\infty$ is conformally diffeomorphic to $\Sigma_r$, and $\Sigma_0\subset \mathcal{U}$, then $\Sigma_\ell\subset \mathcal{U}$.
\end{lemma}
\begin{proof}This is a restatement of Lemma 2.1 of \cite{IsSv00}.
\end{proof}

\begin{lemma}[quasi-conformal estimate]\label{lem:quasiConformalEstimate}
Let $(S,g)$ be a two dimensional Riemannian manifold, and let $\Sigma_\mu$ be equipped with the standard metric and conformal structure as above. Suppose furthermore there is an annular region $\mathcal{A}\subset S$ and a diffeomorphism (but not necessarily conformal) $\psi:\Sigma_\mu\to \mathcal{A}$ for which
\begin{equation*}
\sup_{\rho\in \Sigma_\mu}\|T_\rho\psi\|\|T_{\psi(\rho)}\psi^{-1}\|\leq C <\infty;
\end{equation*}
here $\|\cdot\|$ denotes the norm of a linear map between normed vector spaces. Then, letting $mod_{[g]}(\mathcal{A})$ denote the modulus of the cylinder $(\mathcal{A},[g])\subset (S,[g])$, we have
\begin{equation*}
C^{-1}\mu\leq mod_{[g]}(\mathcal{A})\leq C\mu.
\end{equation*}
\end{lemma}
\begin{proof}
By the uniformization theorem, it is sufficient to prove the result for $S=\mathbb{R}^2$ with the standard conformal structure.  The result then follows from Lemmas 2.3.1 and 2.3.2 in \cite{SmSt92}.
\end{proof}

The following lemma follows from a straight-forward computation.  We state the result here for convenient reference later.

\begin{lemma}[model hyperbolic cylinders]\label{lem:modelHyperbolicCylinders}
The metric $h:=\cosh^2(t)ds^2 + dt^2$ on $\mathbb{R}^2$, satisfies the following properties.
\begin{enumerate}
\item $h$ is a hyperbolic metric; i.e. $h$ has constant Gauss curvature equal to $-1$.
\item If $(S,\tilde{h})$ is a two-dimensional Riemannian manifold equipped with a hyperbolic metric $\tilde{h}$, and $\alpha:(s_0,s_1)\to S$ is an $\tilde{h}$-unit speed geodesic, and $\nu$ is a continuous $\tilde{h}$-unit normal vector field along $\alpha$, then $h=\phi^*\tilde{h}$ where
    \begin{equation*}
    \phi(s,t)=\exp_{\alpha(s)}^{\tilde{h}}(t\nu).
    \end{equation*}
\item The metric $h$ descends to $(\mathbb{R}/\ell\mathbb{Z})\times \mathbb{R}$. Furthermore
the map given by
\begin{align*}
&\tilde{\phi}:(\mathbb{R}/\ell\mathbb{Z})\times \mathbb{R}\to \big(\mathbb{R}/2\pi\mathbb{Z}\big)\times (-\pi^2\ell^{-1},\pi^2\ell^{-1})\\
&\tilde{\phi}(s,t)=\big([2\pi\ell^{-1}s],2\pi\ell^{-1}\arctan(\sinh(t))\big)
\end{align*}
is a conformal diffeomorphism for the conformal structures associated to $h$ on the domain of $\tilde{\phi}$ and the standard Euclidean metric on the range of $\tilde{\phi}$.
\end{enumerate}
\end{lemma}

Let us now proceed with the proof of Lemma \ref{prop:GromovWithNiceBdry}.

\begin{lemma}[convenient marked points]\label{lem:convenientMarkedPoints}
With $(S,j_k)$ as in Proposition \ref{prop:GromovWithNiceBdry}, and $S$ connected, there exists a sequence of finite sets $\Gamma_k\subset S\setminus \partial S$ with the following properties.
\begin{enumerate}
\item $\#\Gamma_k = N$
\item $\chi(S)-\#\Gamma_k <0$
\item The hyperbolic metrics $h_k$ on $S\setminus\Gamma_k$ guaranteed by Lemma \ref{lem:uniformization}, have the property that each connected component of $\partial S$ has length uniformly bounded away from $0$ and $\infty$.
\end{enumerate}

\end{lemma}
\begin{proof}
Observe that if $\partial S=\emptyset$, then Lemma \ref{lem:convenientMarkedPoints} is trivially true, so henceforth we assume $\partial S\neq \emptyset$. Letting $\psi_{i,k}:\Sigma_\epsilon=\mathbb{S}^1\times [0,\epsilon)\to S$ be the maps as in the assumptions of Proposition \ref{prop:GromovWithNiceBdry}, we denote the metrics $\tilde{g}_{i,k}:=(u_{k}\circ\psi_{i,k})^*g_{i,k}$, which converge in $C^\infty$ to the metrics $\tilde{g}_{i,\infty}$ on $\Sigma_\epsilon$.  By the uniformization theorem (and possibly restricting the domains of the $\psi_{i,k}$) we may assume that for each $i$ we have $\tilde{g}_{i,\infty}=e^{f_i}(dx^2 + dy^2)$; in other words the limiting conformal structures on $\Sigma_\epsilon$ are standard.  We then define $\tilde{\Gamma}_{i,k}:= \{\psi_{i,k}(0,\epsilon/4),\psi_{i,k}(0,\epsilon/3),\psi_{i,k}(0,\epsilon/2)\}$, and $\tilde{\Gamma}_k:=\cup_i \tilde{\Gamma}_{i,k}$.  Since $S$ is connected with at least one boundary component, it follows that $\chi(S)-\#\tilde{\Gamma}_k < 0$; we then let $\tilde{h}_k$ be the hyperbolic metrics on  $S\setminus \tilde{\Gamma}_k$ guaranteed by Lemma \ref{lem:uniformization}.

We now claim that the $\tilde{h}_k$ length of each boundary component of $S$ is uniformly bounded. To prove this, we observe that since the $[\tilde{g}_{i,\infty}]=[\bar{g}]$ on $\Sigma_\epsilon$, it follows from Lemma \ref{lem:quasiConformalEstimate} that each boundary component of $S$ has an annular neighborhood of modulus at least $\epsilon/5$.  Then Lemma \ref{lem:modelHyperbolicCylinders} guarantees that the $\tilde{h}_k$-lengths of the boundary components of $S$ are uniformly bounded.

Before completing the proof of Lemma \ref{lem:convenientMarkedPoints}, we will need to make use of the following result.

\begin{lemma}[hyperbolic neighborhood]\label{lem:hyperbolicNeighborhood}
Let $S$ and  $\psi_{i,k}$ be as in Proposition \ref{prop:GromovWithNiceBdry}, and let $\tilde{\Gamma}_k$, $\tilde{h}_k$, $\Sigma_\epsilon$ be as above.  Let $\partial_i S$ denote the $i^{th}$ boundary component of $S$, and define the open set
\begin{equation*}
\mathcal{O}_\delta^{\tilde{h}_k}(\partial S):=\{\zeta\in S:\dist_{h_k}(\zeta,\partial S)< \delta\}.
\end{equation*}
Then there exists $\mu>0$ and $\delta_{i,k}>0$ such that the
$\mathcal{O}_{\delta_{i,k}}^{\tilde{h}_k}(\partial_i S)$
are annular neighborhoods of the $\partial_i S$ with modulus equal to $\mu$, and are contained in every annular
neighborhood of $\partial_i S$ with modulus at least $\epsilon/5$.
\end{lemma}
\begin{proof}
We begin by defining $\dot{S}_k:=S\setminus\tilde{\Gamma}_k$, and the doubled surfaces
\begin{equation*}
2\dot{S}_k:=\big(\dot{S}_k\sqcup\dot{S}_k\big)/\sim
\end{equation*}
where $\sim$ is the identification via the identity map along\footnote{To be clear, $\partial \dot{S}_k=\partial S$; or in other words $\partial \dot{S}_k$ does not contain the "degenerate boundary components" $\tilde{\Gamma}_k$.} $\partial \dot{S}_k$. Observe that since the components of $\partial S$ are $\tilde{h}_k$-geodesics, it follows that the $\tilde{h}_k$ extend via reflection to smooth hyperbolic metrics on $2\dot{S}_k$; we abuse notation by also denoting these metrics $\tilde{h}_k$.  Next note that there is a natural inclusion $\partial \dot{S}_k\hookrightarrow 2\dot{S}_k$, and the image of the $i^{th}$ boundary component is a simple $\tilde{h}_k$-geodesic of length $\tilde{\ell}_{i,k}\leq C<\infty$. Next, it is straight-forward to show that the maps
\begin{equation*}
\phi_{i,k}:\big((\mathbb{R}/\tilde{\ell}_{i,k}\mathbb{Z})\times \mathbb{R},h\big)\to (2\dot{S}_k,\tilde{h}_k),
\end{equation*}
defined as in property 2 of Lemma \ref{lem:modelHyperbolicCylinders} (and associated to the simple closed geodesics $\partial_i \dot{S}_k$ of length $\tilde{\ell}_{i,k}$), are isometric covering maps and hence conformal.  Combining these maps with the conformal diffeomorphisms
\begin{equation*}
\varphi_{i,k}:\big((\mathbb{R}/\tilde{\ell}_{i,k}\mathbb{Z})\times[0,\delta),[h]\big)\to\big(\Sigma_{2\pi\tilde{\ell}_{i,k}^{-1}\arctan(\sinh(\delta))},[\bar{g}]\big)
\end{equation*}
given by property 3 of Lemma \ref{lem:modelHyperbolicCylinders}, we observe that any annular neighborhood of $\partial_i \dot{S}_k$ can be conformally lifted by $\varphi_{i,k}\circ\phi_{i,k}^{-1}$ to annular neighborhoods of $\Sigma_0 \subset \Sigma_\infty$ with the standard conformal structure. Next, observe that as a consequence of Lemma \ref{lem:conformalDistance}, there exists a $\mu>0$ such that any annular neighborhood of $\Sigma_0 \subset \Sigma_\infty$ of modulus at least $\epsilon/5$ contains $\Sigma_\mu$.  Thus to complete the proof of Lemma \ref{lem:hyperbolicNeighborhood}, it is sufficient to show that there exist $\delta_{i,k}>0$ such that
\begin{equation}\label{eq:modulusEstimate}
\frac{2\pi}{\tilde{\ell}_{i,k}}\arctan(\sinh(\delta_{i,k})) = \mu.
\end{equation}
Without loss of generality, we may assume that $\mu< \inf_{i,k}\pi^2\tilde{\ell}_{i,k}^{-1}$, however in this case $\arctan\circ\sinh$ is invertible, and the existence of the $\delta_{i,k}$ that satisfy equation (\ref{eq:modulusEstimate}) follows immediately. This completes the proof of Lemma \ref{lem:hyperbolicNeighborhood}

\end{proof}

With Lemma \ref{lem:hyperbolicNeighborhood} proved, we now finish the proof of of Lemma \ref{lem:convenientMarkedPoints}.  Indeed, let $\Gamma_k:=\tilde{\Gamma}_k\cup\{\zeta_{1,k},\ldots,\zeta_{n,k}\}$ where $\zeta_{i,k}=\phi_{i,k}\circ\varphi_{i,k}^{-1}(0,\mu/2)$, and $\mu$, $\phi_{i,k}$, and $\varphi_{i,k}$ are defined as in the proof of Lemma \ref{lem:hyperbolicNeighborhood}. By construction, properties (1) and (2) of Lemma \ref{lem:convenientMarkedPoints} are satisfied.  To prove property (3), we note that by definition each boundary component of $S\setminus\Gamma_k$ has an annular neighborhood of modulus $\mu/4>0$, so again by property (3) of Lemma \ref{lem:modelHyperbolicCylinders}, it follows that the the $h_k$-lengths of the components of $\partial S$ are uniformly bounded.  All that remains then is to show that the $h_k$-lengths of the components of $\partial S$ are uniformly bounded away from $0$.  Note that if this were not the case, it would follow from property (3) of Lemma \ref{lem:modelHyperbolicCylinders} that for any fixed $\delta>0$, there would exist an $i$ and $k$ such that $\partial_i S$ has a metric annular neighborhood $\mathcal{O}_{\delta}^{h_k}(\partial_i S)\subset S\setminus \Gamma_k$ of modulus as large as we wish (in particular, greater than $\epsilon/5$).  But then
\begin{equation*}
\phi_{i,k}\circ \varphi_{i,k}^{-1}(\Sigma_\mu)=\mathcal{O}_{\delta_{i,k}}^{\tilde{h}_k}(\partial_i S)\subset \mathcal{O}_{\delta}^{h_k}(\partial_i S),
\end{equation*}
where the equality follows by construction of $\phi$ and $\varphi$, and containment follows from Lemma \ref{lem:hyperbolicNeighborhood} since $mod_{[\tilde{h}_k]}\big(\mathcal{O}_{\delta_{i,k}}^{\tilde{h}_k}(\partial_i S)\big) > \epsilon/5$ by assumption. However, this is impossible because $\Gamma_k\cap \phi_{i,k}\circ\varphi_{i,k}^{-1}(\Sigma_{\mu})\neq \emptyset$ by definition of $\Gamma_k$ and $\Gamma_k\cap \mathcal{O}_{\delta}^{h_k}(\partial_iS\setminus\Gamma_k) = \emptyset$ as a consequence of the definition of $h_k$. This contradiction shows that the $h_k$-lengths of the components of $\partial S$ are uniformly bounded away from zero, and thus the proof of Lemma \ref{lem:convenientMarkedPoints} is complete.
\end{proof}

Before proceeding with the proof of Lemma \ref{prop:GromovWithNiceBdry}, we need one more technical result, namely the following.

\begin{lemma}[convergence near the boundary]\label{lem:convergenceNearBoundary}
Let $\mathbf{u}_k=(u_k,S,j_k,J_k)$, $\psi_{i,k}$, and $\Sigma_\epsilon$ be as in the statement of Proposition \ref{prop:GromovWithNiceBdry}; also let $\varphi_{i,k}$ and $\phi_{i,k}$ be the maps defined as in proof of Lemma \ref{lem:hyperbolicNeighborhood}. Then after passing to a subsequence, there exists $\delta>0$ such that the restricted maps
\begin{equation*}
u_k\circ\phi_{i,k}\circ\varphi_{i,k}^{-1}:\Sigma_\delta\to M
\end{equation*}
converge in $C^\infty(\Sigma_\delta,M)$.
\end{lemma}
\begin{proof}
Recall that by assumption the $u_k\circ\psi_{i,k}$ converge in $C^\infty$, so to prove Lemma \ref{lem:convergenceNearBoundary}, it is sufficient to prove the $\psi_{i,k}^{-1}\circ\phi_{i,k}\circ\varphi_{i,k}^{-1}:\Sigma_\delta\to \Sigma_\epsilon$ converge in $C^\infty$.  To that end we treat these maps as pseudo-holomorphic curves with a real one-dimensional Lagrangian boundary condition.  Indeed, we observe that for some $\delta>0$ we must have uniform gradient bounds, since otherwise one could "bubble-off" a non-constant holomorphic map from $\mathbb{C}$ (or the upper half plane) to a compact set in $\mathbb{S}^1\times \mathbb{R}$, which is impossible.  Elliptic regularity then guarantees $C^\infty$ bounds on $\Sigma_\delta$, and thus by passing to a subsequence we have the desired $C^\infty$ convergence.  This completes the proof of Lemma  \ref{lem:convergenceNearBoundary}.
\end{proof}

Finally, we finish the proof of Proposition \ref{prop:GromovWithNiceBdry}. The remainder of this proof is fairly standard, so we simply sketch the argument.  We first note that by construction, the marked $J_k$-curves $(\mathbf{u}_k,\mu_k)$ with $\mu_k:=\Gamma_k$ are all stable, and these marked curves will remain stable even after more marked points are added.  Next we note that either we have uniform $h_k$-gradient bounds on the $u_k$, or else we don't.  If we do, then Deligne-Mumford compactness (and the uniformization theorem) guarantee the existence of a decorated nodal Riemann surface $(S,j,\mu,D,r)$ and diffeomorphisms $\phi_k:S^{D,r}\to S_k$ such that properties (2) - (4) of Definition \ref{def:GromovConvergence} (i.e. Gromov Convergence) are satisfied.  Elliptic regularity and Arzel\`{a}-Ascoli yield property (6), or rather the desired $C_{loc}^\infty$-convergence away from nodes and boundary.  Smooth convergence in boundary neighborhoods then follows from Lemma \ref{lem:convergenceNearBoundary}.
Property (5), in other words $C_0$-convergence on $S^{D,r}$, then follows from Gromov's removable singularity theorem, monotonicity of area, and the uniform gradient bounds with respect to the hyperbolic metric.

On the other hand, we may not have uniform $h_k$-gradient bounds on the $u_k$.  In this case, one applies the usual bubbling analysis to guarantee the existence of a sequence of finite sets $\hat{\mu}_k\supset \mu_k=\Gamma_k$, which satisfy the conclusions of Lemma \ref{lem:convenientMarkedPoints} and for which one indeed has uniform $\hat{h}_k$-gradient bounds.  Note that as a consequence of Lemma \ref{lem:convergenceNearBoundary}, we have $\dist_{h_k}(\hat{\mu}_k,\partial S_k)\geq \epsilon>0$ for some $\epsilon>0$ independent of $k$.  The arguments of the previous paragraph then apply, and we again conclude Gromov convergence.  This completes the proof of Proposition \ref{prop:GromovWithNiceBdry}.
\end{proof}

\subsection{Proof of Proposition \ref{prop:Desingularization}}\label{sec:Desingularization}

\begin{proof}Fix $\epsilon>0$. Recall Lemma \ref{lem:CritPointLocalModel}, which guarantees that for each $z\in \mathcal{Z}_u$ there exists holomorphic coordinate charts $\phi_z :\mathcal{O}_z(z) \subset S \to \mathcal{O}_z(0)\subset \mathbb{C}\simeq\mathbb{R}^2$ and polar geodesic coordinate charts $\Phi_z:\mathcal{O}_z\big(u(z)\big)\subset M \to \mathcal{O}_z(0)\subset \mathbb{C}^n\simeq\mathbb{R}^{2n}$ which satisfy $\phi_z(z)=0$, $\Phi_{z}\big(u(z)\big)=0$, $(\Phi_{z*}J)(0)=i=J_0$, and
\begin{equation*}
\Phi_z\comp u\comp \phi_z^{-1}(\rho)=\big(\rho^{k_z},0,\ldots,0\big)+F_z(\rho)\in \mathbb{C}^{n} \simeq \R^{2n}
\end{equation*}
where $\phi_z(\zeta)=\rho=s+it$, $F_z(\rho)=O_{k_z+1}(|\rho|^{k_z+1})$, $k_z\geq 2$, and the sub-script $z$ denotes dependance on $z\in Z$. Consequently we make the following definition.
\begin{definition} Let $g_0$ be the standard metric on $\R^{2n}$, and let $g$, $J_0$, $F_z$, $\Phi_z$, $\phi_z$, and $u$ be as above. Then define $\epsilon_0>0$ to be a positive constant for which the following hold.
\begin{enumerate}[(e1)]
\item $\epsilon_0 < \min(1,\epsilon)$.\label{en.e1}
\item The
sets $\mathcal{B}_{\epsilon_0}(z):=\{\zeta\in S: \dist_\gamma(z,\zeta)<\epsilon_0\}$
are pair-wise disjoint as $z$ varies over $\mathcal{Z}_u$.\label{en.e2}
\item $\mathcal{B}_{\epsilon_0}(z)\subset \mathcal{O}_z(z)\subset S$ for all $z\in \mathcal{Z}_u$.\label{en.e3}
\item $|F_z(\rho)|< \epsilon_0$ for all $\rho\in \phi_z\big(B_{\epsilon_0}(z)\big)$ and $z\in \mathcal{Z}_u$.\label{en.e4}
\item $\mathcal{B}_{3\epsilon_0}\big(u(z)\big)\subset\mathcal{O}_z\big(u(z)\big)\subset M$ for each $z\in \mathcal{Z}_u$.\label{en.e5}
\item
  $\|dF_z(\rho)\|_{g_0}\leq \frac{1}{2}|\rho|^{k_z-1}$ for all $\rho\in \phi_z\big(\mathcal{B}_{\epsilon_0}(z)\big)$.\label{en.e6}
\end{enumerate}
\end{definition}
\noindent Next define a smooth cut-off function $\beta:\R\to[0,1]$ for which $\beta'\leq 0$ and
\begin{equation*}
\beta(a)= \left\{
  \begin{array}{ll}
    1, & \text{if } a \leq \frac{1}{4} \\
    0, & \text{if } a \geq \frac{3}{4}.
  \end{array}
\right.
\end{equation*}
Also define the following family of perturbed maps for $\delta\in [0,\epsilon_0)$.
\begin{equation*}
\hat{u}(\zeta)=\left\{
  \begin{array}{ll}
    u(\zeta), & \text{if } \zeta\in S\setminus \bigcup_{z\in \mathcal{Z}_u} \mathcal{B}_{\epsilon_0}(z) \\
    \Phi_z^{-1}\circ \hat{v}_z \circ \phi_z(\zeta),
&\text{if } \zeta\in \mathcal{B}_{\epsilon_0}(z),
  \end{array}
\right.
\end{equation*}
where for each $z\in \mathcal{Z}_u$ we have
$\hat{v}_z:\mathcal{O}_z(0)\subset \mathbb{C}\simeq\R^2\to\mathbb{C}^n\simeq\mathbb{R}^{2n}$ given by
\begin{equation*}
\hat{v}_z(\rho) = \Phi_z\circ u \circ \phi_z^{-1}(\rho) + \big(0,\delta^{k_z-1}\beta(|\rho|/r_0)\rho,0,\ldots,0\big)
\end{equation*}
where $r_0\in (0,\epsilon_0)$ has been chosen so that
\begin{equation}\label{eq:constant}
\mathcal{D}_{r_0}:=\{\rho\in \mathbb{C}:|\rho|< r_0\} \subset \bigcap_{z\in \mathcal{Z}_u}\phi_z\big(\mathcal{B}_{\epsilon_0}(z)\big).
\end{equation}
More concisely, our locally defined family of perturbed maps is given by
\begin{equation*}
\hat{v}(\rho)=\big(\rho^k,\delta^{k-1}\beta(|\rho|/r_0)\rho,0,\ldots,0\big) + F(\rho)
\end{equation*}
where $F(\rho)=O_{k+1}(|\rho|^{k+1})$ and we have stopped denoting $z$ dependance. We now take a moment to verify that the $\hat{u}$ are well-defined.  Indeed, each $\hat{v}_z$ is well defined on $\mathcal{O}_z(0)$, so it is sufficient to show that $\hat{v}_z\circ\phi_z(\zeta)\subset \Phi_z\big(\mathcal{O}_z(u(z))\big)$ for each $z\in \mathcal{Z}_u$ and $\zeta\in \phi_z^{-1}(\mathcal{D}_{r_0})$, and that $u(\zeta)=\Phi_z^{-1}\circ \hat{v}_z\circ \phi_z(\zeta)$ for all $\zeta$ near $\partial \mathcal{B}_{\epsilon_0}(z)\subset S$. To that end, define $v_z:=\hat{v}_z\big|_{\delta=0}$ and observe that $\supp(\hat{v}_z-v_z)\subset \mathcal{D}_{r}\subset \phi_z\big(\mathcal{B}_{\epsilon_0}(z)\big)\subset\mathbb{C}$, and thus indeed $u(\zeta)=\Phi_z^{-1}\circ \hat{v}_z\circ \phi_z(\zeta)$ for all $\zeta$ near $\partial \mathcal{B}_{\epsilon_0}(z)\subset S$.  Also observe that since $\supp(\hat{v}_z-v_z)\subset \mathcal{D}_{r_0}$, it follows that for $\rho\in \mathcal{D}_{r_0}$ we have
\begin{align*}
|\hat{v}_z(\rho)|&\leq |\rho^{k_z}| + \delta^{k_z-1}|\rho|\;\big|\beta(|\rho|/r_0)\big| + |F_z(\rho)|\\
&\leq r_0^{k_z} + r_0\delta^{k_z-1} + |F_z(\rho)|\\
&<3\epsilon_0,
\end{align*}
where we have made use of the fact that $r_0,\delta\leq \epsilon_0< 1$, and (e\ref{en.e4}). Consequently, by (e\ref{en.e5}), we have $\hat{v}_z\circ\phi_z(\zeta)\subset \Phi_z\big(\mathcal{O}_z(u(z))\big)$ for each $z\in \mathcal{Z}_u$ and $\zeta\in \phi_z^{-1}(\mathcal{D}_{r_0})$.  This shows that the $\hat{u}$ are well defined perturbations of $u$.

With the above perturbed maps $\hat{u}$ defined, our next goal is to show that for all sufficiently small $\delta>0$, all six properties of Proposition \ref{prop:Desingularization} are satisfied.  To that end, let $\pr_j:\mathbb{C}^n\to \mathbb{C}$ denote the canonical projection to the $j^{th}$ complex coordinate;  by (e\ref{en.e6}) it follows that $d(\pr_1\circ \hat{v})(\rho)=0$ only if $\rho=0$, and by definition of $\hat{v}$ we have $d(\pr_2\circ \hat{v})(0)\neq 0$ provided $\delta\neq 0$.  Consequently the maps $\hat{v}$ are immersions for all $\delta>0$ sufficiently small.

Observe that Property (D\ref{en.D1}) follows from (e\ref{en.e2}) and equation (\ref{eq:constant}). To prove Property (D\ref{en.D2}), we note that $\hat{u}\to  u$ in $C^0(S,M)$ (moreover in $C^\infty(S,M)$) as $\delta\to 0$, and hence (D\ref{en.D2}) is also satisfied for all sufficiently small $\delta>0$.

We now prove Property (D\ref{en.D3}). First observe that for any compact region $\mathcal{K}\subset S\setminus \mathcal{Z}_u$, we have $\hat{u}\to u$ in $C^1(\mathcal{K},M)$ as $\delta\to 0$, and the limit is an immersed $J$-curve.  Since the limit curve has $J$-invariant tangent planes, it follows that there exists a $\delta'>0$ (dependant on $\mathcal{K}$ and $\epsilon_0$) such that the desired estimate holds for all $\delta\in (0,\delta')$ and $X\in \mathcal{T}_\zeta$ with $\zeta\in \mathcal{K}$.  To prove the result on the complement of $\mathcal{K}$ we work locally and fix $z\in \mathcal{Z}_u$ and define $\widetilde{g}:=\Phi_{z*}g$, $\widetilde{J}:=\Phi_{z*} J$, $F:=F_z$, $k:=k_z$, and $\hat{v}:=\hat{v}_z$.  To finish proving property (D\ref{en.D3}), it is then sufficient to prove the following lemma.
\begin{lemma}\label{lem:DesingD3}
There exist constants $\delta'>0$ and $0<r<r_0/4$ (dependant on $\widetilde{g}$, $\widetilde{J}$, and $F$) with the following significance.  If $\delta\in(0,\delta')$, $|\rho|\leq r$, and $X\in \hat{v}_*(T_\rho\R^2)$ with $\|X\|_{\widetilde{g}}=1$, then
\begin{equation*}
\|(\widetilde{J}X)^{\bot}\|_{\widetilde{g}}\leq \epsilon.
\end{equation*}
\end{lemma}
Before providing the proof of Lemma \ref{lem:DesingD3}, we first introduce some notation which will be useful later on. We will let $C$ (resp. $c$) denote any sufficiently large (resp. small) positive constant depending on $F$, $\widetilde{J}$, and $\widetilde{g}$, but not $\delta$. Next, consider a plane $\mathcal{P}\subset T_q \R^{2n}$. Then we can define the $\widetilde{g}$ and $g_0$ orthogonal projections
\begin{align*}
&\Pi_{(q,\mathcal{P})}^{\top_{\widetilde{g}}}:T_q\R^{2n}\to\mathcal{P} &&\Pi_{(q,\mathcal{P})}^{\top_{g_0}}:T_q\R^{2n}\to\mathcal{P}\\
&\Pi_{(q,\mathcal{P})}^{\bot_{\widetilde{g}}}:T_q\R^{2n}\to\mathcal{P}^{\bot_{\widetilde{g}}} &&\Pi_{(q,\mathcal{P})}^{\bot_{g_0}}:T_q\R^{2n}\to\mathcal{P}^{\bot_{g_0}}.
\end{align*}
Note that we may identify each tangent space $T_q \R^{2n}$ with $\R^{2n}$ via coordinate translation, and in this manner we may regard the above projections simply as maps from $\R^{2n}$ to $\mathcal{P}$, $\mathcal{P}^{\bot_{\widetilde{g}}}$ and $\mathcal{P}^{\bot_{g_0}}\subset \R^{2n}$.  In particular, this allows one to add, subtract, compose, etc. these projections \emph{even with different $(q,\mathcal{P})$.} We clarify this last point.  Without a fixed identification of the fibers $T_q \R^{2n}$, the following quantity would be nonsensical:
\begin{equation*}
\Pi_{(q_1,\mathcal{P}_1)}^{\top_{\widetilde{g}}} (X_{q_1}) + \Pi_{(q_2,\mathcal{P}_2)}^{\top_{\widetilde{g}}} (Y_{q_2}).
\end{equation*}
Moreover, even with the above identification defined, the following statements hold in general
\begin{align}
&\Pi_{(q_1,\mathcal{P})}^{\top_{\widetilde{g}}} \circ \Pi_{(q_2,\mathcal{P})}^{\bot_{\widetilde{g}}} \neq 0\label{eq:OrthogonalProjections1}\\
&\Pi_{(q_1,\mathcal{P})}^{\bot_{\widetilde{g}}} \circ \Pi_{(q_2,\mathcal{P})}^{\top_{\widetilde{g}}} = 0.\label{eq:OrthogonalProjections2}
\end{align}
The point of (\ref{eq:OrthogonalProjections1}) is that in general $q_1\neq q_2$, and thus the inner products $\widetilde{g}\big|_{q_1}$ and $\widetilde{g}\big|_{q_2}$ need not be equal, and thus neither do the orthogonal compliments of $\mathcal{P}$.  Of course, if $q_1=q_2$, then the non-equality in (\ref{eq:OrthogonalProjections1}) should be replaced with an equality.  Given this discussion, one may expect that in general (\ref{eq:OrthogonalProjections2}) should be false, however the point here is that $\Pi_{(q,\mathcal{P})}^{\bot_{\widetilde{g}}}\big|_{\mathcal{P}}\equiv 0$, independent of $q$.  We make use of these facts below.

We now abuse this notation for the application we have in mind. Indeed, for smooth immersions $\varphi,\psi: \mathcal{D}_{r_0}\subset \R^2 \to\R^{2n}$ we will use the notation
\begin{equation*}
\Pi_{(\varphi,T\psi)}^\bot = \Pi_{\big(\varphi(\rho),\psi_*(T_\rho \mathcal{D}_{r_0})\big)}^\bot,
\end{equation*}
and similarly for the other projections.  It will also be convenient to define the complex polynomial $P(\rho)=(\rho^k,\delta^k\rho,0,\ldots,0$, so that for $|\rho|< r_0/4$ we have $\hat{v}(\rho)=P(\rho)+F(\rho)$, where $F(\rho)$ where $F(\rho)=O_{k+1}(|\rho|^{k+1})$. With this notation established, we are now prepared to prove Lemma \ref{lem:DesingD3}.
\begin{proof}[Proof of Lemma \ref{lem:DesingD3}]
Let $X$ be a $\widetilde{g}$-unit vector tangent to the image of $\hat{v}$, and define the following.
\begin{align*}
\mathcal{E}_1:&=\big(\Pi_{(\hat{v},T\hat{v})}^{\bot_{g_0}}-\Pi_{(P,TP)}^{\bot_{g_0}}\big)(J_0 X)\\
\mathcal{E}_2:&=\Pi_{(\hat{v},T\hat{v})}^{\bot_{g_0}}\big((\widetilde{J}-J_0) X\big)\\
\mathcal{E}_3:&=\big(\Pi_{(\hat{v},T\hat{v})}^{\bot_{\widetilde{g}}} -\Pi_{(\hat{v},T\hat{v})}^{\bot_{g_0}}\big)(\widetilde{J} X),
\end{align*}
so that
\begin{equation}\label{eq:DesingPropD3a}
\Pi_{(\hat{v},T\hat{v})}^{\bot_{\widetilde{g}}}(\widetilde{J} X)= \Pi_{(P,TP)}^{\bot_{g_0}}(J_0 X) + \mathcal{E}_1 + \mathcal{E}_2 + \mathcal{E}_3.
\end{equation}

Next we recall the estimates $|\hat{v}(\rho)|\leq C|\rho|$, $\|\widetilde{g}(q)-g_0(q)\|_{g_0}\leq C|q|^2$, and also  $\|\widetilde{J}(q) - J_0(q)\|_{g_0}\leq C|q|$, so that
\begin{equation}\label{eq:DesingPropD3b}
\|\mathcal{E}_2\|_{g_0}+\|\mathcal{E}_3\|_{g_0}\leq C|\rho|.
\end{equation}
To estimate $\mathcal{E}_1$, we note that $\|dP\|_{g_0}=(|\rho|^{2k-2}+\delta^{2k-2})^{\frac{1}{2}}$, and since $P$ is a complex polynomial, it follows that the linear maps
\begin{equation*}
\|dP\|_{g_0}^{-1} TP : T_\rho \R^2 \to T_{P(\rho)} \R^{2n}
\end{equation*}
are $g_0$-isometries. Furthermore $\|dF\|_{g_0}\leq C|\rho|^{k} \leq C|\rho| \|dP\|_{g_0}$, so
\begin{equation*}
\Big\|\|dP\|_{g_0}^{-1} TP -\|dP\|_{g_0}^{-1} T\hat{v}  \Big\|_{g_0}\leq C|\rho|,
\end{equation*}
from which it follows that
\begin{equation}\label{eq:DesingPropD3c}
\Big\|\Pi_{(\hat{v},T\hat{v})}^{\bot_{g_0}}-\Pi_{(P,TP)}^{\bot_{g_0}} \Big\| \leq C|\rho|,
\end{equation}
and thus
\begin{equation}\label{eq:DesingPropD3d}
\|\mathcal{E}_1\|_{g_0}\leq C|\rho|.
\end{equation}
Lastly, we observe that since $P$ is a complex polynomial, $J_0$ preserves the tangent and $g_0$-normal bundles along the image of $P$; consequently $J_0$ and $\Pi_{(P,TP)}^{\bot_{g_0}}$ commute. It then follows from (\ref{eq:DesingPropD3c}) that
\begin{equation}\label{eq:DesingPropD3e}
\|\Pi_{(P,TP)}^{\bot_{g_0}}(J_0 X)\|_{g_0} \leq C|\rho|.
\end{equation}
Combining the above inequalities then yields
\begin{equation*}
\|(JX)^{\bot_{\widetilde{g}}}\|_{g_0} \leq C|\rho|,
\end{equation*}
which then proves Lemma \ref{lem:DesingD3}, and completes the proof of Property (D\ref{en.D3}).
\end{proof}


Observe that $u_\delta\to u$ uniformly in $C^1(S,M)$ (in fact, in $C^\infty$), and $S$ is compact, so that (D\ref{en.D4}) follows immediately.


We now move on to the proof of Property (D\ref{en.D5}). Here we consider two cases: compact sets of $\mathcal{D}_{r_0}\setminus\{0\}$, and small neighborhoods of $0\in \mathcal{D}_{r_0}$.  We handle the former case first.

\begin{lemma}
For each compact set $\mathcal{K}\subset \mathcal{D}_{r_0}\setminus\{0\}$, there exists $\delta'>0$ with the following significance.  For each $\delta\in(0,\delta')$, the following estimate holds for all $\rho\in \mathcal{K}$.
\begin{equation*}
K_{\hat{v}^*\widetilde{g}}(\rho) \leq K_{v^*\widetilde{g}}(\rho) + \epsilon.
\end{equation*}
\end{lemma}
\begin{proof}
Observe that $\hat{v}\to v$ in $C^\infty(\mathcal{D}_{r_0},\R^{2n})$ as $\delta\to 0$, and $v$ is immersed on $\mathcal{D}_{r_0}\setminus\{0\}$, so $K_{\hat{v}^*\tilde{g}}\to K_{v^*\tilde{g}}$ in $C_{loc}^\infty(\mathcal{D}_{r_0}\setminus\{0\},\R^{2n})$.  The result is then immediate.
\end{proof}

The proof in the case of neighborhoods of $0$ is more complicated, however we claim it follows quickly from the following technical result.

\begin{lemma}\label{lem:DesingTechnical}
There exist constants $0<r<r_0/4$ and $\delta'>0$, which depend on $F$ and $\widetilde{g}$ (but not $\delta$) with the following significance. For all $\delta\in (0,\delta')$ and $\rho\in \mathcal{D}_{r_0}$, the following inequality holds.
\begin{equation*}
\frac{\la\two(\hat{v}_s,\hat{v}_s),\two(\hat{v}_t,\hat{v}_t)\ra_{\tilde{g}}}{\|\hat{v}_s\wedge\hat{v}_t\|_{\tilde{g}}^2}\leq \epsilon + \sup_{\substack{q\in M\\ J\mathcal{P}_q=\mathcal{P}_q}} {\textstyle \frac{1}{2}}\|\tr_{\mathcal{P}_q} J\nabla J\|_{\tilde{g}}^2,
\end{equation*}
with notation as above.
\end{lemma}
Before proceeding with the proof of Lemma \ref{lem:DesingTechnical}, let us use it to finish the proof of Property (D\ref{en.D5}).  Indeed, recall that Gauss equations for for immersed surfaces guarantees that
\begin{align*}
K_{\hat{v}^*\widetilde{g}}(\rho) &= K_{sec}\big(\hat{v}_*(T_\rho S)\big) + \frac{\la \two(\hat{v}_s,\hat{v}_s),\two(\hat{v}_t,\hat{v}_t)\ra_{\widetilde{g}}}{\|\hat{v}_s\wedge\hat{v}_t\|_{\widetilde{g}}^2}-\frac{\|\two(\hat{v}_s,\hat{v}_t)\|_{\widetilde{g}}^2}{\|\hat{v}_s\wedge\hat{v}_t\|_{\widetilde{g}}^2}\\
&\leq \sup_{q\in M} |K_{sec}(q)| + \sup_{\substack{q\in M\\J\mathcal{P}_q = \mathcal{P}_q}} {\textstyle \frac{1}{2}}\|\tr_{\mathcal{P}_q} J\nabla J\|_{\tilde{g}}^2+\epsilon,
\end{align*}
which is precisely the desired result.  Thus to prove Property (D\ref{en.D5}), all that remains is to prove Lemma \ref{lem:DesingTechnical}.  To that end, we will make use of our notation from the proof of Property (D\ref{en.D3}) concerning the $\widetilde{g}$-orthogonal projections $\Pi_{(\cdot,\cdot)}^\top$ and $\Pi_{(\cdot,\cdot)}^\bot$.  Furthermore, for the remainder of the section we will regard $v_s,\hat{v}_s,v_{ss},\hat{v}_{ss}$, etc. as either vector fields along the image of $v$ or $\hat{v}$ (as appropriate), or else as maps from $\mathcal{D}_{r_0}$ to $\R^{2n}$, with the distinction determined by context.  Consequently, we may now write the following
\begin{equation*}
B_v(v_s,v_s)= \Pi_{(v,Tv)}^{\bot_{\tilde{g}}}(\nabla_{v_s}v_s)\qquad\text{and}\qquad B_{\hat{v}}(\hat{v}_s,\hat{v}_s)= \Pi_{(\hat{v},T\hat{v})}^{\bot_{\tilde{g}}}(\nabla_{\hat{v}_s}\hat{v}_s),
\end{equation*}
and more importantly it will allow us to estimate quantities like the following:
\begin{equation*}
\big\| \Pi_{(v,T\hat{v})}^{\bot_{\tilde{g}}}(\nabla_{v_s}v_s) - \Pi_{(\hat{v},T\hat{v})}^{\bot_{\tilde{g}}}(\nabla_{\hat{v}_s}\hat{v}_s)\big\|_{g_0}.
\end{equation*}
Here, as above, $g_0$ is the Euclidian metric, and $\nabla$ is the Levi-Civita connection associated to $\widetilde{g}$. We locally define the $(1,2)$ tensor $\Gamma$ by the following.
\begin{equation*}
\nabla_X Y = dY(X) + \Gamma(X,Y),
\end{equation*}
where $X,Y$ are vector fields on $\R^{2n}$. As above, it will be important to track the point $q\in\R^{2n}$ at which $\Gamma$ is evaluated, and we denote this $\Gamma_q(X,Y)$. Abusing notation as before, we will also write $\Gamma_v(X,Y)=\Gamma_{v(\rho)}(X,Y)$.

We are nearly ready to prove Lemma \ref{lem:DesingTechnical}, but we need just a few simple estimates, which are collected in the following result.

\begin{lemma}\label{lem:LinearAlgebra}
For all $\rho\in D_c$, the following inequalities hold:
\begin{align}
(1-C|\rho|)\|\hat{v}_s\|_{g_0}&\leq \|\hat{v}_t\|_{g_0}\leq (1+C|\rho|)\|\hat{v}_s\|_{g_0}\label{eq:DesingEst1s}\\
\|v_s\|_{\widetilde{g}}&\leq C\|\hat{v}_s\|_{\widetilde{g}}\label{eq:DesingEst1p}\\
\|v_t\|_{\widetilde{g}}&\leq C\|\hat{v}_t\|_{\widetilde{g}}\label{eq:DesingEst1q}\\
\|\hat{v}_s\|_{\widetilde{g}}\|\hat{v}_t\|_{\widetilde{g}}&\leq (1+C|\rho|)\|\hat{v}_s\wedge\hat{v}_t\|_{\widetilde{g}}.\label{eq:DesingEst1r}
\end{align}
\end{lemma}

\begin{proof}
We begin by observing that $\widetilde{g}(q)=g_0(q)+O(|q|^2)$, and $|\hat{v}(\rho)|+|v(\rho)|\leq C|\rho|$, so it is sufficient to prove the above estimates with $\widetilde{g}$ replaced with $g_0$.  Next recall that $P(\rho)=(\rho^k,\delta^{k-1}\rho,0,\ldots,0)$ and thus $\hat{v}=P + F$, with $F$ defined at the beginning of this section. Observe that
\begin{equation}\label{eq:DesingEst1y}
\|F_t\|_{g_0}+\|F_s\|_{g_0}\leq C|\rho|^k\leq C|\rho|(|\rho|^{2k-2}+\delta^{2k-2})^{\frac{1}{2}}
\end{equation}
and
\begin{equation}\label{eq:DesingEst1z}
\|P_s\|_{g_0} = (|\rho|^{2k-2}+\delta^{2k-2})^{\frac{1}{2}} = \|P_t\|_{g_0}.
\end{equation}
Consequently, for all $\rho\in \mathcal{D}_c$ we have
\begin{align}
(1-C|\rho|)(|\rho|^{2k-2}+\delta^{2k-2})^{\frac{1}{2}}&\leq \min(\|\hat{v}_s\|_{g_0},\|\hat{v}_t\|_{g_0}),\qquad\text{and}\label{eq:DesingEst2a}\\
\max(\|\hat{v}_s\|_{g_0},\|\hat{v}_t\|_{g_0})&\leq (1+C|\rho|)(|\rho|^{2k-2}+\delta^{2k-2})^{\frac{1}{2}};
\end{align}
inequalities (\ref{eq:DesingEst1s}) are immediate.  Also,
\begin{equation}\label{eq:DesingEst1b}
\|v_s\|_{g_0} + \|v_t\|_{g_0} \leq C |\rho|^{k-1} \leq C(|\rho|^{2k-2}+\delta^{2k-2})^{\frac{1}{2}},
\end{equation}
Combining this with (\ref{eq:DesingEst2a}) proves inequalities (\ref{eq:DesingEst1p}) and (\ref{eq:DesingEst1q}). To prove (\ref{eq:DesingEst1r}), we note that $\la P_s,P_t\ra_{g_0}=0$ since $P$ is a complex polynomial; then by using (\ref{eq:DesingEst1y}) - (\ref{eq:DesingEst2a}) we have
\begin{align*}
\big|\la \hat{v}_s,\hat{v}_t\ra_{g_0} \big| &\leq  \|F_s\|_{g_0}\|F_t\|_{g_0}
+ \|P_s\|_{g_0}\|F_t\|_{g_0}+ \|P_t\|_{g_0}\|F_s\|_{g_0}\\
&\leq C|\rho| \|\hat{v}_s\|_{g_0}\|\hat{v}_t\|_{g_0},
\end{align*}
for all $\rho\in \mathcal{D}_c$.  Consequently
\begin{align*}
\|\hat{v}_s\wedge\hat{v}_t\|_{g_0}^2 &= \|\hat{v}_s\|_{g_0}^2\|\hat{v}_t\|_{g_0}^2-\la \hat{v}_s,\hat{v}_t\ra_{g_0}^2\\
&\geq (1-C|\rho|^2)\|\hat{v}_s\|_{g_0}^2\|\hat{v}_t\|_{g_0}^2,
\end{align*}
and inequality (\ref{eq:DesingEst1r}) follows immediately.
\end{proof}

The following result will also be important in the proof of Lemma \ref{lem:DesingTechnical}.

\begin{lemma}\label{lem:ConformallyParameterized}
Let $(u,S,j,J)$ be an immersed $J$-curve in an almost Hermitian manifold $(M,J,g)$, let $\zeta\in S$, and let $(s,t)$ be local complex coordinates around $\zeta$ so that $u_s + J u_t = 0$.  Then
\begin{equation*}
\big(\nabla_{u_s} u_s + \nabla_{u_t} u_t\big)^\top = 0,
\end{equation*}
where $\nabla$ is the Levi-Civita connection associated to $g$, and $X\mapsto X^\top$ is the orthogonal projection to the tangent space of the image $u$.
\end{lemma}
\begin{proof}
We compute
\begin{align*}
\big(\nabla_{u_s} u_s + \nabla_{u_t} u_t\big)^\top &= \big(\nabla_{u_s} u_s + \nabla_{u_t} (Ju_s)\big)^\top = \big(\nabla_{u_s} u_s + J\nabla_{u_t} u_s\big)^\top\\
&= (\nabla_{u_s} u_s)^\top + J(\nabla_{u_s} u_t)^\top = (\nabla_{u_s} u_s)^\top + J^2(\nabla_{u_s} u_s)^\top\\
&=0
\end{align*}
Where to obtain the second equality we have employed the Leibniz rule, together with the fact that $((\nabla J) u_s)^\top = 0.$  Indeed, this result follows from the fact that the tangent planes of the image of $u$ are $J$-invariant, and $\la (\nabla J)X,X\ra = 0 = \la (\nabla J)X, JX \ra$.  This latter result is elementary, and a proof can be found in \cite{Fj09b}.  The remaining equalities are then standard.
\end{proof}

With our preparations completed, we now finish the proof of Property (D\ref{en.D5}).

\begin{proof}[Proof of Lemma \ref{lem:DesingTechnical}]
We begin by defining
\begin{align}
\mathcal{E}_1^s&:=\left(\Pi_{(\hat{v},T\hat{v})}^{\bot_{\tilde{g}}}-\Pi_{(v,T\hat{v})}^{\bot_{\tilde{g}}}
\right)\comp \Pi_{(v,T\hat{v})}^{\bot_{\tilde{g}}}
(v_{ss})\\
\mathcal{E}_2^s&:=\Pi_{(\hat{v},T\hat{v})}^{\bot_{\tilde{g}}}
\big(\Gamma_{\hat{v}}(\hat{v}_s,\hat{v}_s)\big)- \Pi_{(v,T\hat{v})}^{\bot_{\tilde{g}}}
\big(\Gamma_v(v_s,v_s)\big)\label{eq:EstE2a}.
\end{align}
Recall that $\Pi_{(\hat{v},T\hat{v})}^{\bot_{\tilde{g}}}\circ\Pi_{(v,T\hat{v})}^{\bot_{\tilde{g}}}=0$, and thus
\begin{equation}
\Pi_{(\hat{v},T\hat{v})}^{\bot_{\tilde{g}}} (v_{ss})=
\mathcal{E}_1^s+\Pi_{(v,T\hat{v})}^{\bot_{\tilde{g}}} (v_{ss}).
\end{equation}
Recall that for any $0<r<r_0/4$ and $\rho\in \mathcal{D}_{r}$ we have $\beta(|\rho|/r_0)=1$, so it follows that $\hat{v}-v$ is a linear function on $\mathcal{D}_{r_0}$, so $\hat{v}_{ss}=v_{ss}$. Consequently,
\begin{equation*}
\two_{\hat{v}}(\hat{v}_s,\hat{v}_s)=\Pi_{(v,T\hat{v})}^{\bot_{\tilde{g}}}(\nabla_{v_s}v_s) + \mathcal{E}_1^s + \mathcal{E}_2^s.
\end{equation*}
Next we estimate the $\mathcal{E}^s$ terms.  First note that $\max(|v(\rho)|,|\hat{v}(\rho)|)\leq C|\rho|$ for all $\rho\in \mathcal{D}_c$. Combining this with $\|\Gamma_q\|_{g_0}\leq C|q|$ yields
\begin{equation}\label{eq:EstE2b}
\|\Gamma_v(v_s,v_s)\|_{g_0}+\|\Gamma_{\hat{v}}(\hat{v}_s,\hat{v}_s)\|_{g_0} \leq C|\rho| (\|v_s\|_{g_0}^2+\|\hat{v}_s\|_{g_0}^2),
\end{equation}
and thus by (\ref{eq:DesingEst1p}) we have
\begin{equation}
\|\mathcal{E}_2^s\|_{g_0} \leq C |\rho| \|\hat{v}_s\|_{g_0}^2.
\end{equation}
To estimate $\mathcal{E}_1$, we first recall that  $\widetilde{g}(q)= g_0(q) +O(|q|^2)$, and thus for any plane $\mathcal{P}\subset \R^{2n}\simeq T\R^{2n}$, we have
\begin{equation*}
\big\|\Pi_{(q_1,\mathcal{P})}^{\bot_{\tilde{g}}} - \Pi_{(q_2,\mathcal{P})}^{\bot_{\tilde{g}}} \big\|_{g_0} \leq C(|q_1|^2 + |q_2|^2)
\end{equation*}
for any $q_1,q_2\in \R^{2n}$ with $\max(|q_1|,|q_2|)\leq c$.  Combining this with our estimates for $v$, $\hat{v}$, and $\Gamma$ yields the following:
\begin{align*}
\|\mathcal{E}_1^s\|_{g_0}&\leq C|\rho|^2\| \Pi_{(v,T\hat{v})}^{\bot_{\tilde{g}}} (v_{ss})\|_{g_0}\\
&\leq C|\rho|^2 \| \Pi_{(v,T\hat{v})}^{\bot_{\tilde{g}}} (\nabla_{v_s} v_s)\|_{g_0} + C |\rho|^3\|v_s\|_{g_0}^2.
\end{align*}
Combining these inequalities then yields
\begin{equation*}
\|\mathcal{E}_1^s\|_{g_0}+\|\mathcal{E}_2^s\|_{g_0}\leq C|\rho|\big(\|\Pi_{(v,T\hat{v})}^{\bot_{\tilde{g}}}(\nabla_{v_s}v_s)\|_{g_0} +  \|\hat{v}_s\|_{g_0}^2 \big).
\end{equation*}
By replacing $s$ with $t$ above, one may define $\mathcal{E}_1^t$ and $\mathcal{E}_2^t$, and prove  \begin{equation*}
\two_{\hat{v}}(\hat{v}_t,\hat{v}_t)=\Pi_{(v,T\hat{v})}^{\bot_{\tilde{g}}}(\nabla_{v_t}v_t) + \mathcal{E}_1^t + \mathcal{E}_2^t,
\end{equation*}
with similar estimates for the $\mathcal{E}^t$ terms.  We now employ Lemma \ref{lem:ConformallyParameterized}, and the fact that $H:=\tr B = \tr_{\mathcal{T}} J\nabla J$, to obtain
\begin{equation*}
\nabla_{v_s}v_s + \nabla_{v_t} v_t = J(\nabla_{v_s} J) v_s + J(\nabla_{v_t} J) v_t,
\end{equation*}
and consequently
\begin{equation*}
\Pi_{(v,T\hat{v})}^{\bot_{\tilde{g}}}(\nabla_{v_t}v_t) = - \Pi_{(v,T\hat{v})}^{\bot_{\tilde{g}}}(\nabla_{v_s}v_s) + \Pi_{(v,T\hat{v})}^{\bot_{\tilde{g}}}\big(J(\nabla_{v_s}J)v_s + J(\nabla_{v_t}J)v_t\big),
\end{equation*}
and
\begin{equation*}
\|\Pi_{(v,T\hat{v})}^{\bot_{\tilde{g}}}(\nabla_{v_t}v_t)\|_{\tilde{g}}\leq \|\Pi_{(v,T\hat{v})}^{\bot_{\tilde{g}}}(\nabla_{v_s}v_s)\|_{\tilde{g}} + C\|\hat{v}_{s}\|_{\tilde{g}}^2.
\end{equation*}
Combining this with (\ref{eq:DesingEst1s}) and our above inequalities then yields
\begin{equation}\label{eq:MathcalETerms}
\sum_{i=1}^2\big(\|\mathcal{E}_i^s\|_{g_0}+\|\mathcal{E}_i^t\|_{g_0}\big)\leq C|\rho|\big(\|\Pi_{(v,T\hat{v})}^{\bot_{\tilde{g}}}(\nabla_{v_s}v_s)\|_{g_0} +  \|v_s\|_{g_0}^2 \big).
\end{equation}
For clarity, we then define
\begin{align*}
&\mathcal{E}^s:=\mathcal{E}_1^s+\mathcal{E}_2^s,\qquad\qquad\qquad \mathcal{E}^t:=\mathcal{E}_1^t+\mathcal{E}_2^t,\qquad \qquad \mathcal{E}:= \max(\mathcal{E}^s,\mathcal{E}^t)\\
&X^s:= \Pi_{(v,T\hat{v})}^{\bot_{\tilde{g}}}(\nabla_{v_s} v_s),\qquad X^t:=\Pi_{(v,T\hat{v})}^{\bot_{\tilde{g}}}(\nabla_{v_t}v_t)\\
&V:=\Pi_{(v,T\hat{v})}^{\bot_{\tilde{g}}}\big( J(\nabla_{v_s}J)v_s+ J(\nabla_{v_t}J)v_t \big),
\end{align*}
so that $B_{\hat{v}}(\hat{v}_s,\hat{v}_s) = X^s + \mathcal{E}^s$, $B_{\hat{v}}(\hat{v}_t,\hat{v}_t) = X^t + \mathcal{E}^t$, and $V= X^s +X^t$.  Finally, we can estimate
\begin{align*}
\la B_{\hat{v}}(\hat{v}_s,\hat{v}_s), B_{\hat{v}}(\hat{v}_t,\hat{v}_t) \ra_{\tilde{g}} &\leq -\|X^s\|_{\tilde{g}}^2 + \|X^s\|_{\tilde{g}} \| V\|_{\tilde{g}} + 2\|\mathcal{E}\|_{\tilde{g}} \|X^s\|_{\tilde{g}}\\
&\qquad+ \|\mathcal{E}\|_{\tilde{g}} \|V\|_{\tilde{g}} + \|\mathcal{E}\|_{\tilde{g}}^2\\
&\leq {\textstyle \frac{1}{2}}\|V\|_{\tilde{g}}^2 + (-{\textstyle \frac{1}{2}} + C|\rho|) \|X^s\|_{\tilde{g}}^2 + C |\rho| \|\hat{v}_s\|_{\tilde{g}}^4.
\end{align*}
Recall that $(1-C|\rho|) \|\hat{v}_s\|_{\tilde{g}}^4\leq \|\hat{v}_s\wedge \hat{v}_t\|_{\tilde{g}}^2$, so that if $|\rho|$ is sufficiently small (depending only on the $g$, $J$,$F$, and $\epsilon$, but not $\delta$) we find that indeed
\begin{equation*}
\frac{\la B_{\hat{v}}(\hat{v}_s,\hat{v}_s),B_{\hat{v}}(\hat{v}_t,\hat{v}_t)\ra_{\tilde{g}}}{\|\hat{v}_s\wedge \hat{v}_t\|_{\tilde{g}}^2}\leq \epsilon + \sup_{q\in M} \sup_{\substack{e\in  T_q M\\ \|e\|_{\tilde{g}}=1}}{\textstyle\frac{1}{2}}\|J(\nabla_e J)e + J(\nabla_{Je} J) Je\|_{\tilde{g}}^2,
\end{equation*}
which is precisely the desired inequality.  This completes the proof of
Lemma \ref{lem:DesingTechnical}, and hence Property (D\ref{en.D5}) is proven as well.

\end{proof}

We now move on to proving Property (D\ref{en.D6}). We begin by showing that Gaussian curvature $K_{u^*g}:S\setminus \mathcal{Z}_u\to \R$ is integrable.  Since $K_{u^*g}$ is defined and smooth on the compliment of the set of critical points $\mathcal{Z}_u\subset S$, it is sufficient to prove that $K_{v^*\widetilde{g}}$ is integrable on $\mathcal{D}_{r}$.  To that end, recall that the Gauss equations for immersed $J$-curves guarantee that $K_{v^*\widetilde{g}}$ is uniformly bounded from above in terms of $\nabla J$ and the sectional curvature of $M$. Since $\Area_{v^*{\widetilde{g}}}(\mathcal{D}_{r})<\infty$ (more precisely, $\mathcal{D}_{r}$ has finite measure), it follows that a modification of the of the monotone convergence theorem guarantees that $K_{v^*\widetilde{g}}$ is integrable whenever
\begin{equation*}
\lim_{a\to 0}\int_{\mathcal{D}_r\setminus \mathcal{D}_{a}} K_{v^*\widetilde{g}}>-\infty.
\end{equation*}
To show this integral is finite, we define for each $a>0$ the parameterized paths $\alpha,\hat{\alpha}:\R/2\pi\mathbb{Z} \to \R^{2n}$ by
\begin{align*}
\alpha(\theta):&=v(a e^{i\theta})\\
&=(a^k e^{ik\theta},0,\ldots,0) + F(ae^{i\theta})\\
\hat{\alpha}(\theta):&=(a^k e^{ik\theta},0,\ldots,0).
\end{align*}
Along the image of $\alpha$ we define the vector field $\nu(\theta)\in T_{\alpha(\theta)}\R^{2n}$ to be the unique ``inward pointing'' $\widetilde{g}$-unit vector field which is tangent to the image of $v$ and $\widetilde{g}$-orthogonal to $v(\partial \mathcal{D}_a)$.  We similarly define the $g_0$-unit vector field $\hat{\nu}$ along the image of $\hat{v}$, which can be explicitly written as
\begin{equation*}
\hat{\nu}(\theta)=-(e^{ik\theta},0,\ldots,0).
\end{equation*}
Using $\widetilde{g}(q)-g_0(q) = O(|q|^2)$ and $\|\Gamma_q\|_{\widetilde{g}}\leq C|q|$ it is straight forward to show
\begin{equation*}
\Big|\int_0^{2\pi}\frac{\la \nabla_{\alpha'}\alpha',\nu\ra_{\widetilde{g}}}{\|\alpha'\|_{\widetilde{g}}} d\theta - \int_0^{2\pi}\frac{\la \hat{\alpha}'',\hat{\nu}\ra_{g_0}}{\|\hat{\alpha}'\|_{g_0}} d\theta\Big| \leq Ca.
\end{equation*}
Furthermore, letting $\kappa_{v^*\widetilde{g}}$ denote the geodesic curvature of $\partial \mathcal{D}_a$, and applying the Gauss-Bonnet theorem, we find
\begin{align*}
\Big|2\pi k - \Big(\int_{\mathcal{D}_r\setminus \mathcal{D}_a}K_{v^*\widetilde{g}} \Big) -\Big(\int_{\partial \mathcal{D}_r}\kappa_{v^*\widetilde{g}}\Big) \Big|&=\Big| 2\pi k - \int_{\partial \mathcal{D}_a}\kappa_{v^*\widetilde{g}} \Big|\\
&=\Big| 2\pi k - \int_0^{2\pi}\frac{\la \nabla_{\alpha'}\alpha',\nu\ra_{\widetilde{g}}}{\|\alpha'\|_{\widetilde{g}}} d\theta \Big|\\
&\leq  Ca + \Big| 2\pi k - \int_0^{2\pi}\frac{\la \hat{\alpha}'',\hat{\nu}\ra_{g_0}}{\|\hat{\alpha}'\|_{g_0}} d\theta \Big|\\
&= Ca,
\end{align*}
which tends to zero as $a\to 0$, and hence the $K_{u^*g}$ is integrable on $S$.  Moving on with the proof of Property (D\ref{en.D6}), we again let $\mathcal{Z}_u\subset S$ denote the set of critical points of $u$, and $\mathcal{B}_\epsilon(z):=\{\zeta\in S:\dist_{u^*g}(z,\zeta)<\epsilon\}$. Then we break the open set $\mathcal{U}\subset S$ into three regions:
\begin{align*}
\mathcal{U}_1 &= \bigcup_{z\in \mathcal{Z}_u\cap{\mathcal{U}^{\epsilon}}} \mathcal{B}_\epsilon(z)\\
\mathcal{U}_2 &= \mathcal{U}\cap\left(\bigcup_{z\in Z}\mathcal{B}_\epsilon(z)\setminus\mathcal{U}_1\right)\\
\mathcal{U}_3 &= \mathcal{U}\setminus(\mathcal{U}_1\cup\mathcal{U}_2).
\end{align*}
Note that by construction, for all $\zeta\in \mathcal{U}_3$ we have $u(\zeta)=\hat{u}(\zeta)$, and thus
\begin{equation}\label{eq:CurvatureThresh1}
-\int_{\mathcal{U}_3} K_{u^*g} = -\int_{\mathcal{U}_3} K_{\hat{u}^*g}.
\end{equation}
Next, we fix $\epsilon'>0$ and note that since $K_{u^*g}$ is integrable and Property (D\ref{en.D1}) holds, it follows that without loss of generality we may assume $\epsilon>0$ is sufficiently small so that
\begin{equation}\label{eq:CurvatureThresh2}
\max\Big(-\int_{\mathcal{U}_1}K_{u^*g},-\int_{\mathcal{U}_2}K_{u^*g}\Big)\leq \sum_{z\in \mathcal{Z}_u}\int_{\mathcal{B}_\epsilon(z)}|K_{u^*g}| \leq \epsilon'/3.
\end{equation}
Since $\Area_{u^*g}(S)$ is finite, it follows from Property  (D\ref{en.D4}) that without loss of generality we may assume that $\epsilon>0$ is sufficiently small so that the following holds:
\begin{equation}\label{eq:CurvatureThresh3}
(\epsilon+C_{geom})\Area_{\hat{u}^*g}\big(\cup_{z\in \mathcal{Z}_u}
\mathcal{B}_{\epsilon}(z)\big)\leq \epsilon'/3,
\end{equation}
where
\begin{equation*}
C_{geom}:=\sup_{q\in M}|K_{sec}(q)| + \sup_{\substack{q\in M\\ J\mathcal{P}_q = \mathcal{P}_q}}\|\tr_{\mathcal{P}_q}J\nabla J\|_{g}^2.
\end{equation*}
Consequently,
\begin{align}
\int_{\mathcal{U}_2}K_{\hat{u}^*g}&\leq (\epsilon+C_{geom})\Area_{\hat{u}^*g}\big(\cup_{z\in \mathcal{Z}_u}
\mathcal{B}_{\epsilon}(z)\big)\notag\\
&\leq
\epsilon'/3\notag\\
&\leq 2\epsilon'/3 + \int_{\mathcal{U}_2}K_{u^*g}.\label{eq:CurvatureThresh4}
\end{align}
Lastly, recall that the definition of $\hat{u}$ guarantees that for every $\zeta\in \partial \mathcal{B}_\epsilon(z)$ with $z\in \mathcal{Z}_u$ we have $\hat{u}(\zeta)= u(\zeta)$.  Thus we compute
\begin{align}
-\int_{\mathcal{B}_\epsilon(z)}K_{\hat{u}^*g}&=-2\pi + \int_{\partial
\mathcal{B}_\epsilon(z)}\kappa_{\hat{u}^*g}\notag\\
&=-2\pi + \int_{\partial
\mathcal{B}_\epsilon(z)}\kappa_{u^*g}\notag\\
&\to 2\pi(k-1)\label{eq:CurvatureThresh5}
\end{align}
as $\epsilon\to 0$.  Recall that $k-1=\order(z)$, so by integrating over $\mathcal{U}_1$ we find that for sufficiently small $\epsilon>0$ we have
\begin{equation}\label{eq:CurvatureThresh6}
(1-\epsilon')2\pi\sum_{z\in \mathcal{Z}_u\cap \mathcal{U}^\epsilon}\order(z)\leq -\int_{\mathcal{U}_1}K_{\hat{u}^*g}.
\end{equation}
Combining (\ref{eq:CurvatureThresh1}), (\ref{eq:CurvatureThresh2}) (\ref{eq:CurvatureThresh4}), and (\ref{eq:CurvatureThresh6}) yields the desired estimate:
\begin{equation*}
(1-\epsilon')2\pi\sum_{z\in \mathcal{Z}_u\cap\mathcal{U}^\epsilon} -\int_{\mathcal{U}} K_{u^*g} \leq \epsilon'-\int_{\mathcal{U}}K_{\hat{u}^*g}
\end{equation*}
This completes the
proof of Property (D\ref{en.D6}) as well as Proposition \ref{prop:Desingularization}.
\end{proof}

\subsection{Proof of Proposition
\ref{prop:nodesAndCurvature}}\label{sec:nodesAndCurvature}

\begin{proof}
The proof of Proposition \ref{prop:nodesAndCurvature} consists of two main parts.  The first part consists of passing to the desired subsequence and constructing the finite set $\mathcal{S}$ and showing the first part of the proposition holds.  The basic argument here is to construct the subsequence and $\mathcal{S}$ by iteratively passing to further and further subsequences with $\#\mathcal{S}$ getting larger in each subsequent iteration.  We then argue that if $\# \mathcal{S}$ is arbitrarily large, then by monotonicity the $S_k$ must have arbitrarily large area, which is a contradiction.  The second part of the proof consists of a covering argument which reduces the problem to showing that the integral of the Gaussian curvature on disks cannot be arbitrarily negative;  we prove the reduced problem by recalling a differential equation which relates the area of an intrinsic disk to the integral of the Gaussian curvature on said disk, and conclude that since the area is a priori bounded, so too is the desired integral.

Moving on to the actual proof, we note that since the $u_k$ are robustly $\mathcal{K}$-proper, there exists a compact set $\widetilde{\mathcal{K}}_i$ for $i=1,2$ such that
\begin{equation*}
\mathcal{K}:=\widetilde{\mathcal{K}}_0\subset \subset \widetilde{\mathcal{K}}_1\subset \subset \widetilde{\mathcal{K}}_2\subset\subset M,
\end{equation*}
and for which the $u_k$ are robustly $\widetilde{\mathcal{K}}_2$-proper.  Next we fix $\delta_0>0$ such that the following conditions hold.
\begin{enumerate}
\item $\mathcal{O}_{10\delta_0}^{g_k}(\widetilde{\mathcal{K}}_i)\subset \widetilde{\mathcal{K}}_{i+1}$ for all $k\in \mathbb{N}$ and $i\in \{0,1\}$.
\item $10^{10}\delta_0 < \min\big(C_{\infty}^{-1/2},\inf_{q\in \widetilde{\mathcal{K}}_2} \inj_M^{g}(q)\big)$, where
    \begin{equation*}
       C_\infty:=2+\sup_{q\in M} |K_{sec}^g(q)| + \sup_{q\in M} \|\nabla J\|_{g}^2.
    \end{equation*}
\item for each $p\in \widetilde{\mathcal{K}}_2$ and $g_k$-geodesic polar coordinates $(x^1,\ldots,x^{2n})$ centered at $p$ for which $J_k(p) \partial_{x^i} = \partial_{x^{i+n}}$ for $i=1,\ldots,n$ define the differential forms $\check{\omega}_p$ and $\check{\lambda}_p$ on $\mathcal{O}_p:=\mathcal{O}_{10\delta_0}^{g_k}(p)$ by
    \begin{equation*}
    \check{\omega}_p:={\textstyle \sum_{i=1}^{n}}dx^i\wedge dx^{i+n}\qquad\text{and}\qquad\check{\lambda}_p:= {\textstyle \sum_{i=1}^{n}}x^i dx^{i+n}.
    \end{equation*}
    We then require that $\delta_0>0$ is sufficiently small so that
    \begin{equation*}
    \sup_{\substack{p\in \widetilde{\mathcal{K}}_2\\ k\in \mathbb{N}}} \|\check{\lambda}_p\|_{L^\infty(\mathcal{O}_p)} \leq 1, \qquad \sup_{\substack{p\in \widetilde{\mathcal{K}}_2\\ k\in \mathbb{N}}}\|\check{\omega}_p\|_{L^\infty(\mathcal{O}_p)} \leq 2.
    \end{equation*}
    and
    \begin{equation*}
    \inf_{\substack{q\in \mathcal{O}_p\\X\in T_q M\setminus \{0\}}} \big| 1- \check{\omega}_p (X,JX)/\|X\|^2\big| \leq 10^{-10}.
    \end{equation*}
\end{enumerate}
We now aim to prove the following:
\begin{itemize}
\item[$(*)$] \emph{After passing to a subsequence (still denoted with subscripts $k$), there exists a finite set $\mathcal{S}\subset M$ with the property that for each $\delta\in(0,\delta_0)$ there exists an $\epsilon>0$ and $k_0\in\mathbb{N}$ such that if $k\geq k_0$ and $v_k(\zeta)\in \widetilde{\mathcal{K}}_1\setminus \mathcal{O}_\delta^{g_k}(\mathcal{S})$ then $\inj_{S_k}^{v_k^*g_k}(\zeta)>\epsilon$.}
\end{itemize}
Since by construction $\mathcal{O}_{\delta_0}^{g_k}(\mathcal{K})\subset \widetilde{\mathcal{K}}_1$, we see that if $(*)$ holds, then the first part of Proposition \ref{prop:nodesAndCurvature} is true. To find the set $\mathcal{S}$ and the desired subsequence one can argue iteratively in the following way.  Define $\mathcal{S}_0:=\emptyset$, and find a sequence $\zeta_{1,k}\in S_k$ which has the property that a subsequence (denoted with subscripts $k_1$) satisfies
\begin{equation*}
v_{k_1}(\zeta_{1,k_1})\in \widetilde{\mathcal{K}}_1,\quad\inj_{S_{k_1}}^{v_{k_1}^*g_{k_1}}(\zeta_{1,k_1})\to 0,\quad\text{and}\quad\lim_{k_1\to\infty} v_{k_1}(\zeta_{1,k_1})=:\sigma_1\notin\mathcal{S}_0.
\end{equation*}
Define $\mathcal{S}_1:=\mathcal{S}_0\cup\{\sigma_1\}$, and pass to a further subsequence (denoted with subscripts $k_2$) and find a sequence $\zeta_{2,k_2}\in S_{k_2}$ which has the property that
\begin{equation*}
v_{k_2}(\zeta_{2,k_2})\in \widetilde{\mathcal{K}}_1,\quad\inj_{S_{k_2}}^{v_{k_2}^*g_{k_2}}(\zeta_{2,k_2})\to 0\quad \text{and}\quad\lim_{k_2\to\infty} v_{k_2}(\zeta_{2,k_2})=:\sigma_2 \notin \mathcal{S}_1.
\end{equation*}
Define $\mathcal{S}_2:=\mathcal{S}_1\cup\{\sigma_2\}$, and iterate.  Of course if the procedure terminates after a finite number of iterations, then $(*)$ is true; otherwise one can construct a singular set $\mathcal{S}_{n_0}$ with $n_0$ distinct points in $\widetilde{\mathcal{K}}_1$ for $n_0$ arbitrarily large. To derive a contradiction, we assume the latter, in which case we conclude that there exists a point $p\in \widetilde{\mathcal{K}}_1$, a $k\in\mathbb{N}$, a $\delta'$ satisfying $0<\delta'<\min(1,\delta_0)/10^{10}$, and points $\{\zeta_{1},\ldots,\zeta_{n_0}\}\subset S_k$ for which the following hold
\begin{enumerate}
\item $n_0> C_G + 8C_A/(\pi \delta_0^2)$
\item $v_k(\zeta_{i})\in \mathcal{O}_{\delta_0}^{g_k}(p)$ for $i=1,\ldots,n_0$,
\item $\min_{i\neq j} \dist_{g_k}\big(v_k(\zeta_{i}),v_k(\zeta_{j})\big)\geq \delta'$,
\item $10^{10} \epsilon_k < \pi \delta'^2 \leq \min (1,\delta')$ where $v_k$ is the immersed approximation associated to the pair $(u_k,\epsilon_k)$.
\item $\inj_{S_k}^{v_k^*g_k}(\zeta_{i,k})< \frac{1}{4}\min\big(\frac{\delta'^2}{n_0 10^{10}},C_{\infty}^{-1},\delta_0\big)=:\epsilon'$.
\item $\sup_{q \in M} |K_{sec}^{g_k}(q)| + \sup_{q\in M} \|\nabla J_k\|_{g_k}^2 +1\leq C_\infty$
\end{enumerate}
Since $\delta_0$ can be chosen from an open set, we may assume without loss of generality that  the ``trimmed'' $J_k$-curves given by
\begin{equation*}
\widetilde{\mathbf{u}}_k:=(u_k,\widetilde{S}_k,j_k,J_k)\qquad\text{with}\qquad \widetilde{S}_k:=u_k^{-1}\big(\overline{\mathcal{O}_{4\delta_0}^{g_k}(p)}\big)
\end{equation*}
are compact curves with smooth boundary, and $\partial \widetilde{S}_k= u_k^{-1}\big( \partial \mathcal{O}_{4\delta_0}^{g_k}(p)\big)$.  Such a choice is possible since the $\mathbf{u}_k$ are robustly $\widetilde{\mathcal{K}}_2$-proper and $\mathcal{O}_{10\delta_0}^{g_k}(\widetilde{\mathcal{K}}_1)\subset \widetilde{\mathcal{K}}_2$. Recall that by property (D\ref{en.D5}) of Proposition \ref{prop:Desingularization}, it follows that the Gaussian curvature $K_{v_k^*g_k}(\zeta)$ is uniformly bounded from above by $C_\infty$ for all $\zeta\in S_k$; since it's the case that $\inj_{S_k}^{v_k^*g_k}(\zeta_{i})\leq C_{\infty}^{-1}\leq C_\infty^{-1/2}$, it follows that for each $i=1,\ldots,n_0$, there exists a simple $v_k^*g_k$-unit speed geodesic $\alpha_i:[0,\ell_i]\to S_k$ such that $\zeta_{i}=\alpha_i(0)=\alpha_i(\ell_i)$ and $\Length_{v_k^*g_k}(\alpha_i)=\ell_i \leq 2\epsilon'$, and the $\alpha_i$ are pair-wise disjoint. Note that in general $\alpha_i'(0)\neq \alpha_i'(\ell_i)$.  Observe that $\widetilde{S}_k$ need not be connected but it has finitely many connected components; furthermore if we consider the (non-compact) manifold $\widetilde{S}_k\setminus \alpha_1$, and recall Definition \ref{def:genus}, we see that one of the two scenarios must occur:
\begin{enumerate}
\item $\Genus (\widetilde{S}_k) > \Genus (\widetilde{S}_k\setminus \alpha_1)$
\item the number of connected components of $\widetilde{S}_k\setminus \alpha_1$ is strictly larger than the number of connected components of $\widetilde{S}_k$.
\end{enumerate}
Indeed, in general removing a simple loop from a surface either decreases its genus, or increases the number of connected components.  Since the genus of the $S_k$ (and hence $\widetilde{S}_k$) are bounded by $C_G$, and the $\alpha_i$ are pair-wise disjoint and simple, it follows that $\widetilde{S}_k\setminus \cup_{i=1}^{n_0}\alpha_i$ has at least $n_1:=n_0- C_G+1$ connected components of zero genus which have non-trivial intersection with $v_k^{-1}\big(\mathcal{O}_{2\delta_0}^{v_k^*g_k}(p)\big)$; we label these connected components $\check{S}^i$ for $i=1,\ldots,n_1$.  We now make the following claim.
\begin{lemma}\label{lem:componentsRunToBoundary}
For each $i=1,\ldots,n_1$, the connected component $\check{S}^i$ has non-trivial intersection with $v_k^{-1}\big(\partial \mathcal{O}_{4\delta_0}^{v_k^*g_k}(p)\big)$.
\end{lemma}
Before proceeding with the proof of Lemma \ref{lem:componentsRunToBoundary}, we use it to finish the proof of of the first part of Proposition \ref{prop:nodesAndCurvature}.  Indeed, as a consequence of the above lemma, it follows that
for each $i=1,\ldots, n_1$ we have $v_k(\partial \check{S}^i)\subset  \mathcal{O}_{2\delta_0}^{g_k}(p) \cup \partial\mathcal{O}_{4\delta_0}^{g_k}(p)$, and there exists a $\zeta_i'\in \check{S}^i$ with the property that $\dist_{g_k}\big(p,v_k(\zeta_i')\big)=3\delta_0$. We conclude from the monotonicity of area (Proposition \ref{prop:monotonicityOfArea}), that $\pi \delta_0^{2}/8\leq \Area_{u_k^*g_k}(\check{S}^i)$.
It then follows that
\begin{equation*}
(n_0-C_G+1) \pi \delta_0^2/8 \leq {\textstyle\sum_{i=1}^{n_1}} \Area_{u_k^*g_k}(\check{S}^i) \leq \Area_{u_k^*g_k}(S_k)\leq C_A,
\end{equation*}
which is the desired contradiction.  The proof of the first part of Proposition \ref{prop:nodesAndCurvature} will be complete once we prove Lemma \ref{lem:componentsRunToBoundary}.  To that end, we will make use of the following lemma.

\begin{lemma}\label{lem:AreaSympCompare}
Let $(M,J_k,g_k)$, $\mathcal{O}_{10\delta_0}^{g_k}(p)$, $(u_k,S_k,j_k,J_k)$, $C_A$, $v_k$ $\epsilon_k$, $\check{\omega}_p$, and $\check{\lambda}_p$ be as above.  Furthermore, let $\mathcal{U}\subset S_k$ be an open set for which $v_k(\mathcal{U})\subset \mathcal{O}_{10\delta_0}^{g_k}(p)$.  Then
\begin{equation*}
\big|\Area_{v_k^*g_k}(\mathcal{U}) - \int_{\mathcal{U}} v_k^*\check{\omega}_p \big| \leq \frac{1}{2} \Area_{v_k^*g_k}(\mathcal{U}).
\end{equation*}
\end{lemma}

\begin{proof}
Let $E$ be a unit vector tangent to the image of $v$. Define another tangent vector $F:=(JE)^\top/\|(JE)^\top\|$ which is orthonormal to $E$; here $X\mapsto X^\top$ is the orthogonal projection to the plane tangent to the image of $v$.
Recall that $J_k$ is a $g_k$-isometry, and $\|(JE)^\bot\|\leq \epsilon_k$ by property (D\ref{en.D3}) of Proposition \ref{prop:Desingularization}, so it is elementary to show that $1-\epsilon_k\leq \|(JE)^\top\|$ and $\|JE - F\|\leq \frac{2\epsilon_k}{1-\epsilon_k}$. Employing our above estimates for $\check{\omega}$, we then find
\begin{equation*}
\big| 1 - \check{\omega}(E,F) \big| \leq \big| 1 -\check{\omega}(E,JE)\big| + \big|\check{\omega}(E,JE-F)\big|\leq\frac{1}{2}.
\end{equation*}

The desired result then follows immediately by integrating.
\end{proof}

\begin{proof}[Proof of Lemma \ref{lem:componentsRunToBoundary}]
Suppose not. We let $\cl (\check{S}^i)$ denote the metric compactification of $(\check{S}^i,v_k^*g_k)$.  For example, if $\widetilde{S}_k$ is a torus, and $\check{S}^i:=\widetilde{S}_k\setminus \alpha_1 $ is an open cylinder, then $\cl(\check{S}^i)$ is a compact cylinder -- \emph{not a torus} -- with piece-wise smooth boundary. Note that each boundary component of $\cl(\check{S}^i)$ is a copy of the piece-wise smooth geodesic $\alpha_l$ for some $l\in\{1,\ldots, n_1\}$.  There are several cases to consider.

\emph{Case I: $\partial \cl(\check{S}^i)=\emptyset$.} In this case $\check{S}^i$ is closed, and $v_k(\check{S}^i)\subset \mathcal{O}_{4\delta_0}^{g_k}(p)$ on which $\check{\omega}_p=d\check{\lambda}_p$.  As a consequence of Lemma  \ref{lem:AreaSympCompare} we see that $\Area_{v_k^*g_k}(\check{S}^i) = 0$, and hence $v_k:\check{S}^i\to M$ is a constant map.  This is not possible since the $v_k$ are immersions.  Contradiction.

\emph{Case II: $\partial \cl(\check{S}^i)$ has exactly one component.}  By assumption $\check{S}^i$ has zero genus, is connected, and has empty intersection with $v_k^{-1}\big(\partial\mathcal{O}_{4\delta_0}^{g_k}(p)\big)$.  Consequently $\check{S}^i$ is a disk with with boundary component $\alpha_l$.  Let $\theta_0\in[-\pi,\pi]$ denote the exterior angle between $\alpha_l'(0)$ and $\alpha_l'(\ell_l)$.  But then we compute the following.
\begin{align*}
\pi \leq 2\pi - \theta_0 &= \int_{\check{S}^i}K_{v_k^*g_k} \leq C_\infty \Area_{v_k^*g_k}(\check{S}^i)\\
&\leq 2 C_\infty \int_{\check{S}^i}v_k^*\omega = 2C_\infty \int_{\partial\check{S}^i} v_k^* \lambda\\
& \leq 2C_\infty\ell_i \leq 4 C_\infty \epsilon_0 \leq 1
\end{align*}
which is a contradiction.

\emph{Case III: $\partial \cl(\check{S}^i)$ has exactly two components, each of which is a copy of the same $\alpha_l$.} In this case $\cl(\check{S}^i)$ is a compact cylinder, and $\check{S}^i\cup\alpha_l$ is a torus.  As in Case I, it follows that $\mathbf{u}_k$ is not generally immersed, which is a contradiction.

\emph{Case IV: $\partial \cl(\check{S}^i)$ has at least two components $\alpha_l$ and $\alpha_j$ with $l\neq j$.}  In this case we note that there exists a $\zeta\in \check{S}^i$ such that
\begin{equation*}
\min_l \dist_{g_k}\big(u_k(\zeta),u_k(\alpha_l)\big)\geq \delta'/10.
\end{equation*}
By the monotonicity of area, Proposition \ref{prop:monotonicityOfArea}, and property (D\ref{en.D4}) of Proposition \ref{prop:Desingularization}, it follows that
\begin{equation}\label{eq:random}
\Area_{v_k^*g_k}(\check{S}^i) \geq \pi \delta'^2/(2\cdot 10^2)-\epsilon_k\geq  \delta'^2/10^3.
\end{equation}
However, we note that $\cl(\check{S}^i)$ can have at most $2n_0$ geodesic boundary components.  Thus we compute
\begin{align*}
\Area_{v_k^*g_k}(\check{S}^i)\leq 2 \int_{\check{S}^i}v_k^*\omega= 2\int_{\partial \check{S}^i}v_k^*\lambda \leq 8 n_0 \epsilon' <  \delta'^2/10^8,
\end{align*}
but this contradicts (\ref{eq:random}).

Thus we see that all possible cases lead to contradictions, and thus we have completed the proof of Lemma \ref{lem:componentsRunToBoundary}.
\end{proof}

We have completed the proof of the first part of Proposition \ref{prop:nodesAndCurvature} -- indeed, we have proved more, namely the statement $(*)$. We now turn our attention toward proving the second part of Proposition \ref{prop:nodesAndCurvature}, namely we will show that for each $0<\delta<\delta_0/2$ there exists a $C>0$  such that for all sufficiently large $k$ in the subsequence, we have
\begin{equation*}
-\int_{\widehat{S}_k^{{2\delta}}}K_{v_k^*g_k}\leq C,\qquad\text{where}\qquad \widehat{S}_k^{{\delta}}:=v_k^{-1}\big(\Int({\mathcal{K}})\setminus \overline{\mathcal{O}_\delta^{g_k}(\mathcal{S})}\big);
\end{equation*}
To that end, we begin by fixing $\delta\in (0,\delta_0/2)$, and let $\mathcal{S}\subset M$,  $k_0\in\mathbb{N}$, and $\epsilon>0$ be the set and quantities guaranteed by $(*)$; furthermore we henceforth assume that we have passed to an appropriate subsequence. Observe that as a consequence of properties (D\ref{en.D4}) and (D\ref{en.D5}) of Proposition \ref{prop:Desingularization}, it is sufficient to show that for each there exists a $C_K>0$, $n_0'>0$ , $\delta'>0$, and open sets $\mathcal{O}_{\delta'}^{v_k^*g_k}(\zeta_{i,k})\subset S_k$ for $i=1,\ldots, n_k$ (all depending on $\delta$) such that
\begin{align}
\widehat{S}_k^{2\delta} &\subset \bigcup_{i=1}^{n_k}\mathcal{O}_{3\delta'}^{v_k^*g_k}(\zeta_{i,k}),\label{eq:Kover5}\\
n_k&\leq n_0',\label{eq:Kover6}\\
-\int_{\mathcal{O}_{3\delta'}^{v_k^*g_k}(\zeta_{i,k})}K_{v_k^*g_k}&\leq C_K\label{eq:Kover7}
\end{align}
To that end, we fix $\delta'$ such that $0<6\delta' < \min\big(\epsilon, 2\delta, C_{\infty}^{-1/2}\big).$
For each $k\geq k_0$ choose $\zeta_{i,k}\in  \widehat{S}_k^{2{\delta}}$, for $i=1,\ldots,n_k$, so that the sets $\mathcal{O}_{\delta'}^{v_k^*g_k}(\zeta_{i,k})$ are a maximal collection of disjoint sets.  In other words, we choose the $\zeta_{i,k}$ so that the sets $\mathcal{O}_{\delta'}^{v_k^*g_k}(\zeta_{i,k})$ are pairwise disjoint, and so that if $\zeta\in \widehat{S}_k^{2{\delta}}$, then $\mathcal{O}_{\delta'}^{v_k^*g_k}(\zeta_{i,k}) \cap \mathcal{O}_{\delta'}^{v_k^*g_k}(\zeta)\neq \emptyset$ for some $i\in \{1,\ldots,n_k\}$. Note that since $\zeta_{i,k}\in \widehat{S}_k^{2\delta}$ and since $6\delta'< \epsilon$ it follows from $(*)$ that $\mathcal{O}_r^{v_k^*g_k}(\zeta_{i,k})$ is a disk for all $r\in (0,6\delta']$.

We are now ready to show that (\ref{eq:Kover5}) holds. Indeed, suppose not; then there exists $\zeta\in \widehat{S}_k^{2\delta}$ such that $\zeta \notin \cup_{i=1}^{n_k}\mathcal{O}_{3\delta'}^{v_k^*g_k}(\zeta_{i,k})$.  By the triangle inequality, it follows that $\mathcal{O}_{\delta'}^{v_k^*g_k}(\zeta)\cap \mathcal{O}_{\delta'}^{v_k^*g_k}(\zeta_{i,k})=\emptyset$ for $i=1,\ldots,n_k$; however this contradicts the maximality of the $\mathcal{O}_{\delta'}^{v_k^*g_k}(\zeta_{i,k})$.  This proves (\ref{eq:Kover5}).

We now show (\ref{eq:Kover6}) holds.  Recall that since it's the case that the  $\mathcal{O}_{\delta'}^{v_k^*g_k}(\zeta_{i,k})\subset S_k$ are pair-wise disjoint for $i=1,\ldots,n_k$, and since it's the case that  $\Area_{v_k^*g_k}(S_k)\leq C_A$, it is sufficient to show that there exists a constant $c>0$, which is independent of $k$, such that $\Area_{v_k^*g_k}\big(\mathcal{O}_{\delta'}^{v_k^*g_k}(\zeta_{i,k})\big) \geq c$. To that end, we define the functions
\begin{equation*}
\mathcal{A}(r):=\Area_{v_k^*g_k}\big(\mathcal{O}_r^{v_k^*g_k}(\zeta_{i,k})\big)\quad\text{and}\quad \mathcal{L}(r):=\Length_{v_k^*g_k}\big(\partial \mathcal{O}_r^{v_k^*g_k}(\zeta_{i,k})\big).
\end{equation*}
By the co-area formula we recall that $\frac{d}{d r}\mathcal{A}(r)=\mathcal{L}(r)$.  Furthermore the variation of volume
formula and Gauss-Bonnet theorem yield
\begin{equation}\label{eq:dLengthIsGeoCurvature}
\frac{d}{d r}\mathcal{L}(r)=\int_{\partial
\mathcal{O}_{r}^{v_k^*g_k}(\zeta_{i,k})}\kappa_{v_k^*g_k}=2\pi
-\int_{\mathcal{O}_{r}^{v_k^*g_k}(\zeta_{i,k})}K_{v_k^*g_k}.
\end{equation}
Thus
$\mathcal{A}$ satisfies the following differential inequality.
\begin{equation}\label{eq:AreaDifferentialInequality}
\mathcal{A}(0)=0,\qquad\mathcal{A}'(0)=0,\qquad\mathcal{A}''(r)\geq 2\pi - C_\infty \mathcal{A}(r),
\end{equation}
with $C_\infty$ defined as above.
We now claim the for all $r\in [0,6\delta']$, we must have $\mathcal{A}(r)\geq \pi r^2/2$.  Indeed, if this were not the case, then for some $r_0\in (0,6\delta']$ we would have $\mathcal{A}(r)\leq \pi r_0^2/2$ for all $r\in [0,r_0]$. But then the inequality in (\ref{eq:AreaDifferentialInequality}) yields
\begin{equation*}
\mathcal{A}''(r)\geq 2\pi - C_\infty \pi r_0^2/2 \geq 3\pi/2,
\end{equation*}
where we have used the fact that $r_0\leq 6\delta'< C_\infty^{-1/2}$.  Integrating up then yields $\mathcal{A}(r_0)\geq 3\pi r_0^2/4$ which contradicts our assumption that $\mathcal{A}(r_0) \leq \pi r_0^2/2$.  We conclude that $\mathcal{A}(r)\geq \pi r^2/2$ for all $r\in [0,6\delta']$. Moreover,  $\Area_{v_k^*g_k}\big(\mathcal{O}_{\delta'}^{v_k^*g_k}(\zeta_{i,k})\big)\geq \pi \delta'^2/2$, and by our previous discussion, we see that (\ref{eq:Kover6}) holds.

All that remains to complete the proof of Proposition \ref{prop:nodesAndCurvature}, is to show that (\ref{eq:Kover7}) holds.  To that end, we suppose not.  Or in other words, for every $C_K\geq 0$, there exists $i$ and $k$ such that
\begin{equation}\label{eq:Kover11}
-\int_{\mathcal{O}_{3\delta'}^{v_k^*g_k}(\zeta_{i,k})}K_{v_k^*g_k} \geq C_K.
\end{equation}
Our above discussion shows that $\mathcal{A}$ satisfies the following integral equation and subsequent inequality
\begin{align}
\mathcal{A}(r)&=\pi r^2 - \int_0^r \int_0^t \Big( \int_{\mathcal{O}_s^{v_k^*g_k}(\zeta_{i,k})}K_{v_k^*g_k}\Big) ds dt\notag\\
&\geq \pi r^2 -{\textstyle\frac{1}{2}} C_\infty C_A r^2 + \int_0^r \int_0^t \Big(\int_{\mathcal{O}_s^{v_k^*g_k}(\zeta_{i,k})}C_\infty -K_{v_k^*g_k}\Big) ds dt.\label{eq:Kover10}
\end{align}
Since we have the point-wise bound $C_\infty\geq K_{v_k^*g_k}$, the triple integral in
(\ref{eq:Kover10}) is a monotone increasing function in $r$. Consequently
\begin{equation*}
-\int_{\mathcal{O}_{3\delta'}^{v_k^*g_k}(\zeta_{i,k})}K_{v_k^*g_k} \geq C_K
\;\;\Rightarrow \;\;
\int_{\mathcal{O}_{r}^{v_k^*g_k}(\zeta_{i,k})}(C_\infty-K_{v_k^*g_k})\geq\left\{
                                        \begin{array}{ll}
                                          0 & \text{if }r<3\delta', \\
                                          C_K & \text{if }r\geq
3\delta'
                                        \end{array}
                                      \right.
\end{equation*}
and thus
\begin{equation*}
\mathcal{A}(r)\geq \pi r^2-\frac{1}{2}C_\infty C_Ar^2+ \left\{
                                                     \begin{array}{ll}
                                                       0 & \text{if } r<3\delta'\\
                                                       \frac{1}{2}C_K(r-3\delta')^2 & \text{if
}r\geq 3\delta'
                                                     \end{array}
                                                   \right.
\end{equation*}
Evaluating the above inequality at $6\delta'$ shows that if (\ref{eq:Kover11}) holds for arbitrarily large $C_K>0$, then $\mathcal{A}(6\delta')=\Area_{v_k^*g_k}\big(\mathcal{O}_{6\delta'}^{v_k^*g_k}(\zeta_{i,k})\big)$ is also arbitrarily large.  This contradicts our assumption that the areas of the $S_k$ are uniformly bounded.  This contradiction shows that (\ref{eq:Kover7}) must hold, which in turn completes the proof of Proposition \ref{prop:nodesAndCurvature}.

\end{proof}

\bibliographystyle{amsplain}
\bibliography{TLGC}

\providecommand{\bysame}{\leavevmode\hbox to3em{\hrulefill}\thinspace}
\providecommand{\MR}{\relax\ifhmode\unskip\space\fi MR }
\providecommand{\MRhref}[2]{%
  \href{http://www.ams.org/mathscinet-getitem?mr=#1}{#2}
}
\providecommand{\href}[2]{#2}
\begin{thebibliography}{10}

\bibitem{BEHWZ03}
F.~Bourgeois, Y.~Eliashberg, H.~Hofer, K.~Wysocki, and E.~Zehnder,
  \emph{Compactness results in symplectic field theory}, Geom. Topol.
  \textbf{7} (2003), 799--888 (electronic).

\bibitem{ChSr85}
Hyeong~In Choi and Richard Schoen, \emph{The space of minimal embeddings of a
  surface into a three-dimensional manifold of positive {R}icci curvature},
  Invent. Math. \textbf{81} (1985), no.~3, 387--394.

\bibitem{Fj09b}
Joel~W. Fish, \emph{Estimates for {$J$}-curves as submanifolds},
  arXiv:0912.4445.

\bibitem{Fj07}
\bysame, \emph{Compactness results for pseudo-holomorphic curves}, {PhD}
  dissertation, New York University, 2007.

\bibitem{Gm85}
M.~Gromov, \emph{Pseudoholomorphic curves in symplectic manifolds}, Invent.
  Math. \textbf{82} (1985), no.~2, 307--347.

\bibitem{Hc97}
Christoph Hummel, \emph{Gromov's compactness theorem for pseudo-holomorphic
  curves}, Progress in Mathematics, vol. 151, Birkh\"auser Verlag, Basel, 1997.

\bibitem{IsSv00}
S.~Ivashkovich and V.~Shevchishin, \emph{Gromov compactness theorem for
  {$J$}-complex curves with boundary}, Internat. Math. Res. Notices (2000),
  no.~22, 1167--1206.

\bibitem{KsNk96b}
Shoshichi Kobayashi and Katsumi Nomizu, \emph{Foundations of differential
  geometry. {V}ol. {II}}, Wiley Classics Library, John Wiley \& Sons Inc., New
  York, 1996, Reprint of the 1969 original, A Wiley-Interscience Publication.

\bibitem{Mmp94}
Marie-Paule Muller, \emph{Gromov's {S}chwarz lemma as an estimate of the
  gradient for holomorphic curves}, Holomorphic curves in symplectic geometry,
  Progr. Math., vol. 117, Birkh\"auser, Basel, 1994, pp.~217--231.

\bibitem{SmSt92}
Mika Sepp{\"a}l{\"a} and Tuomas Sorvali, \emph{Geometry of {R}iemann surfaces
  and {T}eichm\"uller spaces}, North-Holland Mathematics Studies, vol. 169,
  North-Holland Publishing Co., Amsterdam, 1992.

\bibitem{Ta92}
Anthony~J. Tromba, \emph{Teichm\"uller theory in {R}iemannian geometry},
  Lectures in Mathematics ETH Z\"urich, Birkh\"auser Verlag, Basel, 1992,
  Lecture notes prepared by Jochen Denzler.

\end{thebibliography}

\end{document}